\documentclass{llncs}
\usepackage{makeidx}  

\usepackage{amssymb}
\usepackage{amsmath}
\usepackage{epsfig}
\usepackage{subfigure}
\usepackage{mbboard}
\usepackage{calrsfs}
\usepackage{bbm}
\usepackage{datetime}
\usepackage{mathabx}
\usepackage{harpoon}

\newcommand{\olint}{{\overline{\int}}}
\newcommand{\olsum}{{\overline{\sum}}}

\begin{document}

%
%
\title{
A Unifying Framework for Some Directed Distances in Statistics
\\[-0.2cm]
}
\titlerunning{
Distances in Statistics}  
%
\author{Michel Broniatowski \inst{1} \and Wolfgang Stummer \inst{2}}
\authorrunning{Broniatowski and Stummer} 
%
\tocauthor{Broniatowski and Stummer}
\institute{
LPSM, Sorbonne Universit\'{e}, 
4 place Jussieu, 75252 Paris, 
France. \email{michel.broniatowski@sorbonne-universite.fr.
ORCID 0000-0001-6301-5531.}
\\[0.0cm]
\and
Department of Mathematics, 
University of Erlangen--N\"{u}rnberg,
Cauerstrasse $11$, 91058 Erlangen, Germany, 
as well as Affiliated Faculty Member of the School of Business and Economics,
University of Erlangen--N\"{u}rnberg,
Lange Gasse 20, 90403 N\"{u}rnberg, Germany. \ 
\email{stummer@math.fau.de.
ORCID 0000-0002-7831-4558.
Corresponding Author.
}
}

\maketitle              

\vspace{-0.5cm}

%

\begin{abstract} 
\textit{Density-based} directed distances --- particularly 
known as divergences ---  
between probability distributions
are widely used in statistics
as well as in the adjacent research fields of
information theory, artificial intelligence and machine learning.
Prominent examples are the Kullback-Leibler
information distance (relative entropy)
which e.g. is closely connected to the 
omnipresent maximum likelihood estimation method,
and Pearson's 
$\chi^2-$distance which e.g. is used for 
the celebrated chisquare goodness-of-fit test.
Another line of statistical inference is built upon
\textit{distribution-function-based} divergences    
such as e.g. the prominent
(weighted versions of) Cramer-von Mises test statistics
respectively Anderson-Darling test statistics which are frequently
applied for goodness-of-fit investigations;
some more recent methods deal with (other kinds of) 
cumulative paired divergences 
and closely related concepts.
In this paper, 
we provide a general framework which covers in particular both 
the above-mentioned density-based and distribution-function-based 
divergence approaches; 
the dissimilarity of 
quantiles respectively of other statistical functionals
will be included as well. 
From this framework, we structurally 
extract numerous classical and also state-of-the-art
(including new) procedures. Furthermore, we deduce new
concepts of dependence between random variables,
as alternatives 
to the celebrated mutual information.
Some variational representations are discussed, too.
\end{abstract}

%
\section{Divergences, Statistical Motivations and Connections to Geometry}
%

\subsection{Basic Requirements on Divergences (Directed Distances)}
\label{subsec.1b.divergences}

For a first view, let $P$ and $Q$ be two probability distributions
(probability measures). For those, we would like to employ
real-valued indices $D(P,Q)$ which quantify the ``distance'' 
(respectively dissimilarity, proximity, closeness, discrepancy, discrimination)
between $P$ and $Q$. Accordingly, we require $D(\cdot,\cdot)$ to have the following 
reasonable ``minimal/coarse/wide'' properties
\begin{enumerate}

\item[{\ \ }  (D1)] $D\Big(P,Q\Big) \geq 0$ \quad for all $P$, $Q$ under investigation \hspace{0.7cm} (nonnegativity),
\item[{\ \ }  (D2)] $D\Big(P,Q\Big) = 0$ \ if and only if \ 
$P=Q$ \quad (reflexivity; 
\ identity of indiscernibles\footnote{see e.g. Weller-Fahy et al.~\cite{Wel:15}}),

\end{enumerate}
and such $D(\cdot,\cdot)$ is then called a \textbf{divergence} (in the narrow sense)
or disparity or contrast function. 
Basically,
the divergence $D\Big(P,Q\Big)$ of $P$ and $Q$ can be interpreted as a kind of 
``\textbf{directed distance} from $P$ to $Q$'';
the corresponding directness stems from the fact that
in general one has the asymmetry $D\Big(P,Q\Big) \ne D\Big(Q,P\Big)$.
This can turn out to be especially useful in contexts
where the first distribution $P$
is always/principally of ``more importance'' 
or of ``higher attention''
than the second 
distribution $Q$; moreover, it can technically happen that $D\Big(P,Q\Big) < \infty$
but $D\Big(Q,P\Big) = \infty$, for instance in practically important
applications within a (say) discrete context where
$P$ and $Q$ have different zero-valued probability masses (e.g. zero observations),
see e.g. the discussion in Subsection \ref{subsec.1c.zeros} below.

Notice that we don't assume that the triangle inequality holds for $D(\cdot,\cdot)$. 

\subsection{Some Statistical Motivations}

To start with, let us consider probability distributions $P$ and $Q$ 
having \textit{strictly positive} density functions (densities) $f_P$ and $f_Q$ 
with respect to some measure $\lambda$ on some (measurable) space $\mathcal{X}$. 
For instance, if $\lambda := \lambda_{L}$ is the Lebesgue measure on (some subset of) 
$\mathcal{X}=\mathbb{R}$ then $f_P$ and $f_Q$ are ``classical'' 
(e.g. Gaussian) density functions; in contrast, in the \textit{discrete setup} 
where $\mathcal{X} := \mathcal{X}_{\#}$ has 
countably many elements and is equipped with the counting measure
$\lambda := \lambda_{\#} := \sum_{z \in \mathcal{X}_{\#}} \delta_{z}$ 
(where $\delta_{z}$ is Dirac's 
one-point distribution
$\delta_{z}[A] :=  \boldsymbol{1}_{A}(z)$
(where here and in the sequel $\boldsymbol{1}_{A}(\cdot)$ which stands 
for the indicator function of a set $A$), 
and thus $\lambda_{\#}[\{z\}] =1$ for all $z \in \mathcal{X}_{\#}$),
then $f_P$ and $f_Q$ are probability mass functions (counting-density functions,
relative-frequency functions, frequencies).

\vspace{0.3cm}
\noindent
For such kind of probability measures $P$ and $Q$, let us start with 
the widely used class $D_{\phi}(\cdot,\cdot)$
of \textit{Csiszar-Ali-Silvey-Morimoto (CASM) divergences} (see \cite{Csi:63},\cite{Ali:66},\cite{Mori:63}) 
which are usually abbreviatorily called \textit{$\phi-$divergences} and which are defined by
\begin{eqnarray} 
& & \hspace{-0.2cm} 
0 \leq D_{\phi}(P,Q) := \int_{{\mathcal{X}}} 
f_{Q}(x) \cdot \phi \negthinspace \left( {\frac{f_{P}(x)}{f_{Q}(x)}}\right) 
\, \mathrm{d}\lambda(x) \ , 
\label{BroStuHB22:fo.def.47d.reduced} \\
& & \hspace{-0.2cm} 
=  \int_{{\mathcal{X}}} 
\phi \negthinspace \left( {
\frac{f_{P}(x)}{f_{Q}(x)}}\right) 
\, \mathrm{d}Q(x) \ , 
\label{BroStuHB22:fo.def.47d.reduced2}
\end{eqnarray}
where $\phi: \ ]0,\infty[ \ \mapsto [0,\infty[$ is a convex function which is strictly convex at $1$
and which satisfies $\phi(1)=0$.  
It can be easily seen that this $D_{\phi}(\cdot,\cdot)$
satisfies the above-mentioned requirements/properties/axioms $(D1)$ and $(D2)$.
In the above-mentioned discrete setup with $\mathcal{X}= \mathcal{X}_{\#}$, \eqref{BroStuHB22:fo.def.47d.reduced} turns into
\begin{eqnarray} 
& & \hspace{-0.2cm} 
0 \leq D_{\phi}(P,Q) = \sum_{x \in \mathcal{X}_{\#}} 
f_{Q}(x) \cdot \phi \negthinspace \left( {
\frac{f_{P}(x)}{f_{Q}(x)}}\right)  \ ,
\nonumber
\end{eqnarray}
whereas in the above-mentioned real-valued absolutely-continuous case,
the integral in \eqref{BroStuHB22:fo.def.47d.reduced} reduces 
(except for rare cases) to a classical Riemann integral with integrator
$\mathrm{d}\lambda_{L}(x) = \mathrm{d}x$.
Notice that --- depending on $\mathcal{X}$, $\phi$ etc. --- 
the divergence $D_{\phi}(P,Q)$ in  \eqref{BroStuHB22:fo.def.47d.reduced} may become $\infty$.
For comprehensive treatments
of $\phi-$divergences (CASM divergences), 
the reader is referred to e.g.
Liese \& Vajda~\cite{Lie:87},  Read \& Cressie~\cite{Rea:88}, Vajda~\cite{Vaj:89},
Liese \& Vajda~\cite{Lie:06}, Pardo~\cite{Par:06}, Liese \& Miescke~\cite{Lie:08}, and Basu et al.~\cite{Bas:11}.
Important prominent special cases of \eqref{BroStuHB22:fo.def.47d.reduced}
are 
the omnipresent Kullback-Leibler divergence/distance (relative entropy) with $\phi_{KL}(t) := t \log(t) + 1-t$
and thus
\begin{equation}
D_{\phi_{KL}}(P,Q) = \int_{{\mathcal{X}}} 
f_{P}(x) \cdot \log\negthinspace \left( {\frac{f_{P}(x)}{f_{Q}(x)}}\right) 
\, \mathrm{d}\lambda(x) \ , 
\nonumber
\end{equation}
the reverse Kullback-Leibler divergence/distance with $\phi_{RKL}(t) := - \log(t) + t -1$
and hence
\begin{equation}
D_{\phi_{RKL}}(P,Q) = \int_{{\mathcal{X}}} 
f_{Q}(x) \cdot \log\negthinspace \left( {\frac{f_{Q}(x)}{f_{P}(x)}}\right) 
\, \mathrm{d}\lambda(x) = D_{\phi_{KL}}(Q,P) \ , 
\label{BroStuHB22:fo.RKL.reduced}
\end{equation}
(half of) Pearson's $\chi^2-$distance with $\phi_{PC}(t) := \frac{(t-1)^2}{2}$
and consequently
\begin{equation}
D_{\phi_{PC}}(P,Q) = \frac{1}{2} \int_{{\mathcal{X}}} 
\frac{(f_{P}(x)-f_{Q}(x))^2}{f_{Q}(x)} 
\, \mathrm{d}\lambda(x) \ , 
\label{BroStuHB22:fo.PC.reduced}
\end{equation}
(half of) Neyman's $\chi^2-$distance with $\phi_{NC}(t) := \frac{(t-1)^2}{2}$ and thus
\begin{equation}
D_{\phi_{PC}}(P,Q) = \frac{1}{2} \int_{{\mathcal{X}}} 
\frac{(f_{P}(x)-f_{Q}(x))^2}{f_{P}(x)} 
\, \mathrm{d}\lambda(x) \ , 
\nonumber
\end{equation}
the (double of squared) Hellinger distance 
--- also called (half of) Freeman-Tukey divergence ---
with $\phi_{HD}(t) := 2 (\sqrt{t} \, - \, 1)^2$
and hence
\begin{equation}
D_{\phi_{PC}}(P,Q) = 2 \int_{{\mathcal{X}}} 
\left(\sqrt{f_{P}(x)} \, - \, \sqrt{f_{Q}(x)} \, \right)^2 
\, \mathrm{d}\lambda(x) \ , 
\nonumber
\end{equation}
the total variation distance with $\phi_{TV}(t) := \vert t-1 \vert$
and consequently
\begin{equation}
D_{\phi_{TV}}(P,Q) = \int_{{\mathcal{X}}} 
\left\vert f_{P}(x) \, - \, f_{Q}(x) \, \right\vert 
\, \mathrm{d}\lambda(x) \ , 
\nonumber
\end{equation}
and the power divergences $D_{\phi_{\alpha}}(P,Q)$ 
(also known as alpha-divergences, Cressie-Read measures/distances, and Tsallis
cross-entropies) with 
$\phi_{\alpha}(t) := 
\frac{t^\alpha-\alpha \cdot t+ \alpha - 1}{\alpha \cdot (\alpha-1)}$ 
($\alpha \in \mathbb{R}\backslash\{0,1\}$). 
Notice that (in the current setup of probability distributions with zero-free density functions) 
$D_{\phi_{PC}}(P,Q)$ resp. $D_{\phi_{NC}}(P,Q)$ 
resp. $D_{\phi_{HD}}(P,Q)$ are equal to $D_{\phi_{\alpha}}(P,Q)$ with $\alpha=2$ 
resp. $\alpha=-1$ resp. $\alpha=2$, and that one can prove
$D_{\phi_{KL}}(P,Q) = \lim_{\alpha \uparrow 1} D_{\phi_{\alpha}}(P,Q) =: D_{\phi_{1}}(P,Q)$
as well as $D_{\phi_{KL}}(P,Q) = \lim_{\alpha \downarrow 0} D_{\phi_{\alpha}}(P,Q)
=: D_{\phi_{0}}(P,Q)$; henceforth, we will use 
this comfortable continuous embedding to a divergence family 
$\big(D_{\phi_{\alpha}}(P,Q)\, \big)_{\alpha \in \mathbb{R}}$
which covers important special cases.

\vspace{0.3cm}
\noindent
From a statistical standpoint, the definition \eqref{BroStuHB22:fo.def.47d.reduced} 
finds motivation in the far-reaching approach by Ali \& Silvey \cite{Ali:66}: 
by noting that in a simple model
--- where a random variable $X$ takes values on a finite discrete set $\mathcal{X} = \mathcal{X}_{\#}$
and its 
distribution is either $P$ or $Q$ having probability mass function
$f_{P}$ or $f_{Q}$ ---
the statistics $\frac{f_{P}(X)}{f_{Q}(X)}$ is a sufficient statistics
(meaning that $P\big(X=x \, \big\vert \, \frac{f_{P}(X)}{f_{Q}(X)} =t\Big) = Q\big(X=x \, \big\vert \, \frac{f_{P}(X)}{f_{Q}(X)} =t\Big)$
for all $x$ and $t$)
they argue that any measurement aiming at
inference on the distribution of $X$ should be a function of the likelihood
ratio $LR := \frac{f_{P}(X)}{f_{Q}(X)}$. Thus, a real-valued coefficient $D(P,Q)$ 
of closeness/dissimilarity between $P$ and $Q$ 
should be considered as an aggregation/expectation ---
over some measure (typically 
$P$ or $Q$) ---
of a function $\phi$ of LR,
hence formally leading to \eqref{BroStuHB22:fo.def.47d.reduced}
with not necessarily convex function $\phi$.
This construction is compatible with the following set of four axioms/requirements 
which bear some fundamentals for the construction of a discrimination index between
distributions, and which (amongst other things) imply the convexity of $\phi$:

\begin{enumerate}

\item[(A1)] $D_{\phi}(P_{1},P_{2})$ should be defined for all 
pairs of probability distributions 
$P_{1},P_{2}$ on the same sample space $\mathcal{X}$.

\item[(A2)] Let $x \mapsto t(x)$ a measurable transformation from 
$\left( \mathcal{X}\text{,} \mathcal{F}\right) $ onto a measure space 
$\left( \mathcal{Y}\text{,} \mathcal{G}\right) $ then there should hold
\begin{equation}
\ D_{\phi}(P_{1},P_{2}) \geq \ D_{\phi }(P_{1}t^{-1},P_{2}t^{-1}),
\label{datproc} 
\end{equation}
where $P_{i}t^{-1}$ denotes the induced measure on $\mathcal{Y}$
corresponding to $P_{i}$. Notice that \eqref{datproc}
is called \textit{data processing inequality} or 
\textit{information processing inequality},
and --- as shown in \cite{Ali:66} --- 
it implies that $\phi$ should be a convex function.

\item[(A3)] $D_{\phi }(P_{1},P_{2})$ should take its minimum value when $P_{1}=P_{2}$
and its maximum value when $P_{1}\perp P_{2}$ (i.e., $P_{1}$ and $P_{2}$ are singular,
in the sense that the supports of the distributions $P_{1}$ and $P_{2}$ do not overlap (are disjoint)).

\item[(A4)] A further axiom of statistical nature should be satisfied in relation
with a statistical notion of separability of two distributions in a given
model.\\
Assume that for a given family of parametric
distributions $\left( P_{\theta }\right) _{\theta \in \Theta}$ and for any small risk $
\alpha $ the following property holds: if $P_{\theta_{0}}$ is rejected vs. 
$P_{\theta _{1}}$ with risk $\leq $ $\alpha$ optimally (Neyman-Pearson
approach), then $P_{\theta _{0}}$ is rejected vs. $P_{\theta _{2}}$ with risk 
$\leq $ $\alpha $ 
(meaning $P_{\theta _{3}}$ is further away from $P_{\theta_{1}}$ than $P_{\theta _{2}}$ is).\\
\noindent
Then one should have
\[
D_{\phi }\left( P_{\theta _{0}},P_{\theta _{2}}\right) \geq D_{\phi
}\left( P_{\theta _{0}},P_{\theta _{1}}\right) . 
\]

\end{enumerate}

\noindent
Notice that in (A4) we use a slight extension of the original requirements 
of \cite{Ali:66} (who employ a monotone likelihood ratio concept).

\vspace{0.3cm}
\noindent
As a second use-of-divergence incentive stemming from considerations in statistics
(as well as in the adjacent research fields of information theory, artificial
intelligence and machine learning), 
we mention parameter estimation in terms of $\phi-$divergence minimization.
For this, let $Y$ be a random variable taking values in a finite
discrete space $\mathcal{X} := \mathcal{X}_{\#}$, and let
$f_{P}(x) = P[Y=x]$ be its strictly positive probability mass function 
under an unknown hypothetical law $P$.
Moreover, we assume that $P$ lies in --- respectively can be approximated by ---
a model $\Omega := \{ Q_{\theta}: \, \theta \in \Theta \}$ ($\Theta \subset \mathbb{R}$) 
being a class of finite discrete parametric distributions having strictly positive
probability mass functions $f_{Q_{\theta}}$ on $\mathcal{X}_{\#}$.
Moreover, let $P_{N}^{emp} := \frac{1}{N} \cdot \sum_{i=1}^{N} \delta_{Y_{i}}[\cdot]$
be the well-known data-derived empirical distribution/measure 
of an $N-$size independent and identically distributed (i.i.d.) sample/observations $Y_1, \ldots, Y_N$ of $Y$;
the according probability mass function is $f_{P_{N}^{emp}}(x) = 
\frac{1}{N} \cdot \# \{ i \in \{ 1, \ldots, N\}: Y_i =x \} $
which reflects the underlying (normalized) histogram; here, as usual, $\# A$ denotes the number of elements
in a set $A$. In the following, we assume that the sample size $N$ is large enough
such that $f_{P_{N}^{emp}}$ is strictly positive (see the next subsection for a relaxation).
\\
If the data-generating distribution $P$ lies in $\Omega$, i.e. $P=Q_{\theta_{tr}}$ for
some ``true'' unknown parameter $\theta_{tr} \in \Theta$, then (under some mild technical assumptions) 
it is easy to show that the corresponding \textit{maximum likelihood estimator (MLE)} 
$\widehat{\theta}$ is EQUAL to
\begin{equation}
\widehat{\widehat{\theta}} := \arg \min_{\theta \in \Theta} D_{\phi_{0}}(Q_{\theta}, P_{N}^{emp} )
\nonumber
\end{equation}
where $\phi_{0} := - \log(t) + t -1$ and $ D_{\phi_{0}}(\cdot, \cdot)$ is the 
the reverse Kullback-Leibler divergence already mentioned above.
Due due its construction, $\widehat{\widehat{\theta}}$ is called 
\textit{minimum reverse-Kullback-Leibler divergence (RKLD) estimator},
and $Q_{\widehat{\widehat{\theta}}}$ is the RKLD-projection of $P_{N}^{emp}$
on $\Omega$.
In the other --- also practically important --- case where $P$ does not lie in the model $\Omega$ 
(but is reasonably ``close'' to it),
i.e. the model is \textit{misspecified},  
then $Q_{\widehat{\widehat{\theta}}}$ is still a reasonable proxy of $P$ if the sample
size $N$ is large enough.

\vspace{0.3cm}
\noindent
In the light of the preceding paragraph, it makes sense to consider the
more general \textit{minimum $\phi-$divergence/distance estimation problem} 
\begin{equation}
\widehat{\widehat{\theta}} := \arg \inf_{\theta \in \Theta} D_{\phi}(Q_{\theta}, P_{N}^{emp} )
\label{MDE1}
\end{equation}
where $\phi$ is not necessarily equal to $\phi_{0}$;
for instance, through some comfortably verifiable criteria on $\phi$ one can end up with an
outcoming minimum $\phi-$divergence/distance estimator $\widehat{\widehat{\theta}}$
which is more robust against outliers than the MLE $\widehat{\theta}$ 
(see e.g. the residual-adjustment-function approach of Lindsay~\cite{Lind:94}, 
its comprehensive treatment in Basu et al.~\cite{Bas:11},  
and the corresponding flexibilizations in Ki{\ss}linger \& Stummer~\cite{Kis:16}, 
Roensch \& Stummer~\cite{Roe:17}). Usually, $\widehat{\widehat{\theta}}$ of \eqref{MDE1} 
is called \textit{minimum $\phi-$divergence estimator (MDE)},
and $Q_{\widehat{\widehat{\theta}}}$ is the $phi-$divergence-projection of $P_{N}^{emp}$
on $\Omega$.

\vspace{0.3cm}
\noindent
A further useful generalization is the 
``distribution-outcome type''
\textit{minimum divergence/distance estimation problem}
\begin{equation}
\widehat{Q} := \arg \inf_{Q \in \Omega} D_{\phi}(Q, P_{N}^{emp} )
\label{MDE2}
\end{equation}
where $P_{N}^{emp}$ stems from a general (not necessarily parametric, unknown)
data generating distribution $P$ and
$\Omega$ may be a ``fairly general'' model 
being a class of finite discrete distributions having strictly positive
probability mass functions $f_{Q}$ on $\mathcal{X}_{\#}$
(and, as usual, \eqref{MDE2} can be rewritten as a minimization
problem on the $(\#\Omega-1)-$dimensional probability simplex).
The outcoming $\widehat{Q}$ of \eqref{MDE2} is still called (distribution-type)
\textit{minimum $\phi-$divergence estimator (MDE)}, and can be interpreted as
$phi-$divergence-projection of $P_{N}^{emp}$ on $\Omega$.
Problem \eqref{MDE2} is in particular beneficial in non- and semi-parametric contexts,
where $\Omega$ reflects (partially) non-parametrizable model constraints. 
For instance, $\Omega$ may consist (only) of constraints on 
moments or on L-moments 
(see e.g. Broniatowski \& Decurninge \cite{Bro:16});
alternatively, $\Omega$ may be e.g. a
tubular neighborhood 
of a parametric model
(see e.g. Liu \& Lindsay \cite{Liu:09},
Ghosh \& Basu \cite{Gho:18}).

\vspace{0.3cm}
\noindent
The closeness --- especially in terms of the sample size $N$ --- of the 
data-derived empirical distribution from the model $\Omega$ is quantified
by the corresponding minimum
\begin{equation}
D_{\phi}(\Omega, P_{N}^{emp} ) := \inf_{Q \in \Omega} D_{\phi}(Q, P_{N}^{emp} )
\label{MDE3}
\end{equation}
of \eqref{MDE2}; thus, it carries useful statistical information.
Moreover, under some mild assumptions, $D_{\phi}(\Omega, P_{N}^{emp} )$ converges to
\begin{equation}
D_{\phi}(\Omega, P) := \inf_{Q\in\Omega} D_{\phi}( Q, P )
\label{inf proba new}
\end{equation}
where $P$ is the (unknown) data generating distribution.
In case of $P \in \Omega$ one obtains $D_{\phi}(\Omega, P)=0$, whereas
for $P \notin \Omega$ the $\phi-$divergence minimum $D_{\phi}(\Omega, P)$
--- and thus its approximation $D_{\phi}(\Omega, P_{N}^{emp} )$ ---
quantifies the \textit{adequacy} of the model
$\Omega$ for modeling $P$;
a lower $D_{\phi}(\Omega, P)-$value means 
a better adequacy
(in the sense of a lower departure between the model and the truth,
cf. Lindsay \cite{Lind:04}, Lindsay et al. \cite{Lind:08},
Markatou \& Sofikitou \cite{Mark:18}, 
Markatou \& Chen \cite{Mark:19}). 

\vspace{0.3cm}
\noindent
Hence, 
especially in the context of \textit{model selection/choice} 
(and the related issue of goodness-of-fit testing) 
within complex big-data contexts, for the \textit{search of appropriate models} $\Omega$ and
model elements/members therein,
the (fast and efficient) computation of  
$D_{\phi}(\Omega, P)$ respectively  $D_{\phi}(\Omega, P_{N}^{emp} )$
constitutes a decisive first step, since
if the latter two are \textquotedblleft  too large\textquotedblright\ (respectively, 
\textquotedblleft  much larger than\textquotedblright\
$D_{\phi}(\overline{\Omega}, P)$ respectively  $D_{\phi}(\overline{\Omega}, P_{N}^{emp} )$ 
for some competing
model $\overline{\Omega}$), then the model 
$\Omega$ is \textquotedblleft not adequate enough\textquotedblright\ (respectively
\textquotedblleft  much less adequate than\textquotedblright\ 
 $\overline{\Omega}$).
For tackling the computation of 
$D_{\phi}(\Omega, P)$ respectively  $D_{\phi}(\Omega, P_{N}^{emp} )$
on fairly general (e.g. high-dimensional, non-conex and even highly disconnected) 
constraint sets $\Omega$,
a ``precise bare simulation'' approach has been recently developed by 
Broniatowski \& Stummer \cite{Bro:21a}.

\vspace{0.3cm}
\noindent
For the sake of a compact first glance, 
in this subsection we have mainly dealt with finite discrete distributions $P$ and $Q$
having zeros-free probability mass functions. However, with appropriate technical care,
one can extend the above concepts also to general discrete distributions
with zeros-carrying probability mass functions and even to non-discrete (e.g.
absolutely continuous) distributions with zeros-carrying density functions. 
(Only) The correspondingly necessary
generalization of the basic $\phi-$divergence definition \eqref{BroStuHB22:fo.def.47d.reduced}
is addressed in the next subsection.

\subsection{Incorporating density function zeros }
\label{subsec.1c.zeros}

Recall that in our first basic $\phi-$divergence definition 
\eqref{BroStuHB22:fo.def.47d.reduced},\eqref{BroStuHB22:fo.def.47d.reduced2}
we have employed probability distributions $P$ and $Q$ 
having strictly positive density functions $f_P$ and $f_Q$ 
with respect to some measure $\lambda$ on some (measurable) space $\mathcal{X}$,
and consequently $P$ and $Q$ are equivalent.
However, in many applications one has to allow $f_P$ and/or $f_Q$ to have zero values.
For instance, in the above-mentioned empirical distribution $P_{N}^{emp}$
for small/medium sample size $N$ 
(or even large sample size for rare-events)
one may have $f_{P_{N}^{emp}}(\widetilde{x})=0$ for some $\widetilde{x}$\footnote{which corresponds
to an empty histogram cell at $\widetilde{x}$}, even though
the candidate-model probability mass satisfies $f_{Q_{\theta}}(\widetilde{x}) \ne 0$
for some $\theta \in \Theta$\footnote{if
$f_{Q_{\theta}}(\widetilde{x}) = 0$ for all $\theta \in \Theta$ one
should certainly reduce the space $\mathcal{X}$ by removing $\widetilde{x}$}.

\vspace{0.cm}
\noindent
Accordingly, we employ the following extension: 
for probability distributions $P$ and $Q$ having density functions $f_P$ and $f_Q$ 
with respect to some measure $\lambda$ on some (measurable) space $\mathcal{X}$
one defines the \textit{Csiszar-Ali-Silvey-Morimoto (CASM) divergences}  
--- in short \textit{$\phi-$divergences} --- by 
\begin{eqnarray} 
& & \hspace{-0.2cm} 
0 \leq D_{\phi}(P,Q) := 
\int_{\{f_{P} \cdot f_{Q} > 0\}} 
\phi \negthinspace \left( {
\frac{f_{P}(x)}{f_{Q}(x)}}\right) 
\, \mathrm{d}Q(x) 
\nonumber\\
& & \hspace{2.2cm} 
+  \, \phi(0) \cdot Q[f_{P} =0] \, + \, \phi^{*}(0) \cdot P[f_{Q} =0]
\label{BroStuHB22:fo.def.47d.full} \qquad \ \\[0.1cm]
& & \hspace{-0.2cm}
\textrm{with} \quad \phi(0) \cdot 0 = 0 \quad \textrm{and} \quad
\phi^{*}(0) \cdot 0 = 0 \ 
\label{BroStuHB22:fo.def.47d.full2}
\end{eqnarray}
(see e.g. Liese \& Vajda \cite{Lie:06}).
Here, we have employed
(as above) $\phi: \ ]0,\infty[ \ \mapsto [0,\infty[$ to be a convex function 
which is strictly convex at $1$
and which satisfies $\phi(1)=0$; moreover, we have used
the (always existing) limits 
$\phi(0) := \lim_{t\downarrow 0} \phi(t) \in ]0,\infty]$ and 
$\phi^{*}(0) := \lim_{t\downarrow 0} \phi^{*}(t) = \lim_{t \rightarrow \infty} \frac{\phi(t)}{t}$ 
of the so-called
$*-$adjoint function $\phi^{*}(t):= t \cdot \phi(\frac{1}{t})$ \, ($t \in ]0,\infty[$). 
It can be proved that $D_{\phi}(\cdot,\cdot)$ satisfies the above-mentioned 
requirements/properties/axioms $(D1)$ and $(D2)$; even more, one gets 
the following range-of-value assertion
(cf. Csiszar \cite{Csi:63}, \cite{Csi:67} and Vajda \cite{Vaj:72}, 
see e.g. also Liese \& Vajda \cite{Lie:06}):

\begin{theorem}
\label{BroStuHB22:thm.range}
There holds
\begin{eqnarray} 
& & \hspace{-0.2cm} 
0 \ \leq \ D_{\phi}(P,Q) \, \leq \, \phi(0) \, + \, \phi^{*}(0) 
\qquad \textrm{for all $P$, $Q$}
\nonumber 
\end{eqnarray}
where (i) the left equality holds only for $P=Q$, and 
(ii) the right equality holds always
for $P \perp Q$ (singularity, i.e. the zeros-set of $f_{P}$ is disjoint from  
the zeros-set of $f_{Q}$) and only for $P \perp Q$
in case of $\phi(0) + \phi^{*}(0) < \infty$.
 
\end{theorem}

\noindent
A generalization of Theorem \ref{BroStuHB22:thm.range} 
to the context of finite (not necessarily probability) measures $P$ and $Q$
is given in Stummer \& Vajda \cite{Stu:10}; for instance, in a two-sample test situation 
$P$ and $Q$ may be two generalized empirical distributions which reflect 
non-normalized (rather than normalized) histograms. 

\vspace{0.3cm}
\noindent
As an example, let us illuminate the upper bounds 
$\phi(0) + \phi^{*}(0)$ of --- the zeros-incorpating versions of --- 
of the above-mentioned important power divergence family
$\big(D_{\phi_{\alpha}}(P,Q)\, \big)_{\alpha \in \mathbb{R}}$
with $\phi_{\alpha}(t) := 
\frac{t^\alpha-\alpha \cdot t+ \alpha - 1}{\alpha \cdot (\alpha-1)}$ 
($\alpha \in \mathbb{R}\backslash\{0,1\}$), 
$\phi_{1}(t) := \phi_{KL}(t) = t \log(t) + 1-t$
and $\phi_{0}(t) := \phi_{RKL}(t) := - \log(t) + t -1$. It is easy to see that
for $P \perp Q$ one gets
\begin{eqnarray}
\phi_{\alpha}(0) = \phi^{*}_{1-\alpha}(0) &=&
\begin{cases}
\ \infty , 
\qquad \textrm{if } \alpha \leq 0, \\
\ \frac{1}{\alpha}, \hspace{0.7cm} \textrm{if } \alpha > 0,
\end{cases}
\label{BroStuHB22:fo.range.bound1} 
\end{eqnarray}
and hence
\begin{eqnarray}
D_{\phi_{\alpha}}(P,Q) = \phi_{\alpha}(0) + \phi_{\alpha}^{*}(0) &=&
\begin{cases}
\ \infty , \hspace{0.8cm}
\qquad \textrm{if } \alpha \notin \, ]0,1[, \\
\ \frac{1}{\alpha \cdot (1-\alpha)}, \hspace{0.7cm} \textrm{if } \alpha \in \, ]0,1[.
\end{cases}
\label{BroStuHB22:fo.range.bound2} 
\end{eqnarray}

\noindent
Especially, for $P \perp Q$ one gets for the Kullback-Leibler divergence
$D_{\phi_{KL}}(P,Q) = D_{\phi_{1}}(P,Q) = \infty$ whereas  
$D_{\phi_{0.99}}(P,Q) = \frac{10000}{99}$
one achieves a finite value; thus, 
in order to avoid infinities
it is more convenient to work with the well-approximating divergence generator
$\phi_{0.99}$ of $\phi_{1}$.
Similarly, for the reverse Kullback-Leibler divergence we obtain
$D_{\phi_{RKL}}(P,Q) = D_{\phi_{0}}(P,Q) = \infty$ whereas  
$D_{\phi_{0.01}}(P,Q) = \frac{10000}{99}$.
Furthermore, for $P \perp Q$ one gets for Pearson's $\chi^2-$divergence
$D_{\phi_{2}}(P,Q) = \infty$, for Neyman's $\chi^2-$divergence
$D_{\phi_{-1}}(P,Q) = \infty$ and for the (squared) Hellinger distance
$D_{\phi_{1/2}}(P,Q) = 4$.

\vspace{0.3cm}
\noindent
Returning to the general context,
notice that the upper bound $\phi(0) + \phi^{*}(0)$ in Theorem \ref{BroStuHB22:thm.range} 
is independent of $P$ and $Q$, and thus $D_{\phi}(P,Q)$ is
of no discriminative use in statistical situations where $P$ and $Q$ are singular
(i.e. $P \perp Q$). This is the case, for instance, in the following commonly encountered 
``crossover'' context:

\begin{enumerate}

\item[(CO1)] 
$Y$ is an univariate (absolutely continuous) random variable with 
unknown hypothetical probability distribution $P$
having strictly positive 
density function $f_{P}$ 
with respect to the Lebesgue measure 
$\lambda_{L}$ on $\mathcal{X} = \mathbb{R}$
(recall that this means that $f_{P}$ is a 
``classical'' (e.g. Gaussian) probability density function),

\item[(CO2)] the corresponding model
$\Omega := \{ Q_{\theta}: \, \theta \in \Theta \}$ ($\Theta \subset \mathbb{R}$) 
is a class of parametric distributions having strictly positive
probability density functions $f_{Q_{\theta}}$ with respect to $\lambda_{L}$, and

\item[(CO3)] $P_{N}^{emp} := \frac{1}{N} \cdot \sum_{i=1}^{N} \delta_{Y_{i}}[\cdot]$
is the data-derived empirical distribution
of an $N-$size independent and identically distributed (i.i.d.) sample/observations $Y_1, \ldots, Y_N$ of $Y$;
recall that the according probability mass function is $f_{P_{N}^{emp}}(x) = 
\frac{1}{N} \cdot \# \{ i \in \{ 1, \ldots, N\}: Y_i =x \} $
which is the density function with respect to the counting measure $\lambda_{\#}$
on the distinct values of the sample. 

\end{enumerate}

\noindent
This contrary density-function behaviour can be put in an encompassing framework 
by employing the joint density-building (i.e. dominating) measure 
$\lambda := \lambda_{L} + \lambda_{\#}$.
Clearly, one always has the singularity $P_{N}^{emp} \perp Q_{\theta}$
and thus, due to Theorem \ref{BroStuHB22:thm.range} one gets
\begin{eqnarray} 
& & \hspace{-0.2cm} 
D_{\phi}(Q_{\theta},P_{N}^{emp}) = \phi(0) + \phi^{*}(0)
\ \textrm{for all $\theta \in \Theta$},
\qquad  
\inf_{\theta \in \Theta} D_{\phi}(Q_{\theta}, P_{N}^{emp} ) = \phi(0) + \phi^{*}(0) \ . 
\qquad \ 
\label{BroStuHB22:fo.def.47d.reduced3}
\end{eqnarray}
Accordingly, in such a situation one can not obtain a corresponding
minimum $\phi-$divergence estimator.

\vspace{0.3cm}
\noindent
Also notice that for power divergences $D_{\phi_{\alpha}}(P,Q)$ 
with $\alpha \notin \, ]0,1[$ it can happen that
$D_{\phi_{\alpha}}(P,Q) = \infty$
even though $P$ and $Q$ are not singular
(which due to \eqref{BroStuHB22:fo.range.bound2} 
 is consistent with Theorem \ref{BroStuHB22:thm.range}).
For instance, consider a situation with two different i.i.d. samples 
$Y_1, \ldots, Y_N$ of $Y$ having distribution $P$ and $\widetilde{Y}_1, \ldots, \widetilde{Y}_M$ of 
$\widetilde{Y}$ having distribution $Q$ with (say) $Q \sim P$ (equivalence);
in terms of the corresponding empirical distributions
$P_{N}^{emp} := \frac{1}{N} \cdot \sum_{i=1}^{N} \delta_{Y_{i}}[\cdot]$
and $\widetilde{P}_{M}^{emp} := \frac{1}{N} \cdot \sum_{i=1}^{N} \delta_{\widetilde{Y}_{i}}[\cdot]$
one obtains $D_{\phi_{\alpha}}(P_{N}^{emp},\widetilde{P}_{M}^{emp}) = \infty$
if the set of zeros of the corresponding probability mass function $f_{P_{N}^{emp}}$ is strictly 
larger (for $\alpha \leq 0$) respectively smaller (for $\alpha \geq 1$)
than the set of zeros of $f_{\widetilde{P}_{M}^{emp}}$
(i.e. $\widetilde{P}_{N}^{emp}[f_{P_{M}^{emp}} =0 \, ] >0$
respectively $P_{N}^{emp}[f_{\widetilde{P}_{M}^{emp}} =0 \, ] >0$),  
to be seen by applying \eqref{BroStuHB22:fo.def.47d.full}, \eqref{BroStuHB22:fo.def.47d.full2},
\eqref{BroStuHB22:fo.range.bound1}.
As above, in such a non-singular situation it is e.g. better to use the
(in fact, even sample-dependent !) power divergence
$D_{\phi_{0.99}}(P_{N}^{emp},\widetilde{P}_{M}^{emp})$ instead of the
Kullback-Leibler divergence $D_{\phi_{1}}(P_{N}^{emp},\widetilde{P}_{M}^{emp}) = \infty$.
Similar infinity-effects can be constructed for the above-mentioned other important
special cases 
$\alpha=0$ (reverse Kullback-Leibler divergence), $\alpha=2$ (Pearson's $\chi^2-$divergence),  
$\alpha=-1$ (Neyman's $\chi^2-$divergence) whereas
for the case $\alpha=1/2$ (square Hellinger distance) everything works out well.
Such an approach serves as an alternative to the approach of 
``lifting/unzeroing/adjusting'' (from sampling randomly appearing) zero probability 
masses\footnote{e.g. which correspond to empty cells in sampled histograms, 
e.g. for rare events and small-or medium-sized sample sizes}  
by pseudo-counts or ``smoothing (in a discrete sense)'', see 
e.g. Fienberg \& Holland \cite{Fie:70},
as well as e.g. Section 4.5 (respectively Section 3.5) in
Jurafsky \& Martin \cite{Jur:09and22} and the references therein.

\vspace{0.3cm}
\noindent
Next, we briefly indicate two ways to circumvent
the problem described in the above-mentioned crossover context 
(CO1),(CO2),(CO3):

\begin{enumerate}

\item[(GR)] grouping (partitioning, quantization) of data: \  
convert\footnote{
in several situations, such a conversion can appear in a natural way;
e.g. an institution may generate/collect data of ``continuous value'' but 
mask them for external data analysts to group-frequencies, for reasons of confidentiality 
(information asymmetry)
}
the model $\Omega$
into a purely discrete context, by
subdividing the data-point-set $\mathcal{X} = \bigcup_{j=1}^{s} A_{j}$
into countably many --
(say) $s \in \mathbb{N}\cup\{\infty\} \backslash\{1\}$ -- 
(measurable) disjoint classes $A_{1}, \ldots, A_{s}$ with the property 
$\lambda_{L}[A_{j}] >0$ 
(``essential partition'');
proceed as in above general discrete subsetup with 
$\mathcal{X}^{new} := \{A_{1}, \ldots, A_{s}\}$ 
and thus the $i$-th data observation $Y_{i}(\omega)$ 
and the corresponding running variable $x$ 
manifest (only) the
corresponding class-membership
(see e.g. Vajda \& van der Meulen \cite{Vaj:11} for a survey on different choices).
Some corresponding
thorough statistical investigations (such as efficiency, robustness, types of grouping, 
grouping-error sensitivity, etc.) 
of the corresponding minimum-$\phi-$divergence-estimation can be found e.g.
in Victoria-Feser \& Ronchetti~\cite{Vic:97},
Menendez et al.~\cite{Men:98,Men:01a,Men:01b},
Morales et al.~\cite{Mor:04,Mor:06},
Lin \& He~\cite{Lin2:06}. 

\vspace{0.1cm}

\item[(SM)] smoothing of the empirical density function:
convert everything to a purely continuous context,
by keeping the original data-point-set $\mathcal{X}$
and by ``continuously modifying'' (e.g. with the help of kernels)
the empirical density function $f_{P_{N}^{emp}}(\cdot)$
to a function $f_{P_{N}^{emp,smo}}(\cdot) > 0$ (a.s.) such that 
$\int_{\mathcal{X}} f_{P_{N}^{emp,smo}}(x) \, \mathrm{d}\lambda_{L}(x) =1$.
Some corresponding thorough statistical investigations 
(such as efficiency, robustness, information loss, etc.) 
of the corresponding
minimum-$\phi-$divergence-estimation can be found e.g. in 
Beran \cite{Beran:77}, 
Basu \& Lindsay~\cite{Bas:94}, Park \& Basu~\cite{Park2:04},
Chapter 3 of Basu et al.~\cite{Bas:11},
Kuchibhotla \& Basu~\cite{Kuc:15},  
Al Mohamad~\cite{AlM:18},
and the references therein.

\end{enumerate}

\vspace{0.3cm}
\noindent
In contrast to the above, let us now encounter a crossover situation 
where (CO1) and (CO3) still hold, but the parametric-model-assumption (CO2) is replaced by

\begin{enumerate}

\item[(CO2')] the corresponding model
$\Omega := \{ Q: \, Q \textrm{ satisfies some nonparametric constraints} \}$ 
is a class of distributions $Q$
which contains both (i) distributions $Q$
having strictly positive
probability density functions $f_{Q}$ with respect to $\lambda_{L}$,
as well as (ii) 
all ``context-specific appropriate'' finite discrete distributions $Q$
(having ideally the same (or at least, smaller or equal) support as $P_{N}^{emp}$).

\end{enumerate}

\noindent
The subclasses of $Q \in \Omega$ which satisfy (i) respectively (ii)
are denoted by $\Omega^{ac}$ respectively $\Omega^{dis}$.
Widely applied special cases of (CO2') are nonparametric contexts where $\Omega$ 
is the class of all distributions on $\mathcal{X}=\mathbb{R}$ 
satisfying pregiven moment conditions.
Suppose, that we are interested in the 
corresponding model-adequacy problem (cf. \eqref{inf proba new})
\begin{equation}
D_{\phi}(\Omega^{ac}, P) := \inf_{Q\in\Omega^{ac}} D_{\phi}( Q, P )
\label{BroStuHB22:fo.emplik1}
\end{equation}
where $P$ is the (unknown) data generating distribution (cf. (CO1)).
Recall that in case of $P \in \Omega^{ac}$ one obtains 
$D_{\phi}(\Omega^{ac}, P)=0$, whereas
for $P \notin \Omega^{ac}$ the $\phi-$divergence minimum $D_{\phi}(\Omega^{ac}, P)$
quantifies the adequacy of the model $\Omega^{ac}$ for modeling $P$;
a lower $D_{\phi}(\Omega^{ac}, P)-$value means a better adequacy.
Since in the current setup the empirical distribution $P_{N}^{emp}$ of (CO3) satisfies
$P_{N}^{emp} \perp Q$ for all $Q \in \Omega^{ac}$ we obtain
(analogously to \eqref{BroStuHB22:fo.def.47d.reduced3})

\begin{eqnarray} 
& & \hspace{-0.2cm} 
D_{\phi}(Q,P_{N}^{emp}) = \phi(0) + \phi^{*}(0)
\qquad \textrm{for all $Q \in \Omega^{ac}$},
\nonumber\\ 
& & \hspace{-0.2cm} 
D_{\phi}(\Omega^{ac}, P_{N}^{emp}) := \inf_{Q \in \Omega^{ac}} D_{\phi}(Q, P_{N}^{emp} ) 
\label{BroStuHB22:fo.emplik2}
\\ 
& & \hspace{2.3cm}
= \phi(0) + \phi^{*}(0) \ . 
\qquad \ 
\nonumber
\end{eqnarray}
Hence, statistically it makes no sense to approximate \eqref{BroStuHB22:fo.emplik1}
by \eqref{BroStuHB22:fo.emplik2}. Let us discuss an appropriate alternative,
e.g. for the case 
of the reverse Kullback-Leibler divergence $D_{\phi_{0}}( Q, P )$
with generator $\phi_{0}(t) = \phi_{RKL}(t) = - \log(t) + t -1$ (cf.\eqref{BroStuHB22:fo.RKL.reduced}).
By \eqref{BroStuHB22:fo.range.bound1}, we have $\phi_{0}(0)= \infty$ as well as $\phi_{0}^{*}(0)=1$
and thus $\phi_{0}(0) + \phi_{0}^{*}(0) = \infty$ as well as 
(by \eqref{BroStuHB22:fo.def.47d.full}, \eqref{BroStuHB22:fo.def.47d.full2})
\begin{eqnarray} 
& & \hspace{-0.2cm} 
D_{\phi_{0}}(Q,P_{N}^{emp}) = 
\int_{\{f_{Q} \cdot f_{P_{N}^{emp}} > 0\}} 
\phi_{0} \negthinspace \left( {\frac{f_{Q}(x)}{f_{P_{N}^{emp}}(x)}}\right) 
\, \mathrm{d}P_{N}^{emp}(x) 
\, + \, \infty \cdot P_{N}^{emp}[f_{Q} =0]
\qquad \ 
\nonumber 
\\[0.1cm]
& & \hspace{-0.2cm}
= \frac{1}{N} \cdot 
\sum_{\{i\in\{1,\ldots,N\}: f_{Q}(Y_{i}) \cdot f_{P_{N}^{emp}}(Y_{i}) > 0\}} 
\phi_{0} \negthinspace \left( {
\frac{f_{Q}(Y_{i})}{f_{P_{N}^{emp}}(Y_{i})}}\right) 
\, + \, + \, \infty \cdot P_{N}^{emp}[f_{Q} =0] \, < \, \infty
\nonumber
\end{eqnarray}
for all $Q$ in $\Omega_{N}^{dis}$ which is defined as the 
class of distributions in $\Omega^{dis}$ such that $Q << P_{N}^{emp}$ 
(and thus $Q[f_{P_{N}^{emp}} =0] =0$); also recall that the last term becomes
$\infty \cdot 0= 0$ in case that $Q$ and $P_{N}^{emp}$ have the same support.
Hence, under the assumption that $\Omega_{N}^{dis}$ is non-void, one can
approximate the $\phi=\phi_{0}-$version of \eqref{BroStuHB22:fo.emplik1} by
\begin{equation}
D_{\phi_{0}}(\Omega_{N}^{dis}, P_{N}^{emp}) := 
\inf_{Q\in\Omega_{N}^{dis}} D_{\phi_{0}}( Q, P_{N}^{emp} )
\nonumber
\end{equation}
This is the basic idea of the divergence-minimization formulation of
the so-called ``empirical likelihood'' principle of Owen \cite{Owe:88}, \cite{Owe:90},
\cite{Owe:01}, which leads to many variations
according to the choice of the divergence generator $\phi$; see e.g.
Baggerly \cite{Bag:98},
Judge \& Mittelhammer \cite{Jud:12}, Bertail et al. \cite{Bertail:14},
and Broniatowski \& Keziou \cite{Bro:12}, and references therein.

\vspace{0.3cm}
\noindent
Other ways to circumvent the 
crossover problem (CO1),(CO2),(CO3) respectively 
(CO1),(CO2'),(CO3) can be found e.g. 
in Section VIII of Liese \& Vajda \cite{Lie:06} and Section 4 of Broniatowski \& Stummer \cite{Bro:19a};
moreover, some variational-representation-method approaches will be discussed in Section \ref{sec.7new} below.

\vspace{0.3cm}
\noindent
As a third statistical incentive, let us mention that
with the help of $\phi-$divergence 
minimization one can build \textit{generalizations} of exponential families
with pregiven sufficient statistics (see e.g. Pelletier \cite{Pel:11}, 
Gayen \& Kumar \cite{Gay:21}).
In the special case of Kullback-Leiber divergence (i.e., the divergence
generator $\phi$ is taken to be $\phi_{1}(t) = \phi_{KL}(t) = t \log(t) + 1-t$)
one ends up with classical exponential families.

\subsection{Some Motivations From Probability Theory}

Another environment where $D_{\phi }(Q,P)$ appears in a natural
way is probability theory, in the area of the large deviation paradigm;
the celebrated Sanov theorem states that, up to technicalities, 
\[
\lim_{n\rightarrow \infty }\frac{1}{n}\log P\left( P_{n}\in \Omega \right)
=-D_{\phi _{1}}(\Omega ,P) 
\]
where $P_{n}$ is the empirical distribution of a sample of $n$ independent
copies under $P$, and $\Omega $ is a class of probability dstributions on 
$\left( \mathcal{X}\text{,}\mathcal{B}\right) $, and $D_{\phi_{1}}(\Omega
,P):=\inf_{Q\in \Omega }D_{\phi _{1}}(Q,P)$. Therefore, the
Kullback-Leibler divergence measures the rate of decay of the chances for 
$P_{n}$ to belong to $\Omega $ as $n$ increases, in case that $P$ does not belong
to $\Omega .$ Other divergences inherit of the same character: assume that
the function $\phi $ is the Fenchel-Legendre transform of a moment
generating function $\Lambda (t)$, namely 
\[
\phi (x)=\sup_{t}tx-\Lambda (t) 
\]
where $\Lambda (t):=\log E[e^{tW}]$ for some random variable $W$ defined on
some arbitrary space.\ With $\left( X_{1},..,X_{n}\right) $ being an i.i.d. sample
under $P$ and $\left( W_{1},..,W_{n}\right) $ being an i.i.d. sample of copies of $W$, 
we define the associated weighted empirical distribution as
\[
P_{n}^{W}:=\frac{1}{n}\sum_{i=1}^{n}W_{i}\delta _{X_{i}} .
\]
The following type of conditional Sanov theorem holds:
\[
\lim_{n\rightarrow \infty }\frac{1}{n}\log P\left( \left. P_{n}^{W}\in
\Omega \right\vert X_{1},..,X_{n}\right) =-D_{\phi }(\Omega ,P) , 
\]
where $\Omega $ is a class of signed measures on $\left( \mathcal{X}\text{,}
\mathcal{B}\right) $ satisfying some regularity assumptions. This result
characterizes $D_{\phi }(\Omega ,P)$ as a rate of escape of $P_{n}^{W}$
from $\Omega $ when $P$ does not belong to $\Omega$. 
We refer to Najim \cite{Naj:02}, Trashorras \& Wintenberger \cite{Tras:14}, and 
Broniatowski \& Stummer \cite{Bro:21a} where the latter consider several applications 
of this result for (deterministic as well as statistical) optimization procedures 
by bootstrap.

\vspace{0.3cm}
\noindent
Of course, there are connections between statistical inferences
and $\phi-$divergence-based large deviations results.
For instance, 
the large deviations properties of (types of) the empirical distribution of a sample
from its parent distribution is the cornerstone for the \textit{asymptotic}
study of tests. In this realm, the $\phi-$divergences
play a significant role while testing between some parametric null 
hypothesis $\theta \in \Theta _{0}$ vs. an alternative
$\eta \in \Theta _{1}$; the corresponding \textit{Bahadur slope} 
of a given test statistics
indicates the decay of its $p-$value under the alternative.
In ``standard'' setups, this is connected to the Kullback-Leibler divergence 
$\inf_{\theta \in \Theta_{0}} D_{\phi_{1}}(P_{\eta},P_{\theta})$ 
(between the alternative $\eta$ and the set of all null hypotheses) 
which qualifies the asymptotic efficiency of the statistics at use; 
see Bahadur \cite{Bah:67} \cite{Bah:71}, Hoadley \cite{Hoa:67},
and also e.g. Groeneboom \& Oosterhoff \cite{Gro:77}, 
and Nikitin \cite{Nik:95}. 
As far as other setups is concerned,
Efron \& Tibshirani \cite{Efr:93} generally suggest the weighted bootstrap as a valuable approach
for testing. 
In some concrete frameworks,
it can be proved that testing in parametric models based
on appropriate weighted-bootstrapped $\phi-$divergence test statistics enjoys 
maximal Bahadur efficiency with
respect to any other weighted-bootstrapped test statistics
(see Broniatowski \cite{Bro:21b}); the
corresponding Bahadur slope is related to the specific weighting procedure,
and substitutes the Kullback-Leibler divergence by  
some other $\phi-$divergence, specific of
the large deviation properties of the weighted empirical distribution.

\subsection{Divergences and Geometry}

For this section, we return to the general framework of Section \ref{subsec.1b.divergences}
where we have defined divergences to satisfy the two properties 
the two properties
\begin{enumerate}

\item[{\ \ }  (D1)] $D\Big(P,Q\Big) \geq 0$ \quad for all $P$, $Q$ under investigation \hspace{0.7cm} (nonnegativity),
\item[{\ \ }  (D2)] $D\Big(P,Q\Big) = 0$ \ if and only if \ 
$P=Q$ \quad (reflexivity; 
\ identity of indiscernibles).

\end{enumerate}

\noindent
Being interpreted as ``directed'' distances, the divergences $D(\cdot,\cdot)$
can be connected to geometric issues in various different ways.
For the sake of brevity, we mention here only a few of those.

\vspace{0.3cm}
\noindent
To start with an ``all-encompassing view'', 
following the lines of e.g. Birkhoff~\cite{Bir:32}
and Millmann \& Parker~\cite{Mil:91},
one can build 
from 
any set $\mathcal{S}$, 
whose elements can be interpreted as ``points'',
together with a collection $\mathcal{L}$ of non-empty subsets
of $\mathcal{S}$, interpreted as ``lines'' (as a manifestation of a 
principle sort
of structural connectivity
between points), 
and an arbitrary \textit{symmetric distance} 
$\mathfrak{d}(\cdot,\cdot)$
on $\mathcal{S} \times \mathcal{S}$,
an axiomatic constructive framework of geometry which can be
of far-reaching nature; therein, $\mathfrak{d}(\cdot,\cdot)$ plays basically
the role of a marked ruler. Accordingly, each triplet $(\mathcal{S}, \mathcal{L},
\mathfrak{d}(\cdot,\cdot))$ forms a distinct ``quantitative geometric system'';
the most prominent classical case is certainly 
$\mathcal{S} = \mathbb{R}^2$ with $\mathcal{L}$
as the collection of all vertical and non-vertical lines, equipped 
with the Euclidean distance $\mathfrak{d}(\cdot,\cdot)$, 
hence generating the usual Euclidean geometry in the two-dimensional space.
In the case that $\mathfrak{d}(\cdot,\cdot)$ is only an \textit{asymmetric distance (divergence)}
but not a distance anymore, 
we propose that some of the outcoming geometric building blocks
have to be interpreted in a direction-based way
(e.g. the use of $\mathfrak{d}(\cdot,\cdot)$ as a marked directed ruler,
the construction of points of equal divergence from a center viewed as distorted directed spheres,
etc.). For $D(\cdot, \cdot)$ one has to work with $\mathcal{S}$
being a family of real-valued functions on $\mathcal{X}$.

Secondly, from any \textit{symmetric distance} $\mathfrak{d}(\cdot,\cdot)$ on a ``sufficiently rich''
set $\mathcal{S}$  
and a finite number of (fixed or adaptively flexible) distinct ``reference points''
$s_{i}$ ($i=1, \ldots, n$) one can construct the corresponding
Voronoi cells $V(s_i)$ by
\vspace{-0.23cm}
$$
V(s_i) := \{ z \in \mathcal{S} : \ \mathfrak{d}(z,s_i) \leq \mathfrak{d}(z,s_j) 
\ \textrm{for all $j=1, \ldots, n$} \, \} .
$$

\vspace{-0.23cm}
\noindent
This produces a tesselation (tiling) of $\mathcal{S}$ which is very useful
for classification purposes. Of course, the geometric shape of these tesselations
is of fundamental 
importance.
In the case that $\mathfrak{d}(\cdot,\cdot)$ is only an \textit{asymmetric distance (divergence)},
then $V(s_i)$ has to be interpreted
as a directed Voronoi cell and then there is also the ``reversely directed'' alternative
\vspace{-0.23cm}
$$
\widetilde{V}(s_i) := \{ z \in \mathcal{S} : \ \mathfrak{d}(s_i,z) \leq \mathfrak{d}(s_j,z) 
\ \textrm{for all $j=1, \ldots, n$} \, \} .
$$ 

\vspace{-0.23cm}
\noindent  
Recent applications where $\mathcal{S} \subset \mathbb{R}^d$ and $\mathfrak{d}(\cdot,\cdot)$ is a Bregman divergence 
or a more general conformal divergence,
can be found e.g. in Boissonnat et. al~\cite{Boi:10},
Nock et al.~\cite{Noc:16b} (and the 
references therein),
where they also deal with the corresponding adaption of
k-nearest neighbour classification methods. 

Moreover, with each (say) asymmetric distance (divergence) $\mathfrak{d}(\cdot,\cdot)$ one can associate a 
\textit{divergence-ball} 
$\mathbb{B}_{\mathfrak{d}}(s,\rho)$ with ``center'' 
$s \in \mathcal{S}$ and ``radius'' $\rho \in ]0,\infty[$,
defined by 
$\mathbb{B}_{\mathfrak{d}}(s,\rho) := \{ 
s \in \mathcal{S} : \ \mathfrak{d}(s,z) \leq \rho\}$, whereas
the corresponding 
\textit{divergence-sphere}
is given by
$\mathbb{S}_{\mathfrak{d}}(s,\rho) := \{ 
s \in \mathcal{S} : \ \mathfrak{d}(s,z) = \rho\}$; \  
see e.g. Csiszar \& Breuer~\cite{Csi:16} for a use of some 
divergence balls 
as a constraint in financial-risk related decisions.
Of course, the ``geometry/topology'' induced by divergence balls
and spheres is generally quite non-obvious; see for instance
Roensch \& Stummer~\cite{Roe:17}, who describe and visualize different 
effects in a 3D-setup of scaled Bregman divergences (which will be covered below).
Moreover, the generalization of
\begin{equation}
D(\Omega, P_{N}^{emp} ) := \inf_{Q \in \Omega} D(Q, P_{N}^{emp} ),
\qquad
\widehat{Q} := \arg \inf_{Q \in \Omega} D(Q, P_{N}^{emp} )
\nonumber
\end{equation}
of the above-mentioned statistical minimum divergence/distance 
estimation problems \eqref{MDE2}, \eqref{MDE3}
can e.g. be (loosely) achieved by blowing up the divergence sphere
$\mathbb{S}_{D}(s,P_{N}^{emp})$ through increasing the radius $\rho$
until it first touches the model $\Omega$. Accordingly,
there may be an interesting interplay between the geometric/topological properties
of both $\mathbb{S}_{D}(s,P_{N}^{emp})$ and the (e.g. non-convex, respectively non-smooth,
respectively non-intersection-of-hyperplanes type, respectively complicated-manifold-type) 
boundary $\partial\Omega$ of $\Omega$
(see e.g. Roensch \& Stummer \cite{Roe:17}).

\vspace{0.1cm}
\noindent
Thirdly, consider a framework where 
$P:= \widetilde{P}_{\theta_{1}}$
and
$Q:= \widetilde{P}_{\theta_{2}}$
depend on some parameters $\theta_1 \in \Theta$,
 $\theta_2 \in \Theta$.
The way of dependence of the function (say) $S_{\cdot}(\widetilde{P}_{\theta})$ 
on the underlying parameter $\theta$ from an appropriate space $\Theta$ of e.g. manifold type,
may show up  
directly e.g. via its operation/functioning as a relevant system-indicator,
or it may be manifested implicitly e.g. such that $S_{\cdot}(\widetilde{P}_{\theta})$ is the
solution of an optimization problem with $\theta$-involving constraints.
In such a framework,
one can induce divergences $D\big(S(\widetilde{P}_{\theta_{1}}),S(\widetilde{P}_{\theta_{2}})\big) 
=: f(\theta_{1},\theta_{2})$
and -- under sufficiently smooth dependence --
study their corresponding differential-geometric behaviour of $f(\cdot,\cdot)$ on 
$\Theta$. 
An example is provided by the Kullback-Leibler divergence
between two distributions of the same exponential family of distributions,
which defines a Bregman divergence on the parameter space.
This and related issues are subsumed in the research field of 
``information geometry''; for comprehensive overviews
see e.g. Amari~\cite{Ama:00}, Amari~\cite{Ama:16}, Ay et al.~\cite{Ay:17}.
Moreover, for recent connections between divergence-based information geometry and optimal transport
the reader is e.g. referred to Pal \& Wong~\cite{Pal:16,Pal:18}, Karakida \& Amari~\cite{Kar:17}, Amari et al.~\cite{Ama:18},
Peyre \& Cuturi~\cite{Pey:19}, and the 
literature
therein.

\vspace{0.1cm}
\noindent
Further relations of divergences with other approaches to geometry 
can be over\-viewed e.g. from
the wide-range-covering research-article collections in 
Nielsen \& Bhatia~\cite{Nie:13b}, 
Nielsen \& Barbaresco~\cite{Nie:13a},~\cite{Nie:15a},~\cite{Nie:17a},~\cite{Nie:19a},~\cite{Nie:21b},~\cite{Barb:21} 
and Nielsen~\cite{Nie:21a}.

\vspace{0.1cm}
\noindent
Moreover, geometry also enters as a tool for visualizing
quantitative effects on divergences. A more detailed discussion 
(including also other approaches) on the interplay between
statistics and geometry is beyond the scope of this chapter;
they will appear in other parts of this book.

\subsection{Some Incentives for Extensions}

\subsubsection{$\phi-$Divergences Between Other Statistical Objects:}

Recall that for probability distributions $P$ and $Q$ 
having strictly positive density functions 
$f_P$ and $f_Q$ 
with respect to some measure $\lambda$ on a data space $\mathcal{X}$
(which covers as special cases both the classical density functions respectively
the probability mass functions),
we have defined the $\phi-$divergences 
(Csiszar-Ali-Silvey-Morimoto (CASM) divergences) by
\begin{eqnarray} 
& & \hspace{-0.2cm} 
0 \leq D_{\phi}(P,Q) := \int_{{\mathcal{X}}} 
f_{Q}(x) \cdot \phi \negthinspace \left( 
{\frac{f_{P}(x)}{f_{Q}(x)}}\right) 
\, \mathrm{d}\lambda(x) 
=: D_{\phi,\lambda}(f_{P},f_{Q})
\, , 
\label{BroStuHB22:fo.def.47d.reduced4} 
\end{eqnarray}
where the last notation-type term in \eqref{BroStuHB22:fo.def.47d.reduced4}
indicates the interpretation as $\phi-$divergence between
\textit{density functions}, measuring their similarity.
However, e.g. for $\mathcal{X} \subset \mathbb{R}$ and the Lebesgue measure $\lambda=\lambda_{L}$
(and hence almost always $\mathrm{d}\lambda_{L}(x) = \mathrm{d}x$),
it makes also sense to quantify the dissimilarity --- in terms
of $\phi-$divergences --- between other related ``statistical objects'',
most notably between the information-aggregating cumulative distribution functions $F_{P}$ and $F_{Q}$
of $P$ and $Q$. For instance, formally,  
$D_{\phi_{PC},Q}(F_{P},F_{Q}) = 
\frac{1}{2} \int_{{\mathcal{X}}} 
\frac{(F_{P}(x)-F_{Q}(x))^2}{F_{Q}(x)} 
\, \mathrm{d}Q(x), 
$
(cf. \eqref{BroStuHB22:fo.PC.reduced} with $f_{P}$,$f_{Q}$ replaced by $F_{P}$,$F_{Q}$
and $\lambda =Q$)
is --- in case of employing the empirical measure $P=P_{N}^{emp}$ ---
a special member of the family of weighted Cramer-von Mises test statistics
(in fact it is a \textit{modified Anderson-Darling test statistics}
of e.g. Ahmad et al. \cite{Ahm:88} and Scott \cite{Sco:99}, see also Shin et al. \cite{Shin:12}
for applications in environmental extreme-value theory).  

\vspace{0.2cm}
\noindent
As another incentive, let us mention the use of $\phi-$divergences 
between quantile functions respectively between ``transformations'' thereof. 
For instance, they can be employed in situations
where the above-mentioned
classical minimum $\phi-$divergence/distance estimation problem 
\eqref{MDE2} and \eqref{MDE3} --- which involves $\phi-$divergences between density functions ---
is theoretically and practically intractable; this is e.g. the case when the model $\Omega$ is defined
by constraints on the expectation of a $L-$statistics (e.g. describing a tubular
neighborhood of a distribution with prescribed number of given quantiles;
such constraints are not linear with respect to the underlying distribution of
the data, but merely with respect to their quantile measure).
In such a situation, one can transpose everything to a minimization problem
for the $\phi-$divergence between the corresponding empirical quantile measures 
where the constraint can also be stated in terms of quantile measures
(see Broniatowski \& Decurninge \cite{Bro:16}).

\vspace{0.2cm}
\noindent
Further examples of $\phi-$divergences between other statistical objects can be found in 
Subsection 2.5.1.2 below.

\subsubsection{Some Non-$\phi-$Divergences between Probability Distributions:}

\vspace{0.1cm}
In contrary to the preceding subsection, 
instead of replacing the probability distributions $P$ and $Q$, let us keep the latter two
but consider now some \textit{other} divergences $D(P,Q)$ 
(of non$-\phi-$divergence type) of statistical interest.
For instance, there is a substantially growing amount of applications of the so-called
\textit{(ordinary/classical) Bregman distances/divergences OBD} 
\begin{eqnarray} 
& & \hspace{-0.9cm} \textstyle
0 \leq D_{\phi}^{OBD}(P,Q) 
\nonumber\\ 
& & \hspace{-0.9cm} 
: =
\int_{{\mathcal{X}}} 
\bigg[ \phi \negthinspace \left( f_{P}(x)\right) -\phi \negthinspace \left( f_{Q}(x) \right)
- \phi^{\prime} \negthinspace
\left( f_{Q}(x) \right) \cdot \left( f_{P}(x) - f_{Q}(x) \right) 
\bigg] 
\, \mathrm{d}\lambda(x) \ , 
\label{BroStuHB22:fo.OBD0a}
\end{eqnarray}
(see e.g. Csiszar~\cite{Csi:91}, Pardo \& Vajda~\cite{Par:97},\cite{Par:03},
Stummer \& Vajda~\cite{Stu:12})
where $\phi^{\prime}$ is the derivative of the supposedly differentiable $\phi$. 
The class \eqref{BroStuHB22:fo.OBD0a} includes as important special cases e.g. 
the \textit{density power divergences}
(also known as Basu-Harris-Hjort-Jones distances, 
cf.~\cite{Bas:98}) with the squared $L_2-$norm as a subcase.
The principal types of statistical applications of OBD are basically the 
same as for the $\phi-$divergences (minimum divergence estimation, robustness etc.);
however, the corresponding technical details may differ substantially.

\vspace{0.2cm}
\noindent
Concerning some recent progress of divergences,
Stummer~\cite{Stu:07} 
as well as Stummer \& Vajda~\cite{Stu:12} 
introduced the concept of \textit{scaled Bregman divergences/distances} SBD

\begin{eqnarray} 
& & \hspace{-0.2cm} \textstyle
0 \leq D_{\phi}^{SBD}\left(P,Q\right ) := D_{\phi,\lambda,m}^{SBD}\left(P,Q\right )
\nonumber\\ 
& & \hspace{-0.2cm} 
:= \int_{{\mathcal{X}}} 
\Bigg[ \phi \negthinspace \left( {\frac{f_{P}(x)}{m(x)}}\right) 
-\phi \negthinspace \left( {\frac{f_{Q}(x)}{m(x)}}\right) 
- \phi^{\prime} \negthinspace
\left( {\frac{f_{Q}(x)}{m(x)}}\right) \cdot \left( \frac{f_{P}(x)}{m(x)}-\frac{f_{Q}(x)}{m(x)}\right) 
\Bigg] 
m(x)  \, \mathrm{d}\lambda(x) 
\nonumber 
\end{eqnarray}
which (by using a scaling function $m(\cdot)$) generalizes all the above-mentioned (nearly disjoint) 
density-based
$\phi-$divergences \eqref{BroStuHB22:fo.def.47d.reduced4}
and OBD divergences \eqref{BroStuHB22:fo.OBD0a} at once. 
Hence, the SBD divergence class constitutes
a quite general framework for dealing with a wide range of data analyses, 
in a well-structured way.

\subsubsection{Some Non-CASM Divergences between Other Statistical Objects:}

\vspace{0.1cm}
Of course, for statistical applications it also makes sense 
to the combine the extension-ideas  ,
of the two preceding subsections. For instance,
\begin{eqnarray}
& & \hspace{-0.2cm}   \textstyle 
0 \leq D_{\phi_{PC},Q,m}^{SBD}\left(F_{P},F_{Q}\right )
= \frac{1}{2} \cdot \int_{{\mathbb{R}}} 
\frac{( F_{P}(x) - F_{Q}(x) )^2}{m(x)}
\, \mathrm{d}Q(x) \  
\nonumber
\end{eqnarray} 
constitutes --- in case of employing the empirical measure $P=P_{N}^{emp}$ ---
the family of weighted Cramer-von Mises test statistics
(see \cite{Cra:28}, \cite{vMi:31}, as well as Smirnov \cite{Smi:36}).

\vspace{0.3cm}
\noindent
In the following, for the rest of the paper we work out an extensive
toolkit of divergences between statistical objects, which goes
far beyond the above-mentioned concepts.

%
\section{The Framework}
\label{sec.2new}
%


\subsection{Statistical functionals $S$ and their dissimilarity}
\label{subsec.2.0averynew}


\noindent
Let us assume that the modeled respectively observed 
random data take
values in a state space 
$\mathcal{Y}$ 
(with at least two distinct values), 
which is equipped with
a system  $\mathcal{
A}$ of admissible events 
($\sigma-$algebra). 
On this, we consider two probability distributions (probability measures) 
$P$ and $Q$ of interest.
By appropriate choices of $(\mathcal{Y},\mathcal{A})$,
such a general context also covers modeling of
series of observations, functional data as well as 
stochastic process data (the latter by choosing $\mathcal{Y}$ as an appropriate
space of paths, i.e., whole scenarios along a set of times).

\vspace{0.3cm}
\noindent
In this paper, we deal with situations where -- 
e.g. in face of the 
dichotomous uncertainty
$P$ versus $Q$ -- the \textbf{statistical decision (inference) goal} 
can be e.g. expressed by means of 
\textbf{``dissimilarity-expressing relations''} 
$\mathfrak{R}
\Big(S(P),S(Q)\Big) $
between univariate
real-valued ``statistical functionals''
 $S(\cdot)$ of the form
$S(P) := \left\{S_{x}(P)\right\}_{x \in \mathcal{X}}$
and $S(Q) := \left\{S_{x}(Q)\right\}_{x \in \mathcal{X}}$
for the two distributions $P$ and $Q$
\footnote{
the statistical functional $S(\cdot)$ can also be thought of as
a function-valued ``plug-in statistics'' respectively 
as a real-valued function on $\mathcal{X}$ which carries a
probability-distribution-valued parameter $\mathbf{\cdot}$ \ ;
accordingly $S(P)$, $S(Q)$ are two different functions 
corresponding to the two different parameter constellations $P$,$Q$;
accordingly $S_{x}(P)$, $S_{x}(Q)$ are the corresponding function values
at $x \in \mathcal{X}$
},
where $\mathcal{X}$ is a set of (at least two different) ``functional indices''. 
As corresponding preliminaries, in this section we broadly discuss 
examples of statistical functionals which we shall employ later on
to recover known --- respectively create new --- divergences between them.

\vspace{0.3cm}
\noindent
In principal, one can distinguish between \textit{unit-free} (e.g.\ ``percentage-type'') functionals
$S(\cdot)$
and \textit{unit-dependent} (e.g.\ monetary) functionals 
$S(\cdot)$.
For the real line $\mathcal{Y}=\mathcal{X}=\mathbb{R}$, the most prominent examples for the former are 
the cumulative distribution functions (cdf) 
$\left\{S_{x}(P)\right\}_{x \in \mathbb{R}} := \left\{F_{P}(x)\right\}_{x \in \mathbb{R}}
:= \left\{P[(-\infty,x]]\right\}_{x \in \mathbb{R}} =: S^{cd}(P)$,
the survival functions (suf)
$\left\{S_{x}(P)\right\}_{x \in \mathbb{R}} := \left\{1-F_{P}(x)\right\}_{x \in \mathbb{R}}
:= \left\{P[(x,\infty)]\right\}_{x \in \mathbb{R}} =: S^{su}(P)$
(which are also called reliability functions or tail functions),  
the ``classical'' probability density functions (pdf)
$\left\{S_{x}(P)\right\}_{x \in \mathbb{R}} := \left\{f_{P}(x)\right\}_{x \in \mathbb{R}}
:= \left\{\frac{\mathrm{d} F_{P}(x)}{\mathrm{d} x}\right\}_{x \in \mathbb{R}}  =: S^{pd}(P)$,
the moment generating functions (mgf)
$\left\{S_{x}(P)\right\}_{x \in \mathbb{R}} := \left\{M_{P}(x)\right\}_{x \in \mathbb{R}}
:= \left\{\int_{\mathcal{Y}} e^{x y} \, \mathrm{d}P(y) \right\}_{x \in \mathbb{R}}  =: S^{mg}(P)$,
and for finite/countable $\mathcal{Y}=\mathcal{X} \subset \mathbb{R}$   
the probability mass functions (pmf)
$ \left\{S_{x}(P)\right\}_{x \in \mathcal{X}} := \left\{p_{P}(x)\right\}_{x \in \mathcal{X}}
:= \left\{P[\{x\}]\right\}_{x \in \mathcal{X}}
=: S^{pm}(P)$; furthermore, we also cover
the centered rank function (cf. e.g. Serfling~\cite{Ser:10a}, Serfling \& Zuo ~\cite{Ser:10b},
also called ``center-outward distribution function'' in e.g. Hallin~\cite{Hal:17}, Hallin et al.~\cite{Hal:21})
$
\left\{S_{x}(P)\right\}_{x \in \mathbb{R}} 
:= \left\{2 \cdot F_{P}(x)-1 \right\}_{x \in \mathbb{R}}
=: S^{cr,1}(P)
$.

\vspace{0.2cm}
\noindent
Continuing on the real line, in contrast to the above discussion on \textit{unit-free} statistical functionals, 
let us now turn our intention to \textit{unit-dependent} statistical functionals.
For the latter,
in case of $\mathcal{Y}=\mathbb{R}$ and $\mathcal{X}=]0,1[$, the most prominent examples 
are the univariate quantile functions
$$
\left\{S_{x}(P)\right\}_{x \in ]0,1[}  
:= \left\{F_{P}^{\leftarrow}(x)\right\}_{x \in ]0,1[}
:= \left\{\inf\{z \in \mathbb{R}: F_{P}(z) \geq x \} \right\}_{x \in ]0,1[}
 =: S^{qu}(P);
 $$ 
for $\mathcal{Y}=[0,\infty)$ we take
$$\left\{S_{x}(P)\right\}_{x \in ]0,1[}  
:= \left\{F_{P}^{\leftarrow}(x)\right\}_{x \in ]0,1[}
:= \left\{\inf\{z \in [0,\infty): F_{P}(z) \geq x \} \right\}_{x \in ]0,1[}
 =: S^{qu}(P).
$$   
Of course, if the underlying cdf $z \rightarrow F_{P}(z)$ is strictly increasing, then
$x \rightarrow F_{P}^{\leftarrow}(x)$ is nothing but its ``classical'' inverse function. 
Let us also mention that in quantitative finance and insurance, 
the quantile $F_{P}^{\leftarrow}(x)$ (e.g. quoted in US dollars units) is called the value-at-risk for confidence level $x \cdot 100\%$. A detailed discussion on properties and pitfalls of univariate quantile
functions can be found e.g. in Embrechts \& Hofert~\cite{Emb:13};
see also e.g. Gilchrist \cite{Gil:00}
for a comprehensive survey on quantile functions for practitioners of statistical
modelling.\\
\indent
Similarly, the generalized inverse of the 
centered rank function amounts to 
so-called median-oriented quantile function (cf. Serfling~\cite{Ser:06})
$
\left\{S_{x}(P)\right\}_{x \in ]-1,1[} 
:= \left\{(2 \cdot F_{P}(\cdot)-1)^{-1}(x) \right\}_{x \in ]-1,1[}
=\left\{F_{P}^{\leftarrow}\negthinspace\left(\frac{1+x}{2} \right) \right\}_{x \in ]-1,1[}
=: S^{mqu}(P)
$. The sign of $x$ indicates the direction from the median $M_{P} := F_{P}^{\leftarrow}\big(\frac{1}{2}\big)$.

\vspace{0.2cm}
\noindent
If the distribution $P$ is generated by 
some univariate real-valued random variable, say $Y$, then (with a slight abuse of
notation) one has the obvious interpretations 
$F_{P}(x) =
P[Y \leq x]$, $p_{P}(x) = P[Y=x]$ and   
$F_{P}^{\leftarrow}(x) = \inf\{z \in \mathbb{R}: P[Y \leq z] \geq x \}$.

\vspace{0.2cm}
\noindent
Let us mention that for $\mathcal{X}=\mathcal{Y}= \mathbb{R}$
we also cover ``integrated statistical functionals'' of the form
$S(P) := \left\{S_{x}(P)\right\}_{x \in \mathbb{R}} := 
\left\{\int_{-\infty}^{x} \breve{S}_{z}(P)  \, \mathrm{d}\breve{\lambda}(z) \right\}_{x \in \mathbb{R}}
=: S^{\breve{\lambda},\breve{S}}(P)$
where $\breve{\lambda}$ is a $\sigma-$finite measure on $\mathbb{R}$
and $\breve{S}(P) := \left\{\breve{S}_{z}(P)\right\}_{z \in \mathbb{R}}$ is a 
non-negative respectively $\breve{\lambda}-$integrable statistical
functional. For special cases 
$S^{Q,S^{cd}}(P)$ (i.e. $\breve{\lambda} = Q$ and $\breve{S} = S^{cd}$) 
as well as $S^{Q,S^{cd}}(Q)$
in a goodness-of-fit testing context, see e.g. Henze \& Nikitin~\cite{Hen:00}.

\vspace{0.2cm}
\noindent
For the multidimensional Euclidean space $\mathcal{Y}=\mathcal{X}=\mathbb{R}^d$ ($d\in \mathbb{N}$), 
unit-free-type examples 
are the ``classical'' cumulative distribution functions (cdf) 
$\left\{S_{x}(P)\right\}_{x \in \mathbb{R}^d} := \left\{F_{P}(x)\right\}_{x \in \mathbb{R}^d}
:= \left\{P[\, ]-\infty,x] \, ]\right\}_{x \in \mathbb{R}^d} =: S^{cd}(P)$
(which are based on marginal orderings),
the ``classical'' probability density functions (pdf)
$\left\{S_{x}(P)\right\}_{x \in \mathbb{R}^d} := \left\{f_{P}(x)\right\}_{x \in \mathbb{R}^d}
=: S^{pd}(P)$ (such that $P[\cdot] := \int_{\cdot} f_{P}(x) \, \mathrm{d}\lambda_{L}(x)$
with $d-$dimensional Lebesgue measure $\lambda_{L}$),
the moment generating functions (mgf)
$\left\{S_{x}(P)\right\}_{x \in \mathbb{R}^d} := \left\{M_{P}(x)\right\}_{x \in \mathbb{R}^d}
:= \left\{\int_{\mathcal{Y}} e^{<x,y>} \, \mathrm{d}P(y) \right\}_{x \in \mathbb{R}^d}  =: S^{mg}(P)$,
and for finite/countable $\mathcal{Y}=\mathcal{X} \subset \mathbb{R}^d$   
the probability mass functions (pmf)
$ \left\{S_{x}(P)\right\}_{x \in \mathcal{X}} := \left\{p_{P}(x)\right\}_{x \in \mathcal{X}}
:= \left\{P[\{x\}]\right\}_{x \in \mathcal{X}}
=: S^{pm}(P)$. Furthermore, we cover statistical depth functions
$ \left\{S_{x}(P)\right\}_{x \in \mathbb{R}^d} := \left\{\mathfrak{D}_{P}(x)\right\}_{x \in \mathbb{R}^d}
=: S^{de}(P)$
and statistical outlyingness functions
$ \left\{S_{x}(P)\right\}_{x \in \mathbb{R}^d} := \left\{O_{P}(x)\right\}_{x \in \mathbb{R}^d}
=: S^{ou}(P)$
e.g. in the sense of 
Zuo \& Serfling~\cite{Zuo:00a} (see also Chernozhukov et al.~\cite{Che:17}):
basically,
$x \mapsto \mathfrak{D}_{P}(x) \geq 0$
\footnote{there are also version allowing for negative values, not discussed here}
provides a $P-$based center-outward ordering of points $x \in \mathbb{R}^d$
(in other words, it measures how deep (central) a point $x \in \mathbb{R}^d$
is with respect to $P$),
where the point $M_{P}$ of maximal depth (deepest point, if unique) is interpreted as multidimensional median
and the depth decreases monotonically as $x$ moves away from $M$ along any straight line 
running through the deepest point; moreover, $\mathfrak{D}_{P}(\cdot)$ should be affine invariant 
(in particular, independent on the underlying coordinate system) and vanishing at infinity;
in practice, $\mathfrak{D}_{P}(\cdot)$ is typically bounded. In essence,
higher depth values represent greater ``centrality''.
A corresponding outlying function $O_{P}(\cdot)$ is basically 
$O_{P}(\cdot) := f_{OD}(\mathfrak{D}_{P}(\cdot))$ for some strictly 
decreasing (but not necessarily bounded) nonnegative function 
$f_{OD}$ of $\mathfrak{D}_{P}(\cdot)$, such as
$O_{P}(\cdot) := \frac{1}{\mathfrak{D}_{P}(\cdot)} -1$ or 
$O_{P}(\cdot) := c \cdot (1-\frac{\mathfrak{D}_{P}(\cdot)}{\sup_{z \in \mathbb{R}^d} \mathfrak{D}_{P}(z)})$ 
for some constant $c>0$ (in case that $\mathfrak{D}_{P}(\cdot)$ is bounded).
Accordingly, $O_{P}(\cdot)$ provides a $P-$based center-inward ordering of points $x \in \mathbb{R}^d$:
higher values represent greater ``outlyingness''.
Since $f_{OD}$ is invertible, one can always ``switch equivalently'' between 
$\mathfrak{D}_{P}(\cdot)$ and $O_{P}(\cdot)$. Several examples for $\mathfrak{D}_{P}(\cdot)$ respectively $O_{P}(\cdot)$
can be found e.g. in Liu et al.~\cite{LiuR:99}, Zuo \& Serfling~\cite{Zuo:00a,Zuo:00b}, Serfling~\cite{Ser:02}.

\vspace{0.2cm}
According to the ``D-O-Q-R paradigm'' of Serfling~\cite{Ser:10a},
one can link to the \textit{univariate/one-dimensional} 
$P-$characteristics $\mathfrak{D}_{P}(\cdot)$, $O_{P}(\cdot)$ 
two \textit{multivariate/$d-$dimensional} $P-$characteristics,
namely a centered rank function $R_{P}(\cdot)$ 
(also called center-outward distribution function in Hallin~\cite{Hal:17}, Hallin et al.~\cite{Hal:21})
and a quantile function (also called center-outward quantile surface in Liu et al.~\cite{LiuR:99},
and center-outward quantile function in Hallin~\cite{Hal:17}, Hallin et al.~\cite{Hal:21}),
which are inverses of each other. Such a linkage works e.g. basically as follows:
firstly, one chooses some bounded 
set $\mathbb{B} \subset \mathbb{R}^d$ of ``indices'',
often the $d-$dimensional 
unit ball $\mathbb{B}:= \mathbb{B}_{d}(0)$
which we henceforth use for the following explanations.
Secondly, a $P-$based quantile function 
$\mathfrak{Q}_{P}: \mathbb{B}_{d}(0) \mapsto \mathbb{R}^d$
with ``full'' range $\mathcal{R}\left(\mathfrak{Q}_{P}\right) = \mathbb{R}^d$
is such that it generates contour sets (level sets)
 $\mathcal{C}_{c}:= \{\mathfrak{Q}_{P}(u): ||u|| =c\}$, $0 \leq c < 1$
(where $|| \cdot ||$ denotes the Euclidean norm on $\mathbb{R}^d$)
which are \textit{nested} 
(as $c$ varies increasingly). 
The most central point $M_{P} := \mathfrak{Q}_{P}(0)$ is interpreted as $d-$dimensional median. 
The magnitude $c$ represents a degree of outlyingness for all data points 
in $\mathcal{C}_{c}$,
and higher $c-$values corresponding to ``more extreme data points'' 
\footnote{notice that this kind of outlyingness concept is intrinsic (with respect to
$P$), as opposed to the ``relative outlyingness'' defined as a degree of mismatch between
the frequency of certain data-observation points compared to the corresponding (very much lower) modelling frequency;
see e.g. Lindsay~\cite{Lind:94}, Basu et al.~\cite{Bas:11}, 
and the corresponding flexibilization in Ki{\ss}linger \& Stummer~\cite{Kis:16}
}. 
Thirdly, $R_{P}$ 
is taken to be the (possibly multi-valued) inverse of $\mathfrak{Q}_{P}$.
For technical purposes, one attempts to use quantile functions $\mathfrak{Q}_{P}$
such that the contour sets $\mathcal{C}_{c}$ are ``strictly nested''
(in the sense that the do not intersect for different $c$'s)
such that the inverse function $R_{P}: \mathbb{R}^d \mapsto \mathbb{B}_{d}(0)$ 
is determined by uniquely solving the equation
$y=\mathfrak{Q}_{P}(u)$ for $u \in \mathbb{B}_{d}(0)$,
for all $y \in \mathbb{R}^d$.
Finally, as a naturally corresponding outlyingness function one can e.g. take
the magnitude
$O_{P}(y) := || R_{P}(y) ||$ (i.e. the $c$ for which $y \in \mathcal{C}_{c}$)
and derive the associated depth function $\mathfrak{D}_{P}(y) = f_{OD}^{\leftarrow}(O_{P}(y))$.
Since our divergence framework deals with univariate statistical functionals,
we shall work with the $i-$th components  
$ \left\{S_{x}(P)\right\}_{x \in \mathbb{R}^d} := \left\{\mathfrak{Q}_{P}^{(i)}(x)\right\}_{x \in \mathbb{B}}
=: S^{cqu,i}(P)$ and
$ \left\{S_{x}(P)\right\}_{x \in \mathbb{R}^d} := \left\{R_{P}^{(i)}(x)\right\}_{x \in \mathbb{R}^d}
=: S^{cr,i}(P)$
($i \in \{ 1,\ldots, d\}$) and finally aggregate the results by adding up
the correspondingly outcoming $d$ divergences over $i$ 
(see e.g. \eqref{BroStuHB22:fo.centrank} and \eqref{BroStuHB22:fo.centquant} below).

\vspace{0.2cm}
\noindent
There are several ways to build up concrete ``D-O-Q-R'' setups.
A recent one which generates centered $d-$dimensional analogues of the 
univariate quantile-transform mapping
and the reciprocal probability-integral transformation --
and which uses Brenier-McCann techniques connected to 
the Monge-Kantorovich theory of optimal mass transporation  --
is constructed by 
Chernozhukov et al.~\cite{Che:17} and Hallin~\cite{Hal:17}, Hallin et al.~\cite{Hal:21}
(see also Figalli~\cite{Fig:18}, Faugeras \& R\"uschendorf~\cite{Fau:17}):
indeed, for absolutely continuous distributions $P$ on $\mathbb{R}^{d}$ with nonvanishing
(Lebesgue) density functions they define
$R_{P}$ 
as the unique gradient $\nabla \psi$ of a convex function $\psi$ -- mapping $\mathbb{R}^{d}$ to
$\mathbb{B}_{d}(0)$ and -- ``pushing
$P$ forward to'' 
the uniform measure $\mathcal{U}(\mathbb{B}_{d}(0))$ on $\mathbb{B}_{d}(0)$
(i.e., the distribution of $\nabla \psi$ under $P$ is $\mathcal{U}(\mathbb{B}_{d}(0))$);
as corresponding quantile function they take the inverse 
$\mathfrak{Q}_{P} := R_{P}^{\leftarrow}$ of $R_{P}$.
As indicated above, this implies the transformations $Z \sim P$  if and only if $R_{P}(Z) \sim  
\mathcal{U}(\mathbb{B}_{d}(0))$ as well as 
$U \sim \mathcal{U}(\mathbb{B}_{d}(0))$  if and only if $\mathfrak{Q}_{P}(U) \sim P$. 
Depth functions for $P$ can be generated from depth functions $D_{\mathcal{U}(\mathbb{B}_{d}(0))}(\cdot)$
by $\mathfrak{D}_{P}(x) := D_{\mathcal{U}(\mathbb{B}_{d}(0))}(R_{P}(x))$ ($x \in \mathbb{R}^d$).
For $d=1$, one arrives at the univariate
$R_{P}(x) = R_{P}^{(1)}(x) = 2 \cdot F_{P}(x)-1$,
$
\mathfrak{Q}_{P}^{(1)}(x) = F_{P}^{\leftarrow}\negthinspace\left(\frac{1+x}{2}\right) 
$,
and thus there are the consistencies  
$S^{cr,1}(P) := \left\{R_{P}^{(1)}(x)\right\}_{x \in \mathbb{R}}
= S^{cr}(P)$, 
$S^{cqu,1}(P) := \left\{\mathfrak{Q}_{P}^{(1)}(x)\right\}_{x \in [-1,1]}
= S^{cqu}(P)
$.

\vspace{0.2cm}
\noindent
There are also several other different approaches to define multidimensional 
analogues of quantile functions, see e.g.
Serfling~\cite{Ser:02,Ser:10a}, Galichon \& Henry~\cite{Gal:12}, 
Faugeras \& R\"uschendorf~\cite{Fau:17}.
All those multivariate quantile functions are also covered by our divergence toolkit,
componentwise.

\vspace{0.2cm}
\noindent
Let us finally mention that
for general state space $\mathcal{Y}$, 
as unit-free statistical functionals one can also take for instance 
families $\left\{S_{x}(P)\right\}_{x \in \mathcal{X}} := \left\{P[E_x]\right\}_{x \in \mathcal{X}}$
of probabilities of some particularly selected concrete events $E_x \in \mathcal{A}$ 
of purpose-driven interest, where $\mathcal{X}$ is some set of indices.

\vspace{0.2cm}
\noindent
As needed later on, notice that
these statistical functionals 
$S(P) = \left\{S_{x}(P)\right\}_{x \in \mathcal{X}}$
have the following different ranges 
$\mathcal{R}\left(S(P)\right)$:
$\mathcal{R}\left(S^{cd}(P)\right) = \mathcal{R}\left(S^{su}(P)\right) 
= \mathcal{R}\left(S^{pm}(P)\right)
\subset [0,1]$,
$\mathcal{R}\left(S^{pd}(P)\right)  
\subset [0,\infty[$,
$\mathcal{R}\left(S^{mg}(P)\right)  
\subset [0,\infty]$,
$\mathcal{R}\left(S^{qu}(P)\right)  
\subset ]-\infty,\infty[$ (respectively $\mathcal{R}\left(S^{qu}(P)\right)  
\subset[0,\infty[$ for non-negative random variable $Y\geq 0$),
$\mathcal{R}\left(S^{de}(P)\right)  
\subset [0,\infty]$,
$\mathcal{R}\left(S^{ou}(P)\right)  
\subset [0,\infty]$,
$\mathcal{R}\left(S^{cr,i}(P)\right)  
\subset [-1,1]$,
$\mathcal{R}\left(S^{cqu,i}(P)\right)  
\subset ]-\infty,\infty[$ ($i \in \{1,\ldots,d\}$),
and 
$\mathcal{R}\left(S^{\breve{\lambda},\breve{S}}(P)\right)$
depends on the choice of $\breve{\lambda}$ and $\breve{S}$.

\vspace{0.2cm}
\noindent
The above-mentioned \textbf{``dissimilarity-expressing functional relations''} \\ 
$\mathfrak{R} \Big(S(P),S(Q)\Big)$ 
can be typically of (i) \textit{numerical} nature or (ii) \textit{graphical/plotting} nature,
or hybrids thereof.
As far as (i) is concerned, for fixed $x \in \mathcal{X}$ the dissimilarity between
the real-valued $S_{x}(P)$
and $S_{x}(Q)$ can be expressed by (weighted) ratios close to $1$, 
(weighted) differences close to $0$, and combinations thereof;
these informations on ``pointwise'' dissimilarities can then be compressed to a
single real number e.g. by means of aggregation (weighted summation, 
weighted integration, etc.) over $x$
or by taking the maximum respectively minimum value with respect to $x$.
In contrast, for $\mathcal{X}=\mathbb{R}$ one widespread tool for (ii)
is to draw a two-dimensional scatterplot $\Big(S_{x}(P),S_{x}(Q)\Big)_{x \in \mathcal{X}}$ 
and evaluate -- visually by eyeballing or quantitatively -- the dissimilarity in terms of sizes of deviations from the
equality-expressing diagonal $(t,t)$. In the above-mentioned special case of
$S_{x}(P)=F_{P}(x)=P[(-\infty,x]]$, $S_{x}(Q)=F_{Q}(x)=Q[(-\infty,x]]$
this leads to the well-known ``Probability-Probability-Plot'' ($PP-Plot$),
whereas the choice
$S_{x}(P) = F_{P}^{\leftarrow}(x) = \inf\{z \in \mathbb{R}: F_{P}(z) \geq x \}$,
$S_{x}(Q) = F_{Q}^{\leftarrow}(x) = \inf\{z \in \mathbb{R}: F_{Q}(z) \geq x \}$,
amounts to the very frequently used ``Quantile-Quantile-Plot'' ($QQ-$Plot).
Moreover, the choice
$S_{x}(P) = \mathfrak{D}_{P}(x)$ and 
$S_{x}(Q) = \mathfrak{D}_{Q}(x)$ for some $P-$based respectively $Q-$based depth function 
generates the $DD-$Plot in the sense of Liu et al.~\cite{LiuR:99}.


\subsection{The divergences (directed distances) $D$}


\vspace{0.3cm}
\noindent
Let us now specify the details of the 
divergences (directed distances)
$D\Big(S(P),S(Q)\Big)$
which we are going to employ henceforth as dissimilarity measures between the statistical functionals
$S(P) := \left\{S_{x}(P)\right\}_{x \in \mathcal{X}}$
and $S(Q) := \left\{S_{x}(Q)\right\}_{x \in \mathcal{X}}$.
To begin with, we equip the index space $\mathcal{X}$ with a $\sigma-$algebra
$\mathcal{F}$ and a $\sigma-$finite measure $\lambda$
(e.g. a probability measure, the Lebesgue measure, a counting measure, etc.);
furthermore, we 
assume that $x \rightarrow S_{x}(P) \in [-\infty,\infty]$ and 
$x \rightarrow S_{x}(Q) \in [-\infty,\infty]$ are correspondingly measurable functions
which satisfy $S_{x}(P) \in ]-\infty,\infty[$, 
$S_{x}(Q) \in ]-\infty,\infty[$ for $\lambda$-almost all (abbreviated as $\lambda$-a.a.) 
$x \in \mathcal{X}$.
For such a context, we quantify the (aggregated) divergence
$D(S(P),S(Q)) := D^{c}_{\beta}(S(P),S(Q))$
between the two statistical functionals 
$S(P)$ and $S(Q)$
in terms of the ``parameters'' $\beta = (\phi,m_{1},m_{2},m_{3},\lambda)$ and $c$ by

\vspace{-0.3cm}

\begin{eqnarray} 
& & \hspace{-0.2cm} \textstyle
0 \leq D^{c}_{\phi,m_{1},m_{2},m_{3},\lambda}(S(P),S(Q)) 
\nonumber\\ 
& & \hspace{-0.2cm} 
: =
\olint_{{\mathcal{X}}} 
\Bigg[ \phi \negthinspace \left( {\frac{S_{x}(P)}{m_{1}(x)}}\right) 
-\phi \negthinspace \left( {\frac{S_{x}(Q)}{m_{2}(x)}}\right)
- \phi_{+,c}^{\prime} \negthinspace
\left( {\frac{S_{x}(Q)}{m_{2}(x)}}\right) \cdot \left( \frac{S_{x}(P)}{m_{1}(x)}-\frac{S_{x}(Q)}{m_{2}(x)}\right) 
\Bigg] 
m_{3}(x) \, \mathrm{d}\lambda(x),
\nonumber \\ 
\label{BroStuHB22:fo.def.1}
\end{eqnarray}
where the meaning of the integral symbol $\olint$ -- as a shortcut of the 
integral over an appropriate extension of the integrand -- 
will become clear in \eqref{BroStuHB22:fo.def.3} below.
Here, in accordance with the \textit{BS distances} 
of Broniatowski \& Stummer~\cite{Bro:19a} --- who
flexibilized/widened the concept of \textit{scaled Bregman distances} of 
Stummer~\cite{Stu:07} and Stummer \& Vajda~\cite{Stu:12} 
--- we use the following ingredients:

\begin{itemize}

\item[(I1)] 
(measurable) \textit{scaling functions}
$m_{1}: \mathcal{X} \rightarrow [-\infty, \infty]$
and $m_{2}: \mathcal{X} \rightarrow [-\infty, \infty]$
as well as a
nonnegative 
(measurable) \textit{aggregating function}
$m_{3}: \mathcal{X} \rightarrow [0,\infty]$
such that $m_{1}(x) \in ]-\infty, \infty[$, 
$m_{2}(x) \in ]-\infty, \infty[$,
$m_{3}(x) \in [0, \infty[$
for $\lambda-$a.a.
$x \in \mathcal{X}$. 
In analogy with the above notation, we 
use the symbols $m_{i} := \big\{m_{i}(x)\big\}_{x \in \mathcal{X}}$
to refer to the whole functions.
Let us emphasize that we also allow for 
adaptive situations
in the sense that all three functions $m_1(x)$, $m_{2}(x)$, $m_{3}(x)$
(evaluated at $x$)
may also depend on $S_{x}(P)$ and $S_{x}(Q)$, see below.
In the following, $\mathcal{R}\big(G\big)$ denotes the range (image) of 
a function $G := \big\{G(x)\big\}_{x \in \mathcal{X}}$.

\item[(I2)]
the so-called ``divergence-generator'' $\phi$ 
which is a continuous,
convex  (finite) function $\phi: E \rightarrow ]-\infty,\infty[$
on some appropriately chosen open interval $E = ]a,b[$
such that $[a,b]$ 
covers (at least) the union 
$\mathcal{R}\left(\frac{S(P)}{m_{1}}\right) \cup \mathcal{R}\left(\frac{S(Q)}{m_{2}}\right)$ of both ranges 
$\mathcal{R}\left(\frac{S(P)}{m_{1}}\right)$ of 
$\left\{\frac{S_{x}(P)}{m_{1}(x)}\right\}_{x \in \mathcal{X}}$
and $\mathcal{R}\left(\frac{S(Q)}{m_{2}}\right)$ of $\left\{\frac{S_{x}(Q)}{m_{2}(x)}\right\}_{x \in \mathcal{X}}$;
for instance, $E=]0,1[$, $E=]0,\infty[$ or $E=]-\infty,\infty[$;
the class of all such functions will be
denoted by $\Phi(]a,b[)$.
Furthermore, we assume that $\phi$ is continuously extended 
to $\overline{\phi}: [a,b] \rightarrow [-\infty,\infty]$ by setting 
$\overline{\phi}(t) := \phi(t)$ for $t\in ]a,b[$ as well as  
$\overline{\phi}(a):= \lim_{t\downarrow a} \phi(t)$,
$\overline{\phi}(b):= \lim_{t\uparrow b} \phi(t)$
on the two boundary points $t=a$ and $t=b$. The latter two 
are the the only points at which infinite values may appear
(e.g. because of division by $m_{1}(x)=0$ for some $x$).
Moreover, for any fixed $c \in [0,1]$
the (finite) function 
$\phi_{+,c}^{\prime}: ]a,b[ \rightarrow ]-\infty,\infty[$
is well-defined by
$\phi_{+,c}^{\prime}(t) := c \cdot \phi_{+}^{\prime}(t)
+ (1- c) \cdot \phi_{-}^{\prime}(t)$,  
where $\phi_{+}^{\prime}(t)$ denotes
the (always finite) right-hand derivative of $\phi$ at the point $t \in ]a,b[$
and $\phi_{-}^{\prime}(t)$ the (always finite) left-hand derivative of $\phi$ at $t \in ]a,b[$. 
If $\phi \in \Phi(]a,b[)$ is
also continuously differentiable -- which we denote by $\phi \in \Phi_{C_{1}}(]a,b[)$ --
then for all $c \in [0,1]$
one gets $\phi_{+,c}^{\prime}(t) = \phi^{\prime}(t)$ ($t \in ]a,b[$)
and in such a situation we also suppress $+$ as well as $c$ in all
the following expressions.
We also employ the continuous continuation
$\overline{\phi_{+,c}^{\prime}}: [a,b] \rightarrow [-\infty,\infty]$
given by 
$\overline{\phi_{+,c}^{\prime}}(t) := \phi_{+,c}^{\prime}(t)$ ($t \in ]a,b[$), 
$\overline{\phi_{+,c}^{\prime}}(a) := \lim_{t\downarrow a} \phi_{+,c}^{\prime}(t)$, 
$\overline{\phi_{+,c}^{\prime}}(b) := \lim_{t\uparrow b} \phi_{+,c}^{\prime}(t)$.
To explain the precise meaning of \eqref{BroStuHB22:fo.def.1}, we also make 
use of the (finite, nonnegative) function
$\psi_{\phi,c}: ]a,b[ \times ]a,b[ \rightarrow [0,\infty[$
given by
$\psi_{\phi,c}(s,t) := \phi(s) - \phi(t) - \phi_{+,c}^{\prime}(t) \cdot (s-t) \geq 0$
($s,t \in ]a,b[$). To extend this to a lower semi-continuous function
$\overline{\psi_{\phi,c}}: [a,b] \times [a,b] \rightarrow [0,\infty]$
we proceed as follows:
firstly, we set $\overline{\psi_{\phi,c}}(s,t) := \psi_{\phi,c}(s,t)$ for all $s,t \in ]a,b[$.
Moreover, since for fixed $t \in ]a,b[$, the function $s \rightarrow \psi_{\phi,c}(s,t)$
is convex and continuous, the limit 
$\overline{\psi_{\phi,c}}(a,t) := \lim_{s \rightarrow a} \psi_{\phi,c}(s,t)$
always exists and (in order to avoid 
overlines in \eqref{BroStuHB22:fo.def.1}) will be 
interpreted/abbreviated as $\phi(a) - \phi(t) - \phi_{+,c}^{\prime}(t) \cdot (a-t)$.
Analogously, for fixed $t \in ]a,b[$ we set 
$\overline{\psi_{\phi,c}}(b,t) := \lim_{s \rightarrow b} \psi_{\phi,c}(s,t)$
with corresponding short-hand notation
$\phi(b) - \phi(t) - \phi_{+,c}^{\prime}(t) \cdot (b-t)$.
Furthermore, for fixed $s\in ]a,b[$ we interpret
$\phi(s) - \phi(a) - \phi_{+,c}^{\prime}(a) \cdot (s-a)$
as 
\begin{eqnarray} 
& & \hspace{-0.2cm}
\overline{\psi_{\phi,c}}(s,a) := \left\{
\phi(s) - \overline{\phi_{+,c}^{\prime}}(a) \cdot s
+ \lim_{t \rightarrow a} \Big(t \cdot \overline{\phi_{+,c}^{\prime}}(a) - \phi(t) \Big) 
\right\}
\cdot \boldsymbol{1}_{]-\infty,\infty[}\left(\overline{\phi_{+,c}^{\prime}}(a)\right)
\nonumber\\ 
& & \hspace{1.6cm} 
+ \ \infty \cdot \boldsymbol{1}_{\{-\infty\}}\left(\overline{\phi_{+,c}^{\prime}}(a)\right) \, ,
\nonumber
\end{eqnarray}
where the involved limit always exists but may be infinite.
Analogously, for fixed $s\in ]a,b[$ we interpret
$\phi(s) - \phi(b) - \phi_{+,c}^{\prime}(b) \cdot (s-b)$
as 
\begin{eqnarray} 
& & \hspace{-0.2cm}
\overline{\psi_{\phi,c}}(s,b) := \left\{
\phi(s) - \overline{\phi_{+,c}^{\prime}}(b) \cdot s
+ \lim_{t \rightarrow b} \Big(t \cdot \overline{\phi_{+,c}^{\prime}}(b) - \phi(t) \Big)
\right\}
\cdot \boldsymbol{1}_{]-\infty,\infty[}\left(\overline{\phi_{+,c}^{\prime}}(b)
\right)
\nonumber\\ 
& & \hspace{1.6cm} 
+ \ \infty \cdot \boldsymbol{1}_{\{+\infty\}}\left(\overline{\phi_{+,c}^{\prime}}(b)\right) \, ,
\nonumber
\end{eqnarray}
where again the involved limit always exists but may be infinite.
Finally, we always set $\overline{\psi_{\phi,c}}(a,a):= 0$, $\overline{\psi_{\phi,c}}(b,b):=0$,
and $\overline{\psi_{\phi,c}}(a,b) := \lim_{s \rightarrow a} \overline{\psi_{\phi,c}}(s,b)$,
$\overline{\psi_{\phi,c}}(b,a) := \lim_{s \rightarrow b} \overline{\psi_{\phi,c}}(s,a)$.
Notice that $\overline{\psi_{\phi,c}}$ is lower-semicontinuous but 
not necessarily continuous.
Since ratios are ultimately involved, we also consistently take 
$\overline{\psi_{\phi,c}}\left(\frac{0}{0},\frac{0}{0}\right) := 0$.

\vspace{0.2cm}
\noindent
With (I1) and (I2), we define the \textit{BS divergence (BS distance)} of
\eqref{BroStuHB22:fo.def.1} precisely as

\begin{eqnarray} 
& & \hspace{-0.2cm} 
\textstyle
0 \leq D^{c}_{\phi,m_{1},m_{2},m_{3},\lambda}\big(S(P),S(Q)\big)
= \olint_{{\mathcal{X}}}  
\psi_{\phi,c}\Big(\frac{S_{x}(P)}{m_{1}(x)},
\frac{S_{x}(Q)}{m_{2}(x)}\Big) \cdot
m_{3}(x) \, \mathrm{d}\lambda(x) \qquad \ 
\label{BroStuHB22:fo.def.2}
\\ 
& & \hspace{-0.2cm} 
\textstyle
: =
\int\displaylimits_{{\mathcal{X}}}  
\overline{\psi_{\phi,c}}\Big(\frac{S_{x}(P)}{m_{1}(x)},
\frac{S_{x}(Q)}{m_{2}(x)}\Big) \cdot
m_{3}(x) \, \mathrm{d}\lambda(x) ,
\label{BroStuHB22:fo.def.3}
\end{eqnarray}

\vspace{-0.3cm}

\noindent
but mostly use the less clumsy notation with $\olint$ given in \eqref{BroStuHB22:fo.def.1},
\eqref{BroStuHB22:fo.def.2} henceforth,
as a shortcut for the implicitly involved boundary behaviour. \quad $\square$

\end{itemize}

\vspace{0.2cm}
\noindent
As a side remark let us mention that, we could further generalize
\eqref{BroStuHB22:fo.def.1} by  
adapting a wider divergence (e.g. non-convex generators $\phi$ covering) 
concept of Stummer \& Ki{\ss}linger~\cite{Stu:17a}
who also deal even with 
nonconvex nonconcave divergence generators $\phi$;
for the sake of brevity, this is omitted here. 

\vspace{0.2cm}
\noindent
Notice that by construction 
one has the following important
assertion (cf. Broniatowski \& Stummer~\cite{Bro:19a}):

\begin{theorem}
Let $\phi \in \Phi(]a,b[)$ and $c \in [0,1]$.\\
Then there holds
$D_{\phi,m_{1},m_{2},m_{3}^{c},\lambda}^{c}(S(P),S(Q)) \geq 0$
(i.e. the above-mentioned desired property (D1) is satisfied).\\
Moreover, $D_{\phi,m_{1},m_{2},m_{3}^{c},\lambda}^{c}(S(P),S(Q)) = 0$
if
$\frac{S_{x}(P)}{m_1(x)}=\frac{S_{x}(Q)}{m_2(x)}$
for $\lambda-$almost all $x \in \mathcal{X}$.\\
Depending on the concrete situation,  
$D_{\phi,m_{1},m_{2},m_{3},\lambda}^{c}(S(P),S(Q))$ may take infinite value.
 
\end{theorem}

\noindent
To get a ``sharp identifiability'', i.e. the correspondingly adapted version 
of the above-mentioned desired reflexivity property (D2) in the form of
\begin{eqnarray} 
& & \hspace{-0.2cm}
D^{c}_{\phi,m_{1},m_{2},m_{3},\lambda}(S(P),S(Q)) = 0 \quad
\textrm{if and only if} \quad 
\frac{S_{x}(P)}{m_1(x)}=\frac{S_{x}(Q)}{m_2(x)}
\ \textrm{for $\lambda-$a.a. $x \in \mathcal{X}$} , 
\nonumber \\  
\label{BroStuHB22:fo.reflex1}
\end{eqnarray} 
one needs further requirements on $\phi \in \Phi(]a,b[)$ and $c \in [0,1]$;
for the rest of the paper, we assume the validity of \eqref{BroStuHB22:fo.reflex1}
holds. 

For instance, the latter is satisfied in a setup where 
$m_{3}(x)= w\left(x,\frac{S_{x}(P)}{m_{1}(x)},\frac{S_{x}(Q)}{m_{2}(x)} \right)$
for some (measurable) function $w: \mathcal{X} \times [a,b] \times [a,b] \rightarrow [0,\infty]$,
and the (correspondingly adapted) Assumptions 2 respectively
Assumptions 3 of Broniatowski \& Stummer~\cite{Bro:19a} hold (cf. Theorem 4 respectively
Corollary 1 therein); in particular, 
this means that $\mathcal{R}\big(\frac{S(P)}{m_{1}}\big) 
\cup \mathcal{R}\big(\frac{S(Q)}{m_{2}}\big) \subset [a,b]$
and that for all 
$s \in \mathcal{R}\big(\frac{S(P)}{m_{1}}\big)$
and   
all 
$t \in \mathcal{R}\big(\frac{S(Q)}{m_{2}}\big)$
the following conditions hold:

\vspace{-0.3cm}

\begin{itemize}

\item 
$\phi$ is strictly convex at $t$;

\item
if $\phi$ is differentiable at $t$ and $s \ne t$, then 
$\phi$ is not affine-linear on the interval $[\min(s,t),\max(s,t)]$
(i.e. between $t$ and $s$);

\item
if $\phi$ is not differentiable at 
$t$, $s > t$ and 
$\phi$ is affine linear on $[t,s]$, then we 
exclude $c=1$ for the (``globally/universally chosen'') subderivative
$\phi_{+,c}^{\prime}(\cdot) = c \cdot \phi_{+}^{\prime}(\cdot)
+ (1- c) \cdot \phi_{-}^{\prime}(\cdot)$;

\item
if $\phi$ is not differentiable at 
 $t$, $s < t$ and 
$\phi$ is affine linear on $[s,t]$, then we exclude $c=0$ for $\phi_{+,c}^{\prime}(\cdot)$.
\end{itemize}

\vspace{0.3cm}
\noindent
In the following, we discuss several important 
(classes of) special cases 
of $\beta = (\phi,m_{1},m_{2},m_{3},\lambda)$
in a well-structured way. Let us start with the latter. 


\subsection{The reference measure $\lambda$}
\label{subsec.2.1verynew}


In \eqref{BroStuHB22:fo.def.1}, 
$\lambda$ governs the
\textit{principle} aggregation structure.
For instance, if one chooses $\lambda = \lambda_{L}$ as
the Lebesgue measure on $\mathcal{X} \subset \mathbb{R}$, then the integral
in \eqref{BroStuHB22:fo.def.1} turns out to be of Lebesgue-type and (with some rare exceptions)
consequently of Riemann-type with $\mathrm{d}\lambda(x) = \mathrm{d}x$. 
In contrast, in the \textit{discrete setup} 
where the index set $\mathcal{X} = \mathcal{X}_{\#}$ has 
countably many elements and is equipped with the counting measure
$\lambda := \lambda_{\#} := \sum_{z \in \mathcal{X}_{\#}} \delta_{z}$ 
(where $\delta_{z}$ is Dirac's 
one-point 
distribution
$\delta_{z}[A] :=  \boldsymbol{1}_{A}(z)$, 
and thus $\lambda_{\#}[\{z\}] =1$ for all $z \in \mathcal{X}_{\#}$), 
then \eqref{BroStuHB22:fo.def.1} simplifies to
\begin{eqnarray} 
& & \hspace{-0.2cm} \textstyle
0 \leq D_{\phi,m_{1},m_{2},m_{3},\lambda_{\#}}(S(P),S(Q)) 
\nonumber\\ 
& & \hspace{-0.2cm} 
: =
\olsum_{z \in \mathcal{X}} 
\Bigg[ \phi \negthinspace \left( {\frac{S_{z}(P)}{m_{1}(z)}}\right) 
-\phi \negthinspace \left( {\frac{S_{z}(Q)}{m_{2}(z)}}\right) 
- \phi_{+,c}^{\prime} \negthinspace
\left( {\frac{S_{z}(Q)}{m_{2}(z)}}\right) \cdot \left( \frac{S_{z}(P)}{m_{1}(z)}-\frac{S_{z}(Q)}{m_{2}(z)}\right) 
\Bigg] 
m_{3}(z) \, ,
\nonumber \\ 
\label{BroStuHB22:fo.def.5}
\end{eqnarray}

\noindent
which we interpret as
$\sum_{{z \in \mathcal{X}}} 
\overline{\psi_{\phi,c}}\big(\frac{S_{z}(P)}{m_{1}(z)},\frac{S_{z}(Q)}{m_{2}(z)}\big) \cdot
\mathbbm{m}_{3}(z)$
with the same conventions and limits as in the paragraph right after \eqref{BroStuHB22:fo.def.1}.


\subsection{The divergence generator $\phi$}
\label{subsec.2.2verynew}


We continue with the inspection of interesting special cases
of $\beta = (\phi,m_{1},m_{2},m_{3},\lambda)$
by dealing with the first component.  
For divergence generator $\phi \in \Phi_{C_1}(]a,b[)$
(recall that then we suppress the obsolete $c$ and subderivative index $+$), 
the formula \eqref{BroStuHB22:fo.def.1} turns into
\begin{eqnarray} 
& & \hspace{-0.2cm} \textstyle
0 \leq D_{\phi,m_{1},m_{2},m_{3},\lambda}(S(P),S(Q)) 
\nonumber\\ 
& & \hspace{-0.2cm} 
: =
\olint_{{\mathcal{X}}} 
\Bigg[ \phi \negthinspace \left( {\frac{S_{x}(P)}{m_{1}(x)}}\right) 
-\phi \negthinspace \left( {\frac{S_{x}(Q)}{m_{2}(x)}}\right) 
- \phi^{\prime} \negthinspace
\left( {\frac{S_{x}(Q)}{m_{2}(x)}}\right) \cdot \left( \frac{S_{x}(P)}{m_{1}(x)}-\frac{S_{x}(Q)}{m_{2}(x)}\right) 
\Bigg] 
m_{3}(x) \, \mathrm{d}\lambda(x) \ ,
\nonumber \\ 
\label{BroStuHB22:fo.def.10}
\end{eqnarray}
whereas \eqref{BroStuHB22:fo.def.5} becomes
\begin{eqnarray} 
& & \hspace{-0.2cm} \textstyle
0 \leq D_{\phi,m_{1},m_{2},m_{3},\lambda_{\#}}(S(P),S(Q)) 
\nonumber\\ 
& & \hspace{-0.2cm} 
: =
\olsum_{x \in \mathcal{X}}  
\Bigg[ \phi \negthinspace \left( {
\frac{S_{x}(P)}{m_{1}(x)}}\right) -\phi \negthinspace \left( {\frac{S_{x}(Q)}{m_{2}(x)}}\right)
- \phi^{\prime} \negthinspace
\left( {\frac{S_{x}(Q)}{m_{2}(x)}}\right) \cdot \left( \frac{S_{x}(P)}{m_{1}(x)}-\frac{S_{x}(Q)}{m_{2}(x)}\right) 
\Bigg] 
m_{3}(x) .
\nonumber 
\end{eqnarray}
\noindent
Formally, by defining the integral functional 
$g_{\phi,m_{3},\lambda}(\xi) := \int_{\mathcal{X}} \phi(\xi(x)) \cdot m_{3}(x) \, \mathrm{d}\lambda(x)$
and plugging in e.g. 
$g_{\phi,m_{3},\lambda}\negthinspace \negthinspace \left( {\frac{S(P)}{m_{1}}}\right) 
= \int_{\mathcal{X}} \phi\negthinspace \left( {\frac{S_{x}(P)}{m_{1}(x)}}\right) 
\cdot m_{3}(x) \, \mathrm{d}\lambda(x)$,
the divergence in \eqref{BroStuHB22:fo.def.10} can be interpreted as 
\begin{eqnarray} 
& & \hspace{-0.2cm} \textstyle
0 \leq D_{\phi,m_{1},m_{2},m_{3},\lambda}(S(P),S(Q)) 
\nonumber\\ 
& & \hspace{-0.2cm} 
= g_{\phi,m_{3},\lambda}\negthinspace \negthinspace \left( {\frac{S(P)}{m_{1}}}\right)
- g_{\phi,m_{3},\lambda}\negthinspace \negthinspace \left( {\frac{S(Q)}{m_{2}}}\right)
- g_{\phi,m_{3},\lambda}^{\prime} \negthinspace \negthinspace \left( {\frac{S(Q)}{m_{2}}},
{\frac{S(P)}{m_{1}}} - {\frac{S(Q)}{m_{2}}}\right)
\label{BroStuHB22:fo.def.17}
\end{eqnarray}
where $g_{\phi,m_{3},\lambda}^{\prime} \negthinspace \negthinspace \left( 
\eta, \, \cdot \, \right)$ denotes the corresponding directional derivate
at $\eta = \frac{S(Q)}{m_{2}}$.

\vspace{0.2cm}
\noindent
An important special case is the following: consider the ``nonnegativity-setup''
$$\textrm{ 
(NN0) \quad
$\frac{S_{x}(P)}{m_{1}(x)} \geq 0$
\ and \ $\frac{S_{x}(Q)}{m_{2}(x)}\geq 0$ \ \ for all $x \in \mathcal{X}$;}
$$
for instance, this always holds for nonnegative scaling functions $m_{1}$,  
$m_{2}$, in combination with $S^{cd}$, $S^{pd}$, $S^{pm}$, $S^{su}$, $S^{mg}$,
$S^{de}$, $S^{ou}$, 
and for nonnegative real-valued random variables also with $S^{qu}$.  
Under (NN0), one can take $a=0$, $b=\infty$, i.e. $E=]0,\infty[$, 
and employ the strictly convex power functions

\begin{eqnarray} 
& & \hspace{-0.2cm}
\tilde{\phi}(t): = 
\tilde{\phi}_{\alpha}(t) := \frac{t^\alpha-1}{\alpha(\alpha-1)} \ \in ]-\infty,\infty[ , 
\qquad t \in ]0,\infty[, \ \alpha \in \mathbb{R}\backslash\{0,1\} \ , 
\nonumber
\\ 
& & \hspace{-0.2cm} 
\phi(t): = \phi_{\alpha}(t) := \tilde{\phi}_{\alpha}(t) - \tilde{\phi}_{\alpha}^{\prime}(1) \cdot (t-1) =
\frac{t^\alpha-1}{\alpha(\alpha-1)}-\frac{t-1}{\alpha-1} \ \in [0,\infty[ , 
\quad t \in ]0,\infty[, 
\nonumber \\
& & \hspace{9.2cm}
\ \alpha \in \mathbb{R}\backslash\{0,1\} \ ,
\label{BroStuHB22:fo.def.32}
\end{eqnarray}

\vspace{0.2cm}
\noindent
The perhaps most important special case is
 $\alpha=2$, for which \eqref{BroStuHB22:fo.def.32} turns into
\begin{eqnarray} 
& & \hspace{-0.2cm} 
\phi_{2}(t) := \frac{(t-1)^2}{2},
\quad t \in ]0,\infty[ = E.  
\label{BroStuHB22:fo.def.34a}
\end{eqnarray}

\noindent
Also notice that the divergence-generator $\phi_{2}$ of \eqref{BroStuHB22:fo.def.34a} 
can be trivially extended to
\begin{eqnarray} 
& & \hspace{-0.2cm} 
\bar{\phi}_{2}(t) := \frac{(t-1)^2}{2},
\quad t \in ]-\infty,\infty[ = \bar{E},  
\label{BroStuHB22:fo.def.34b}
\end{eqnarray}
which is useful in the general setup 
$$\textrm{ 
(GS) \quad
$\frac{S_{x}(P)}{m_{1}(x)} \in [-\infty, \infty]$
\ and \ $\frac{S_{x}(Q)}{m_{2}(x)} \in [-\infty, \infty]$ \ \ for all $x \in \mathcal{X}$;}
$$
which appears for nonnegative scaling functions $m_{1}$,  
$m_{2}$ in combination with $S^{qu}$
for real-valued random variables.

\vspace{0.2cm}
Further examples of everywhere strictly convex divergence generators $\phi$ for the 
nonnegativity-setup (NN0) (i.e. $a=0$, $b=\infty$, $E=]0,\infty[$) 
can be obtained by taking the $\alpha-$limits
\begin{eqnarray} 
& & \hspace{-0.2cm}
\tilde{\phi}_{1}(t) := \lim_{\alpha \rightarrow 1} \phi_{\alpha}(t) =
t \cdot \log t \ \in [- e^{-1},\infty[ , 
\qquad t \in ]0,\infty[, 
\label{BroStuHB22:fo.def.120a} 
\\ 
& & \hspace{-0.2cm} 
\phi_{1}(t) := \lim_{\alpha \rightarrow 1} \phi_{\alpha}(t) =
\tilde{\phi}_{1}(t) - \tilde{\phi}_{1}^{\prime}(1) \cdot (t-1) =
t \cdot \log t + 1 - t \ \in [0, \infty[, 
\quad t \in ]0,\infty[,
\nonumber\\  
\label{BroStuHB22:fo.def.120b} 
\\
& & \hspace{-0.2cm}
\tilde{\phi}_{0}(t) := \lim_{\alpha \rightarrow 0} \phi_{\alpha}(t) =
- \log t \ \in ]-\infty,\infty[ , 
\qquad t \in ]0,\infty[, 
\nonumber
\\ 
& & \hspace{-0.2cm} 
\phi_{0}(t) := \lim_{\alpha \rightarrow 0} \phi_{\alpha}(t) =
\tilde{\phi}_{0}(t) - \tilde{\phi}_{0}^{\prime}(1) \cdot (t-1) =
- \log t + t - 1 \ \in [0, \infty[, 
\quad t \in ]0,\infty[.  
\nonumber \\
\label{BroStuHB22:fo.def.120d}
\end{eqnarray}

\noindent
A list of extension-relevant (cf. (I2)) properties of the functions $\phi_{\alpha}$ 
with $\alpha \in \mathbb{R}$
can be found in Broniatowski \& Stummer~\cite{Bro:19a}. The latter also discuss in detail the
important but (in our context) technically delicate divergence generator 
\begin{eqnarray} 
& & \hspace{-0.2cm}
\phi_{TV}(t):= |t-1|
\label{BroStuHB22:fo.def.699}
\end{eqnarray} 
which is non-differentiable at $t=1$; the latter is
also the only point of strict convexity.

As demonstrated in \cite{Bro:19a}, 
$\phi_{TV}$ can --  in our context -- only be potentially applied if
$\frac{S_{x}(Q)}{m_{2}(x)} = 1$ for $\lambda-$a.a. $x \in \mathcal{X}$,
and one 
\textit{generally} has to exclude $c=1$ and $c=0$ for $\phi_{+,c}^{\prime}(\cdot)$
(i.e. we choose $c \in ]0,1[$);
the latter two can be avoided under some non-obvious 
constraints on the statistical 
functionals $S(P)$, $S(Q)$, 
see for instance Subsection 2.5.1.2 below.


\subsection{The scaling and the aggregation functions $m_1$, $m_2$, $m_3$}


In the above two Subsections \ref{subsec.2.1verynew} and \ref{subsec.2.2verynew}, we have 
presented special cases of the first and the last component of the ``divergence
parameter'' 
$\beta = (\phi,m_{1},m_{2},m_{3},\lambda)$, whereas now we 
focus on $m_1$, $m_2$, $m_3$.
To start with, in accordance with \eqref{BroStuHB22:fo.def.1}, 
the aggregation function $m_{3}$ tunes 
the \textit{fine} aggregation details (recall that $\lambda$ governs the
\textit{principle} aggregation structure).
Moreover, the function $m_1(\cdot)$ scales the statistical
functional $S_{\cdot}(P)$ evaluated at $P$ and 
$m_2(\cdot)$ the same statistical functional $S_{\cdot}(Q)$ evaluated at $Q$.
From a modeling perspective, 
these two scaling functions can e.g. 

\vspace{-0.3cm}

\begin{itemize}
 
\item ``purely direct'' in the sense that
$m_{1}(x)$, $m_{2}(x)$ are chosen to directly reflect some dependence
on the index-state $x\in \mathcal{X}$ (independent of the choice of $S$),
or 

\item ``purely adaptive'' in the sense that
$m_{1}(x) = w_{1}(S_{x}(P),S_{x}(Q))$, $m_{2}(x) = w_{2}(S_{x}(P),S_{x}(Q))$
for some appropriate (measurable) ``connector functions'' $w_{1}$, $w_{2}$ on the
product $\mathcal{R}(S(P)) \times \mathcal{R}(S(Q))$ of the ranges of 
$\left\{S_{x}(P)\right\}_{x \in \mathcal{X}}$
and $\left\{S_{x}(Q)\right\}_{x \in \mathcal{X}}$,
or 

\item ``hybrids'' $m_{1}(x) = w_{1}(x,S_{x}(P),S_{x}(Q))$
$m_{2}(x) = w_{2}(x,S_{x}(P),S_{x}(Q))$.

\end{itemize}

\vspace{0.2cm} 
\noindent
In the remainder of Section \ref{sec.2new}, we illuminate several important sub-setups 
of $m_{1}$, $m_{2}$, $m_{3}$, and special cases therein.
As a side effect, this also shows that our framework \eqref{BroStuHB22:fo.def.1} generalizes
considerably all the concrete divergences in the below-mentioned references (even
for the same statistical functional such as e.g. $S= S^{cd}$); for the sake of
brevity, we mention that only at this point, collectively.

\vspace{0.3cm}

%
\noindent
\textbf{2.5.1 \ \
$\mathbf{m_{1}(x) = m_{2}(x) := m(x)}$, $\mathbf{m_{3}(x) = r(x) \cdot m(x)\in [0,\infty]}$ for some
(measurable) function $\mathbf{r: \mathcal{X} \rightarrow \mathbb{R}}$
satisfying
$\mathbf{r(x) \in ]-\infty,0[ \cup ]0,\infty[}$ for $\mathbf{\lambda -}$a.a. $\mathbf{x \in \mathcal{X}}$
}
%

\vspace{0.3cm}

\noindent
In such a sub-setup, the scaling functions are 
strongly coupled with the aggregation function.
In order to avoid ``case-overlapping'' and ``uncontrolled boundary effects'', unless 
otherwise stated we assume here that
the function $r(\cdot)$ does not (explicitly) dependent
on the functions $m(\cdot)$, $S_{\cdot}(P)$ and $S_{\cdot}(Q)$, 
i.e. it is not of the adaptive form $r(\cdot)= h(\cdot, m(\cdot), S_{\cdot}(P), S_{\cdot}(Q))$.
From \eqref{BroStuHB22:fo.def.1} one can derive
\begin{eqnarray} 
& & \hspace{-0.2cm} \textstyle
0 \leq D^{c}_{\phi,m,m,r\cdot m,\lambda}(S(P),S(Q)) 
\nonumber\\ 
& & \hspace{-0.2cm} 
: =
\olint_{{\mathcal{X}}} 
\Bigg[ \phi \negthinspace \left( {\frac{S_{x}(P)}{m(x)}}\right) 
-\phi \negthinspace \left( {\frac{S_{x}(Q)}{m(x)}}\right)
- \phi_{+,c}^{\prime} \negthinspace
\left( {\frac{S_{x}(Q)}{m(x)}}\right) \cdot \left( \frac{S_{x}(P)}{m(x)}-\frac{S_{x}(Q)}{m(x)}\right) 
\Bigg] 
m(x) \cdot r(x) \, \mathrm{d}\lambda(x) \ ,
\nonumber \\ 
\label{BroStuHB22:fo.def.20}
\end{eqnarray}
which for the discrete setup 
$(\mathcal{X},\lambda) = (\mathcal{X}_{\#},\lambda_{\#})$ 
(recall $\lambda_{\#}[\{x\}] =1$ for all $x \in \mathcal{X}_{\#}$) 
simplifies to
\begin{eqnarray}
& & \hspace{-0.2cm} \textstyle
0 \leq D_{\phi,m,m,r\cdot m,\lambda_{\#}}(S(P),S(Q)) 
\nonumber\\ 
& & \hspace{-0.2cm} 
 =
\olsum_{{x \in \mathcal{X}}} 
\Bigg[ \phi \negthinspace \left( {\frac{S_{x}(P)}{m(x)}}\right)
 -\phi \negthinspace \left( {\frac{S_{x}(Q)}{m(x)}}\right)
- \phi_{+,c}^{\prime} \negthinspace
\left( {\frac{S_{x}(Q)}{m(x)}}\right) \cdot \left( \frac{S_{x}(P)}{m(x)}-\frac{S_{x}(Q)}{m(x)}\right)  
\Bigg] 
m(x) \cdot r(x) \ .
\nonumber \\ 
\label{BroStuHB22:fo.def.21}
\end{eqnarray}

\begin{remark}
\label{BroStuHB22:rem.40}
(a) \ In a context of ``$\lambda-$probability-density functions'' with general 
$\mathcal{X}$ and
$P[\cdot] := \int_{\cdot} f_{P}(x) \, \mathrm{d}\lambda(x)$,
$Q[\cdot] := \int_{\cdot} f_{Q}(x) \, \mathrm{d}\lambda(x)$
satisfying $P[\mathcal{X}] = Q[\mathcal{X}] =1$,
one can take the statistical functionals $S_{x}^{\lambda pd}(P) := f_{P}(x) \geq 0$, 
$S_{x}^{\lambda pd}(Q) := f_{Q}(x) \geq 0$;
accordingly,  for $r(x) \equiv 1$ (abbreviated as function $\mathbb{1}$ with constant value 1)
and $M[\cdot] := \int_{\cdot} m(x) \, \mathrm{d}\lambda(x)$ 
the divergence \eqref{BroStuHB22:fo.def.20} 
can be interpreted as
\begin{eqnarray} 
& & \hspace{-0.2cm} \textstyle
0 \leq D^{c}_{\phi,m,m,\mathbb{1} \cdot m,\lambda}\left(S^{\lambda pd}(P),S^{\lambda pd}(Q)\right ) 
\nonumber\\ 
& & \hspace{-0.2cm} 
= \olint_{{\mathcal{X}}} 
\Bigg[ \phi \negthinspace \left( {\frac{f_{P}(x)}{m(x)}}\right) 
-\phi \negthinspace \left( {\frac{f_{Q}(x)}{m(x)}}\right) 
- \phi_{+,c}^{\prime} \negthinspace
\left( {\frac{f_{Q}(x)}{m(x)}}\right) \cdot \left( \frac{f_{P}(x)}{m(x)}-\frac{f_{Q}(x)}{m(x)}\right) 
\Bigg] 
m(x) 
\, \mathrm{d}\lambda(x) 
\nonumber \\ 
& & \hspace{-0.2cm} 
=: B_{\phi }\left( P, Q \,|\,M\right) , 
\label{BroStuHB22:fo.def.20b}
\end{eqnarray}
\footnote{
in a context where $P$ and $Q$ are risk distributions (e.g. $Q$ is a pregiven reference one)
the SBD
$B_{\phi }\left( P, Q \,|\,M\right)$ can be interpreted
as risk excess of $P$ over $Q$ (or vice versa),
in contrast to Faugeras \& R\"uschendorf~\cite{Fau:18} who
use hemimetrics rather than divergences
}
where the
scaled Bregman divergence 
$B_{\phi }\left( P, Q \,|\,M\right)$
has been first defined in Stummer~\cite{Stu:07}, Stummer \& Vajda~\cite{Stu:12}, see also
Kisslinger \& Stummer~\cite{Kis:13},~\cite{Kis:15a},~\cite{Kis:16} for the ``purely adaptive'' case 
$m(x) = w\big(f_{P}(x),f_{Q}(x)\big)$ and indications
on non-probability measures.
Notice that this directly subsumes for $\mathcal{X}=\mathcal{Y}=\mathbb{R}$ the 
``classical density'' functional
$S^{\lambda pd}(\cdot) = S^{pd}(\cdot)$ with the choice $\lambda = \lambda_{L}$
(and the Riemann integration $\mathrm{d}\lambda_{L}(x)=  \mathrm{d}x$),
as well as for the discrete setup  
$\mathcal{Y} = \mathcal{X}= \mathcal{X}_{\#}$
the ``classical probability mass'' functional $S^{\lambda pd}(\cdot) = S^{pm}(\cdot)$ with the choice $\lambda = \lambda_{\#}$
(recall $\lambda_{\#}[\{x\}] =1$ for all $x \in \mathcal{X}_{\#}$); for the latter,
the divergence \eqref{BroStuHB22:fo.def.20b} reads as
\begin{eqnarray} 
& & \hspace{-0.2cm} \textstyle
0 \leq D^{c}_{\phi,m,m,\mathbb{1} \cdot m,\lambda_{\#}}\left(S^{\lambda_{\#} pd}(P),S^{\lambda_{\#} pd}(Q)\right ) 
= D^{c}_{\phi,m,m,\mathbb{1} \cdot m,\lambda_{\#}}\left(S^{pm}(P),S^{pm}(Q)\right )
\nonumber\\ 
& & \hspace{-0.2cm} 
= \olsum_{{x \in \mathcal{X}_{\#}}} 
\Bigg[ \phi \negthinspace \left( {
\frac{p_{P}(x)}{m(x)}}\right) -\phi \negthinspace \left( {\frac{p_{Q}(x)}{m(x)}}\right)
- \phi_{+,c}^{\prime} \negthinspace
\left( {\frac{p_{Q}(x)}{m(x)}}\right) \cdot \left( \frac{p_{P}(x)}{m(x)}-\frac{p_{Q}(x)}{m(x)}\right) 
\Bigg] 
m(x) \nonumber \\ 
& & \hspace{-0.2cm} 
=: B_{\phi }^{\#}\left( P, Q \,|\,M\right) . 
\label{BroStuHB22:fo.def.20c}
\end{eqnarray}

\vspace{0.2cm}
\noindent
For the important special case of the above-mentioned power-function-type generator
$\phi(t) := \phi_{\alpha}(t) = 
\frac{t^\alpha-\alpha \cdot t+ \alpha - 1}{\alpha \cdot (\alpha-1)}$ 
($\alpha \in ]0,\infty[\backslash\{1\}$), 
Roensch \& Stummer \cite{Roe:19a} (see also Ghosh \& Basu \cite{Gho:16b}
for the unscaled special case $m(x)=1$)
employed the corresponding 
scaled Bregman divergences \eqref{BroStuHB22:fo.def.20b}
in order to obtain robust minimum-divergence-type parameter estimates
for the setup of sequences of independent random variables whose distributions are 
non-identical but linked by a common (scalar or multidimensional) parameter;
this is e.g. important in the context of generalized linear models (GLM)
which are omnipresent in statistics, artificial intelligence and machine learning.

Returning to the general framework, 
for the important special case $\alpha=2$ leading to the above-mentioned generator
$\phi_{2}(t) := \frac{(t-1)^2}{2}$, the scaled Bregman divergences
\eqref{BroStuHB22:fo.def.20b} respectively \eqref{BroStuHB22:fo.def.20c} turn into
\begin{eqnarray} 
& & \hspace{-0.2cm} \textstyle
0 \leq 
B_{\phi_{2}}\left( P, Q \,|\,M\right)
= \olint_{{\mathcal{X}}} 
\frac{(f_{P}(x)-f_{P}(x))^{2}}{2 \cdot m(x)} 
\, \mathrm{d}\lambda(x) 
\nonumber
\end{eqnarray}
respectively
\begin{eqnarray} 
& & \hspace{-0.2cm} \textstyle
0 \leq 
B_{\phi_{2}}^{\#}\left( P, Q \,|\,M\right)
= \olsum_{{x \in \mathcal{X}_{\#}}} 
\frac{(p_{P}(x)-p_{P}(x))^{2}}{2 \cdot m(x)} .
\label{BroStuHB22:fo.quSBDbdis}
\end{eqnarray}

\vspace{0.2cm}
\noindent
For instance, in \eqref{BroStuHB22:fo.def.20c} and \eqref{BroStuHB22:fo.quSBDbdis}, 
if $Y$ is a random variable taking values in the
discrete space $\mathcal{X}_{\#}$, then
$p_{Q}(x) = Q[Y=x]$ may be its probability mass function under a hypothetical/candidate law $Q$,
and $p_{P}(x) = 
\frac{1}{N} \cdot \# \{ i \in \{ 1, \ldots, N\}: Y_i =x \} 
=: p_{P_{N}^{emp}}(x)$
is the probability mass function of the 
corresponding data-derived ``empirical distribution'' 
$P:= P_{N}^{emp} := \frac{1}{N} \cdot \sum_{i=1}^{N}   \delta_{Y_{i}}[\cdot]$
of an $N-$size independent and identically distributed (i.i.d.) sample $Y_1, \ldots, Y_N$ of $Y$
which is nothing but the probability distribution reflecting the underlying
(normalized) histogram; moreover, $m(\cdot)$ is a scaling/weighting. \\
In contrast, within a context of clustered multinomial data, we can basically rewrite 
the \textit{parametric extension
of the Brier's consistent estimator} of Alonso-Revenga et al. \cite{Alo:17} as
$c \cdot \sum_{\ell=1}^{L} B_{\phi_{2}}^{\#}\left( P_{N}^{emp, \ell}, P_{N}^{emp} \,|\,
P_{\widehat{\theta}}\right)$ where $P_{N}^{emp,\ell}$ is the empirical distribution
of the $\ell-$th cluster, $P_{N}^{emp} = \frac{1}{L} \sum_{\ell=1}^{L} P_{N}^{emp,\ell}$,
$P_{\widehat{\theta}}$ is a (minimum-divergence-)estimated distribution
from a (log-linear) model class, and $c$ is an appropriately chosen multiplier
(under the assumption of equal cluster sizes, which can be relaxed in a straightforward
manner).
\\
(b) \ In contrast to (a),
for the context $\mathcal{Y}=\mathcal{X}=\mathbb{R}$,
$r(x) \equiv 1$, 
one obtains
in terms of the cumulative distribution functions $S_{x}^{cd}(P) = F_{P}(x)$,
$S_{x}^{cd}(Q) = F_{Q}(x)$ the two non-probability measures
$\mu^{\mathbb{1} \cdot \lambda,cd}[\cdot] : = \int_{\cdot} F_{P}(x) \, \mathrm{d}\lambda(x) \leq \lambda[\cdot]$ 
and $\nu^{\mathbb{1} \cdot \lambda,cd}[\cdot] : = \int_{\cdot}  F_{Q}(x) \, \mathrm{d}\lambda(x)  \leq \lambda[\cdot]$
with -- possibly infinite -- 
total masses 
$\mu^{\mathbb{1} \cdot \lambda,cd}[\mathbb{R}]$,
$\nu^{\mathbb{1} \cdot \lambda,cd}[\mathbb{R}]$.
The latter two are 
finite if $\lambda$ is a probability measure or a finite measure;
for the non-finite Lebesgue measure $\lambda = \lambda_{L}$ and for intervals $[x_{1},x_{2}]$ 
one can interpret $\mu^{\mathbb{1} \cdot \lambda_{L}, cd}[\, [x_{1},x_{2}] \, ]$ 
as the corresponding
area (between $x_{1}$ and $x_{2}$)  under the distribution function
$F_{P}(\cdot)$. Analogously to
\eqref{BroStuHB22:fo.def.20b}, one can interpret
\begin{eqnarray} 
& & \hspace{-0.2cm} \textstyle
0 \leq D^{c}_{\phi,m,m,\mathbb{1} \cdot m,\lambda}\left(S^{cd}(P),S^{cd}(Q)\right ) 
\nonumber\\ 
& & \hspace{-0.2cm} 
= \olint_{{\mathcal{X}}} 
\Bigg[ \phi \negthinspace \left( {\frac{F_{P}(x)}{m(x)}}\right)
 -\phi \negthinspace \left( {\frac{F_{Q}(x)}{m(x)}}\right)
- \phi_{+,c}^{\prime} \negthinspace
\left( {\frac{F_{Q}(x)}{m(x)}}\right) \cdot \left( \frac{F_{P}(x)}{m(x)}-\frac{F_{Q}(x)}{m(x)}\right) 
\Bigg] 
m(x) \, \mathrm{d}\lambda(x) 
\nonumber \\ 
& & \hspace{-0.2cm} 
=: B_{\phi }\left( \mu^{\mathbb{1} \cdot \lambda,cd}, \nu^{\mathbb{1} \cdot \lambda,cd} \,|\,M\right)
\nonumber  
\end{eqnarray}
as scaled Bregman divergence between the non-probability measures
$\mu^{\mathbb{1} \cdot \lambda,cd}$ and $\nu^{\mathbb{1} \cdot \lambda,cd}$.\\
(c) In a context of mortality data analytics (which is essential for the calculation
of insurance premiums, financial reserves, annuities, pension benefits, various benefits
of social insurance programs, etc.), the divergence \eqref{BroStuHB22:fo.def.21}
(with $r(x) =1$) has been employed by Kr{\"o}mer \& Stummer \cite{Kro:19}
in order to achieve a realistic representation of mortality rates by smoothing
and error-correcting of crude rates; there, $\mathcal{X}$ is a set of ages (in years), 
$S_{x}(P)$ is the so-called data-based \textit{crude annual mortality rate by age $x$},
$S_{x}(Q)$ is an --- optimally determinable --- candidate model member 
(out of a parametric or nonparametric model) for the
unknown true \textit{annual mortality rate by age $x$},
and $m(x)$ is an appropriately chosen scaling at $x$. 
\\
This concludes the current Remark \ref{BroStuHB22:rem.40}.

\end{remark}

\vspace{0.3cm}
\noindent
In the following, we illuminate two important special cases of the scaling (and aggregation-part) function $m(\cdot)$,
namely $m(x) := 1$ and $m(x):= S_{x}(Q)$:

\vspace{0.3cm}

%
\noindent
\textbf{2.5.1.1 \ \
$\mathbf{m_{1}(x) = m_{2}(x) := 1}$, $\mathbf{m_{3}(x) = r(x)}$ for some
(measurable) function $\mathbf{r: \mathcal{X} \rightarrow [0,\infty]}$
satisfying
$\mathbf{r(x) \in ]0,\infty[}$ for $\mathbf{\lambda -}$a.a. $\mathbf{x \in \mathcal{X}}$
}
%

\vspace{0.3cm}

\noindent
In this sub-setup, \eqref{BroStuHB22:fo.def.20} becomes
\begin{eqnarray} 
& & \hspace{-0.9cm} \textstyle
0 \leq D^{c}_{\phi,\mathbb{1},\mathbb{1},r\cdot \mathbb{1},\lambda}(S(P),S(Q)) 
\nonumber\\ 
& & \hspace{-0.9cm} 
: =
\olint_{{\mathcal{X}}} 
\bigg[ \phi \negthinspace \left( S_{x}(P)\right) -\phi \negthinspace \left( S_{x}(Q) \right) 
- \phi_{+,c}^{\prime} \negthinspace
\left( S_{x}(Q) \right) \cdot \left( S_{x}(P) - S_{x}(Q) \right) 
\bigg] 
r(x) \, \mathrm{d}\lambda(x) \ ,
\label{BroStuHB22:fo.def.40var}
\end{eqnarray}
which for the discrete setup 
$(\mathcal{X},\lambda) = (\mathcal{X}_{\#},\lambda_{\#})$ 
turns into
\begin{eqnarray} 
& & \hspace{-0.9cm} \textstyle
0 \leq D^{c}_{\phi,\mathbb{1},\mathbb{1},r\cdot \mathbb{1},\lambda_{\#}}(S(P),S(Q)) 
\nonumber\\ 
& & \hspace{-0.9cm} 
: =
\olsum_{{x \in \mathcal{X}}} 
\bigg[ \phi \negthinspace \left( S_{x}(P)\right) -\phi \negthinspace \left( S_{x}(Q) \right)
- \phi_{+,c}^{\prime} \negthinspace
\left( S_{x}(Q) \right) \cdot \left( S_{x}(P) - S_{x}(Q) \right) 
\bigg] 
r(x) \ .
\label{BroStuHB22:fo.def.41}
\end{eqnarray}
For reasons to be clarified below, in case of differentiable generator $\phi$
(and thus $\phi_{+,c}^{\prime} = \phi^{\prime}$ is the classical derivative)
one can interpret \eqref{BroStuHB22:fo.def.40var}
and \eqref{BroStuHB22:fo.def.41} as weighted Bregman distances between the
two statistical functionals $S(P)$ and $S(Q)$.

\vspace{0.2cm}
\noindent
Let us first discuss the important special case 
$\phi = \phi_{\alpha}$ ($\alpha \in \mathbb{R}$,
cf. \eqref{BroStuHB22:fo.def.32}, \eqref{BroStuHB22:fo.def.120b},
\eqref{BroStuHB22:fo.def.120d}, \eqref{BroStuHB22:fo.def.34b})
together with 
$S_{x}(P) \geq 0$,
$S_{x}(Q) \geq 0$ 
-- as it is always the case for $S^{cd}$, $S^{pd}$, $S^{pm}$, $S^{su}$, $S^{mg}$,
$S^{de}$, $S^{ou}$, 
and for nonnegative real-valued random variables also with $S^{qu}$.
By incorporating the above-mentioned extension-relevant (cf. (I2)) properties of $\phi_{\alpha}$ 
(see Broniatowski \& Stummer~\cite{Bro:19a}) into \eqref{BroStuHB22:fo.def.40var},
we end up with

\begin{eqnarray} 
& & \hspace{-0.2cm}   \textstyle \textstyle
0 \leq D_{\phi_{\alpha},\mathbb{1},\mathbb{1},r\cdot \mathbb{1},\lambda}(S(P),S(Q))
\nonumber\\ 
& & \hspace{-0.2cm}   \textstyle 
= \olint_{{\mathcal{X}}} 
\frac{r(x)}{\alpha \cdot (\alpha-1)} \cdot
\Big[ 
\left(S_{x}(P)\right)^{\alpha} + (\alpha-1) \cdot \left(S_{x}(Q)\right)^{\alpha} - \alpha 
\cdot S_{x}(P) \cdot \left(S_{x}(Q)\right)^{\alpha-1} 
\Big]
 \, \mathrm{d}\lambda(x) \ ,
\label{BroStuHB22:fo.def.140aol}
\\ 
& & \hspace{-0.2cm}   \textstyle 
= \int_{{\mathcal{X}}} 
\frac{r(x)}{\alpha \cdot (\alpha-1)} \cdot
\big[ 
\left(S_{x}(P)\right)^{\alpha} + (\alpha-1) \cdot \left(S_{x}(Q)\right)^{\alpha} - \alpha 
\cdot S_{x}(P) \cdot \left(S_{x}(Q)\right)^{\alpha-1} 
\big] \nonumber\\
& & \hspace{7.8cm} 
\cdot \boldsymbol{1}_{]0,\infty[}\big(S_{x}(P) \cdot S_{x}(Q) \big)
 \, \mathrm{d}\lambda(x) 
\nonumber \\ 
& & \hspace{-0.2cm}   \textstyle 
+ \negthinspace
\int_{{\mathcal{X}}} \negthinspace
r(x)
\negthinspace \cdot \negthinspace \big[ 
\frac{\left(S_{x}(P)\right)^{\alpha}}{\alpha \negthinspace \cdot \negthinspace (\alpha-1)} 
\negthinspace \cdot \negthinspace \boldsymbol{1}_{]1,\infty[}(\alpha) 
\negthinspace + \negthinspace \infty \negthinspace \cdot \negthinspace \boldsymbol{1}_{]-\infty,0[ \cup ]0,1[}(\alpha) 
\big] \negthinspace
\cdot \negthinspace \boldsymbol{1}_{]0,\infty[}\big(S_{x}(P) 
\big) \negthinspace
\cdot \negthinspace \boldsymbol{1}_{\{0\}}\big(S_{x}(Q)\big)
\, \mathrm{d}\lambda(x)
\nonumber \\ 
& & \hspace{-0.2cm}   \textstyle 
+ \negthinspace
\int_{{\mathcal{X}}} \negthinspace
r(x) \negthinspace
\cdot \negthinspace \big[ 
\frac{\left(S_{x}(Q)\right)^{\alpha}}{\alpha} \negthinspace \cdot \negthinspace \boldsymbol{1}_{]0,1[ \cup ]1,\infty[}(\alpha)
\negthinspace + \negthinspace \infty \negthinspace \cdot \negthinspace \boldsymbol{1}_{]-\infty,0[}(\alpha)
\big] \negthinspace
\cdot \negthinspace \boldsymbol{1}_{]0,\infty[}\big(S_{x}(Q) 
\big) \negthinspace
\cdot \negthinspace \boldsymbol{1}_{\{0\}}\big(S_{x}(P)\big) 
\, \mathrm{d}\lambda(x) \, ,
\nonumber \\ 
& & \hspace{8.2cm} 
\textrm{ for } \alpha \in \mathbb{R}\backslash\{0,1\},
\label{BroStuHB22:fo.def.140a} 
\\[-0.2cm]
& & \hspace{-0.2cm}   \textstyle 
0 \leq D_{\phi_{1},\mathbb{1},\mathbb{1},r\cdot \mathbb{1},\lambda}(S(P),S(Q)) 
\nonumber\\ 
& & \hspace{-0.2cm}   \textstyle 
= \olint_{{\mathcal{X}}} 
r(x)  \cdot \big[ 
S_{x}(P) \cdot \log\big(\frac{S_{x}(P)}{S_{x}(Q)}\big) + S_{x}(Q) - S_{x}(P) 
\big] \, \mathrm{d}\lambda(x) 
\label{BroStuHB22:fo.def.140bol} \\ 
& & \hspace{-0.2cm}   \textstyle 
= \int_{{\mathcal{X}}} 
r(x)  \cdot \big[ 
S_{x}(P) \cdot \log\big(\frac{S_{x}(P)}{S_{x}(Q)}\big) + S_{x}(Q) - S_{x}(P) 
\big]
\cdot \boldsymbol{1}_{]0,\infty[}\big(S_{x}(P) \cdot S_{x}(Q)\big) \, \mathrm{d}\lambda(x) 
\nonumber \\ 
& & \hspace{-0.2cm}   \textstyle 
+ 
\int_{{\mathcal{X}}} 
r(x)
\cdot \infty
\cdot \boldsymbol{1}_{]0,\infty[}\big(S_{x}(P) 
\big)
\cdot \boldsymbol{1}_{\{0\}}\big(S_{x}(Q)\big)
\, \mathrm{d}\lambda(x)
\nonumber \\ 
& & \hspace{-0.2cm}   \textstyle 
+ 
\int_{{\mathcal{X}}} 
r(x)
\cdot S_{x}(Q)
\cdot \boldsymbol{1}_{]0,\infty[}\big(S_{x}(Q) 
\big)
\cdot \boldsymbol{1}_{\{0\}}\big(S_{x}(P)\big) 
\, \mathrm{d}\lambda(x)
\label{BroStuHB22:fo.def.140b} 
\end{eqnarray}
\begin{eqnarray}
& & \hspace{-0.2cm}   \textstyle 
0 \leq D_{\phi_{0},\mathbb{1},\mathbb{1},r\cdot \mathbb{1},\lambda}(S(P),S(Q)) 
\nonumber\\ 
& & \hspace{-0.2cm}   \textstyle 
= \olint_{{\mathcal{X}}} 
r(x)  \cdot \Big[ 
- \log\big(\frac{S_{x}(P)}{S_{x}(Q)}\big) + \frac{S_{x}(P)}{S_{x}(Q)} - 1 
\Big] \, \mathrm{d}\lambda(x)
\label{BroStuHB22:fo.def.140col} \\ 
& & \hspace{-0.2cm}   \textstyle 
= \int_{{\mathcal{X}}} 
r(x)  \cdot \Big[ 
- \log\big(\frac{S_{x}(P)}{S_{x}(Q)}\big) + \frac{S_{x}(P)}{S_{x}(Q)} - 1 
\Big]
\cdot  \boldsymbol{1}_{]0,\infty[}\big(S_{x}(P) \cdot S_{x}(Q)\big) \, \mathrm{d}\lambda(x)
\nonumber \\ 
& & \hspace{-0.2cm}   \textstyle 
+ 
\int_{{\mathcal{X}}} 
r(x)
\cdot \infty
\cdot \boldsymbol{1}_{]0,\infty[}\big(S_{x}(P) 
\big)
\cdot \boldsymbol{1}_{\{0\}}\big(S_{x}(Q)\big)
\, \mathrm{d}\lambda(x)
\nonumber \\ 
& & \hspace{-0.2cm}   \textstyle 
+ 
\int_{{\mathcal{X}}} 
r(x)
\cdot \infty
\cdot \boldsymbol{1}_{]0,\infty[}\big(S_{x}(Q) 
\big)
\cdot \boldsymbol{1}_{\{0\}}\big(S_{x}(P)\big) 
\, \mathrm{d}\lambda(x) \, ,
\label{BroStuHB22:fo.def.140c} \\ 
& & \hspace{-0.2cm} \textstyle
0 \leq D_{\phi_{2},\mathbb{1},\mathbb{1}, r\cdot \mathbb{1},\lambda}(S(P),S(Q)) 
= \int_{{\mathcal{X}}} 
\frac{r(x)}{2} \cdot
\Big[ S_{x}(P) - S_{x}(Q) \Big]^2
 \, \mathrm{d}\lambda(x) \ ; 
\label{BroStuHB22:fo.def.140d}
\end{eqnarray} 
as a recommendation,
one should avoid $\alpha \leq 0$ whenever $S_{x}(P) =0$ for all $x$ 
in some $A$ with $\lambda[A]>0$,
respectively $\alpha \leq 1$ whenever $S_{x}(Q) =0$ for all $x$ 
in some $\tilde{A}$ with $\lambda[\tilde{A}]>0$.
As far as splitting of the integral e.g. in \eqref{BroStuHB22:fo.def.140b} 
resp. \eqref{BroStuHB22:fo.def.140c}
is concerned, notice that $\int_{{\mathbb{R}}} \left[S_{x}(Q) - S_{x}(P) \right]
\cdot r(x)  \, \mathrm{d}\lambda(x)$ 
resp.
$\int_{{\mathcal{X}}} \Big[  \frac{S_{x}(P)}{S_{x}(Q)} - 1 \Big]
\cdot r(x)  \, \mathrm{d}\lambda(x)$ 
may be finite even in cases where
$\int_{{\mathcal{X}}} S_{x}(P) \cdot r(x)  \, \mathrm{d}\lambda(x) = \infty$
and 
$\int_{{\mathcal{X}}} S_{x}(Q) \cdot r(x)  \, \mathrm{d}\lambda(x) = \infty$
(take e.g. $\mathcal{X} = [0,\infty[$, $\lambda = \lambda_{L}$,  
$r(x) \equiv 1$, and the exponential distribution functions
$S_{x}(P)= F_{P}(x) = 1- \exp(- c_{1} \cdot x)$, $S_{x}(Q) = F_{Q}(x) = 1- \exp(- c_{2} \cdot x)$
with $0 < c_{1} < c_{2}$). Notice that \eqref{BroStuHB22:fo.def.140d} can be used
also in cases where 
$S_{x}(P) \in \mathbb{R}$,
$S_{x}(Q) \in \mathbb{R}$, 
and thus e.g. for $S^{qu}$ for arbitrary real-valued random variables. 
\\ 
As before, for the discrete setup 
$(\mathcal{X},\lambda) = (\mathcal{X}_{\#},\lambda_{\#})$ 
all the terms $\int_{{\mathcal{X}}} \ldots \, \mathrm{d}\lambda(x)$ in 
\eqref{BroStuHB22:fo.def.140a} to \eqref{BroStuHB22:fo.def.140d}
turn into $\sum_{x \in \mathcal{X}} \ldots $\ .\\

\noindent
\textbf{Distribution functions.} For $\mathcal{Y}=\mathcal{X}=\mathbb{R}$, 
$S_{x}(P) = S_{x}^{cd}(P) = F_{P}(x)$,
$S_{x}(Q) = S_{x}^{cd}(Q) = F_{Q}(x)$, let us illuminate the case
$\alpha=2$ of \eqref{BroStuHB22:fo.def.140d}.
For instance, if $Y$ is a real-valued random variable and 
$F_{Q}(x) = Q[Y\leq x]$ is its probability mass function under a hypothetical/candidate law $Q$,
one can take $F_{P}(x) = 
\frac{1}{N} \cdot \# \{ i \in \{ 1, \ldots, N\}: Y_i \leq x \} =: F_{P_{N}^{emp}}(x)$
as the distribution function of the 
corresponding data-derived ``empirical distribution'' 
$P:= P_{N}^{emp} := \frac{1}{N} \cdot \sum_{i=1}^{N}   \delta_{Y_{i}}[\cdot]$
of an $N-$size i.i.d. sample $Y_1, \ldots, Y_N$ of $Y$.
In such a set-up, the choice (say) $\lambda = Q$
in \eqref{BroStuHB22:fo.def.140d} and multiplication with $2N$
lead to the weighted Cramer-von Mises test statistics 
(see \cite{Cra:28}, \cite{vMi:31}, Smirnov \cite{Smi:36}, and also Darling \cite{Dar:57} 
for a historic account)
\begin{eqnarray} 
& & \hspace{-0.2cm} \textstyle
0 \leq 2N \cdot D_{\phi_{2},\mathbb{1},\mathbb{1},r\cdot 
\mathbb{1},Q}(S^{cd}\left(P_{N}^{emp}\right),S^{cd}(Q)) 
= N \cdot \int_{{\mathbb{R}}} 
\Big[ F_{P_{N}^{emp}}(x) - F_{Q}(x) \Big]^2
\cdot r(x) \, \mathrm{d}Q(x) \ 
\nonumber \\
\label{BroStuHB22:fo.def.140e}
\end{eqnarray}
which are special ``quadratic EDF statistics'' in the sense of Stephens~\cite{Ste:86}
(who also uses the term ``Cramer-von Mises family'').
The special case $r(x) \equiv 1$ is nothing but the 
prominent (unweighted) Cramer-von Mises test statistics;
for some recent statistical insights on the latter, 
see e.g. Baringhaus \& Henze~\cite{Bar:17}. 
In contrast, if one chooses the Lebesgue measure
$\lambda = \lambda_{L}$ and $r(x) \equiv 1$ in \eqref{BroStuHB22:fo.def.140d},
then one ends up with the $N-$fold of the ``classical'' squared $L^2-$ distance between
the two distribution functions $F_{P_{N}^{emp}}(\cdot)$ and $F_{Q}(\cdot)$,
i.e. with
\begin{eqnarray} 
& & \hspace{-0.2cm} \textstyle
0 \leq 2N \cdot D_{\phi_{2},\mathbb{1},\mathbb{1},\mathbb{1},\lambda_{L}}(S^{cd}\left(P_{N}^{emp}\right),S^{cd}(Q)) 
= N \cdot \int_{{\mathbb{R}}} 
\Big[ F_{P_{N}^{emp}}(x) - F_{Q}(x) \Big]^2
\, \mathrm{d}\lambda_{L}(x) \ 
\nonumber \\
\label{BroStuHB22:fo.def.140f}
\end{eqnarray}
where one can typically identify $\mathrm{d}\lambda_{L}(x) = \mathrm{d}x$ (Riemann-integral).

\vspace{0.2cm}
\noindent
In a similar fashion, 
for the special case $\mathcal{Y}=\mathcal{X}=\mathbb{R}$, 
and the
above-mentioned integrated statistical functionals (cf. Section \ref{subsec.2.0averynew})
$S_{x}(P) = S_{x}^{Q,S^{cd}}(P) = \int_{-\infty}^{x} F_{P}(z)  \, \mathrm{d}Q(z)$,
$S_{x}(Q) = S_{x}^{Q,S^{cd}}(Q) = \int_{-\infty}^{x} F_{Q}(z)  \, \mathrm{d}Q(z)$,
we get from 
\eqref{BroStuHB22:fo.def.40var} (analogously to \eqref{BroStuHB22:fo.def.140d})
\begin{eqnarray} 
& & \hspace{-0.2cm} \textstyle
0 \leq D_{\phi_{2},\mathbb{1},\mathbb{1},r\cdot \mathbb{1},\lambda}\left(S^{Q,S^{cd}}(P),S^{Q,S^{cd}}(Q)\right) 
= \int_{{\mathbb{R}}} 
\frac{r(x)}{2} \cdot
\Big[ S_{x}^{Q,S^{cd}}(P) - S_{x}^{Q,S^{cd}}(Q) \Big]^2
 \, \mathrm{d}\lambda(x) \ , 
 \nonumber 
\end{eqnarray}
for which the choice $r(x) \equiv 2N$, $P:= P_{N}^{emp}$, $\lambda=Q$
leads to the divergence used 
in a goodness-of-fit testing context by Henze \& Nikitin~\cite{Hen:00}.

\vspace{0.2cm}
\noindent
\textbf{$\lambda-$probability-density functions.} 
If for general $\mathcal{X}$ one takes the special case  
$r(x) \equiv 1$ together with the ``$\lambda-$probability-density functions'' context 
(cf. Remark \ref{BroStuHB22:rem.40}(c))
$S_{x}(P) = S_{x}^{\lambda pd}(P) := f_{P}(x) \geq 0$, 
$S_{x}(Q) = S_{x}^{\lambda pd}(Q) = f_{Q}(x) \geq 0$,
then the divergences \eqref{BroStuHB22:fo.def.40var} and \eqref{BroStuHB22:fo.def.41}
become
\begin{eqnarray} 
& & \hspace{-0.9cm} \textstyle
0 \leq D^{c}_{\phi,\mathbb{1},\mathbb{1},\mathbb{1},\lambda}(S^{\lambda pd}(P),S^{\lambda pd}(Q)) 
\nonumber\\ 
& & \hspace{-0.9cm} 
: =
\olint_{{\mathcal{X}}} 
\bigg[ \phi \negthinspace \left( f_{P}(x)\right) -\phi \negthinspace \left( f_{Q}(x) \right)
- \phi_{+,c}^{\prime} \negthinspace
\left( f_{Q}(x) \right) \cdot \left( f_{P}(x) - f_{Q}(x) \right) 
\bigg] 
\, \mathrm{d}\lambda(x) \ ,
\label{BroStuHB22:fo.OBD1}
\end{eqnarray}
and
\begin{eqnarray} 
& & \hspace{-0.9cm} \textstyle
0 \leq D^{c}_{\phi,\mathbb{1},\mathbb{1},\mathbb{1},\lambda_{\#}}(S^{\lambda pd}(P),S^{\lambda pd}(Q)) 
\nonumber\\ 
& & \hspace{-0.9cm} 
: =
\olsum_{{x \in \mathcal{X}}} 
\bigg[ \phi \negthinspace \left( f_{P}(x)\right) -\phi \negthinspace \left( f_{Q}(x) \right)
- \phi_{+,c}^{\prime} \negthinspace
\left( f_{Q}(x) \right) \cdot \left( f_{P}(x) - f_{Q}(x) \right) 
\bigg] 
\ .
\label{BroStuHB22:fo.OBD2}
\end{eqnarray}
In case of differentiable generator $\phi$
(and thus $\phi_{+,c}^{\prime} = \phi^{\prime}$ is the classical derivative),
the divergences in \eqref{BroStuHB22:fo.OBD1}
and \eqref{BroStuHB22:fo.OBD2} are nothing but the classical Bregman distances between the
two probability distributions $P$ and $Q$
(see e.g. Csiszar~\cite{Csi:91}, Pardo \& Vajda~\cite{Par:97},\cite{Par:03},
Stummer \& Vajda~\cite{Stu:12}).
If one further specializes $\phi = \phi_{\alpha}$, 
the divergences \eqref{BroStuHB22:fo.def.140aol}, \eqref{BroStuHB22:fo.def.140bol}, 
\eqref{BroStuHB22:fo.def.140col}
and \eqref{BroStuHB22:fo.def.140d} become
\begin{eqnarray} 
& & \hspace{-0.2cm} \textstyle
0 \leq D_{\phi_{\alpha},\mathbb{1},\mathbb{1},\mathbb{1},\lambda}(S^{\lambda pd}(P),S^{\lambda pd}(Q)) 
\nonumber\\ 
& & \hspace{-0.2cm} 
= \olint_{{\mathcal{X}}} 
\frac{1}{\alpha \cdot (\alpha-1)} \cdot
\Big[ 
\left(f_{P}(x)\right)^{\alpha} + (\alpha-1) \cdot \left(f_{Q}(x)\right)^{\alpha} - \alpha 
\cdot f_{P}(x) \cdot \left(f_{Q}(x)\right)^{\alpha-1} 
\Big]
 \, \mathrm{d}\lambda(x) \ ,
\nonumber \\ 
& & \hspace{8.2cm} 
\textrm{ for } \alpha \in \mathbb{R}\backslash\{0,1\},
\label{BroStuHB22:fo.def.153a}
\\
& & \hspace{-0.2cm} 
0 \leq D_{\phi_{1},\mathbb{1},\mathbb{1},\mathbb{1},\lambda}(S^{\lambda pd}(P),S^{\lambda pd}(Q)) 
\nonumber\\ 
& & \hspace{-0.2cm} 
= \olint_{{\mathcal{X}}} 
\Big[ 
f_{P}(x) \cdot \log\left(\frac{f_{P}(x)}{f_{Q}(x)}\right) + f_{Q}(x) - f_{P}(x) 
\Big]
\, \mathrm{d}\lambda(x) \ ,
\label{BroStuHB22:fo.def.153b}
\\
& & \hspace{-0.2cm} 
0 \leq D_{\phi_{0},\mathbb{1},\mathbb{1},\mathbb{1},\lambda}(S^{\lambda pd}(P),S^{\lambda pd}(Q)) 
\nonumber\\ 
& & \hspace{-0.2cm} 
= \olint_{{\mathcal{X}}} 
\Big[ 
- \log\left(\frac{f_{P}(x)}{f_{Q}(x)}\right) + \frac{f_{P}(x)}{f_{Q}(x)} - 1 
\Big]
\, \mathrm{d}\lambda(x) \ ,
\label{BroStuHB22:fo.def.153c}\\
& & \hspace{-0.2cm} \textstyle
0 \leq D_{\phi_{2},\mathbb{1},\mathbb{1},\mathbb{1},\lambda}(S^{\lambda pd}(P),S^{\lambda pd}(Q)) 
= \int_{{\mathcal{X}}} 
\frac{1}{2} \cdot
\Big[ f_{P}(x) - f_{Q}(x) \Big]^2
 \, \mathrm{d}\lambda(x) \ .  
\label{BroStuHB22:fo.def.153d}
\end{eqnarray}
Analogously to the paragraph after \eqref{BroStuHB22:fo.def.140d}, 
one can recommend here to exclude $\alpha \leq 0$ whenever $f_{P}(x) =0$ 
for all $x$ in some $A$ with $\lambda[A]>0$,
respectively $\alpha \leq 1$ whenever $f_{Q}(x) =0$ for all $x$ in some 
$\tilde{A}$ with $\lambda[\tilde{A}]>0$.
As far as splitting of the integral e.g. in \eqref{BroStuHB22:fo.def.140b} 
resp. \eqref{BroStuHB22:fo.def.140c}
is concerned, notice that the integral 
$\left(\mu^{1\cdot\lambda,\lambda pd} -  
\nu^{1\cdot\lambda,\lambda pd}\right)[\mathcal{X}] = \int_{{\mathcal{X}}} \left[f_{Q}(x) - f_{P}(x) \right]
\, \mathrm{d}\lambda(x) = 1 - 1 =0$ but
$\int_{{\mathcal{X}}} \Big[  \frac{f_{P}(x)}{f_{Q}(x)} - 1 \Big]
\, \mathrm{d}\lambda(x)$ 
may be infinite (take e.g. $\mathcal{X} = [0,\infty[$, $\lambda = \lambda_{L}$,
and the exponential distribution density functions 
$f_{P}(x):= c_{1} \cdot \exp(- c_{1} \cdot x)$,
$f_{Q}(x):= c_{2} \cdot \exp(- c_{2} \cdot x)$ with $0 \leq c_{1} \leq c_{2}$).
The choice $\alpha >0$ in \eqref{BroStuHB22:fo.def.153a} coincides with the ``order$-\alpha$'' 
density power divergences DPD of Basu et al.~\cite{Bas:98};
for their statistical applications see e.g. Basu et al.~\cite{Bas:15a}, Ghosh \& Basu~\cite{Gho:16a},
~\cite{Gho:16b} and the references therein,
and for general $\alpha \in \mathbb{R}$ see e.g. Stummer \& Vajda~\cite{Stu:12}.

\vspace{0.2cm}
\noindent
The divergence
\eqref{BroStuHB22:fo.def.153b}
is the celebrated
``Kullback-Leibler information divergence KL'' 
between $f_{P}$ and $f_{Q}$ (respectively between $P$ and $Q$);
alternatively, instead of KL one often uses the terminology ``relative entropy''.
The divergence \eqref{BroStuHB22:fo.def.153d} (cf. $\alpha=2$)
is nothing but half of the squared $L^2-$ distance between
the two $\lambda-$density functions $f_{P}(\cdot)$ and $f_{Q}(\cdot)$.

\vspace{0.2cm} 
\noindent
Notice that for the classical case $\mathcal{X}=\mathcal{Y}=\mathbb{R}$, 
$r(x) \equiv 1$, $\lambda = \lambda_{L}$ --
where one has $f_{P}(x)=f_{P}(x)$,
$S^{\lambda pd}(P) = S^{pd}(P)$,
and $F_{P}(x) = \int_{-\infty}^{x} f_{P}(z) \mathrm{d}\lambda_{L}(z)$ --
\eqref{BroStuHB22:fo.def.153a} is essentially different from
\eqref{BroStuHB22:fo.def.140aol} with $S(P)=S^{cd}(P)$, $S(Q)=S^{cd}(Q)$ 
which is explicitly of
the ``doubly aggregated form''
\begin{eqnarray} 
& & \hspace{-0.2cm} \textstyle
0 \leq D_{\phi_{\alpha},\mathbb{1},\mathbb{1},\mathbb{1},\lambda}(S^{cd}(P),S^{cd}(Q)) 
\nonumber\\ 
& & \hspace{-0.2cm} 
= \olint_{{\mathbb{R}}} 
\frac{1}{\alpha \cdot (\alpha-1)} \cdot
\Big[ 
\left(\int_{-\infty}^{x} f_{P}(z) \mathrm{d}\lambda_{L}(z)\right)^{\alpha} + 
(\alpha-1) \cdot \left(\int_{-\infty}^{x} f_{Q}(z) \mathrm{d}\lambda_{L}(z)\right)^{\alpha}
\nonumber \\
& & \hspace{-0.2cm} 
 - \alpha 
\cdot \int_{-\infty}^{x} f_{P}(z) \mathrm{d}\lambda_{L}(z) \cdot
 \left(\int_{-\infty}^{x} f_{Q}(z) \mathrm{d}\lambda_{L}(z) \right)^{\alpha-1} 
\Big]
 \, \mathrm{d}\lambda_{L}(x) \ ,
 \quad
\textrm{ for } \alpha \in \mathbb{R}\backslash\{0,1\} ,
\nonumber
\end{eqnarray}
with the usual $\mathrm{d}\lambda_{L}(x)=dx$.

\vspace{0.2cm}
\noindent
In contrast, for the discrete setup 
$(\mathcal{X},\lambda) = (\mathcal{X}_{\#},\lambda_{\#})$ with $\mathcal{X}_{\#} \subset \mathbb{R}$
(recall $\lambda_{\#}[\{x\}] =1$) one has 
$f_{P}(x) = p_{P}(x)$ for all $x \in \mathcal{X}_{\#}$) and
the divergences  
\eqref{BroStuHB22:fo.def.153a} to 
\eqref{BroStuHB22:fo.def.153d}
simplify to
\begin{eqnarray} 
& & \hspace{-0.2cm} \textstyle
0 \leq D_{\phi_{\alpha},1,1,1,\lambda_{\#}}(S^{pm}(P),S^{pm}(Q)) 
\nonumber\\ 
& & \hspace{-0.2cm} 
= \olsum_{{x \in \mathcal{X}}} 
\frac{1}{\alpha \cdot (\alpha-1)} \cdot
\Big[ 
\left(p_{P}(x)\right)^{\alpha} + (\alpha-1) \cdot \left(p_{Q}(x)\right)^{\alpha} - \alpha 
\cdot p_{P}(x) \cdot \left(p_{Q}(x)\right)^{\alpha-1} 
\Big]
\nonumber \\ 
& & \hspace{8.2cm} 
\textrm{ for } \alpha \in \mathbb{R}\backslash\{0,1\},
\nonumber
\\
& & \hspace{-0.2cm} 
0 \leq D_{\phi_{1},1,1,1,\lambda_{\#}}(S^{pm}(P),S^{pm}(Q)) 
= \olsum_{{x \in \mathcal{X}}} 
\Big[ 
p_{P}(x) \cdot \log\left(\frac{p_{P}(x)}{p_{Q}(x)}\right) + p_{Q}(x) - p_{P}(x) 
\Big] \ ,
\nonumber
\\
& & \hspace{-0.2cm} 
0 \leq D_{\phi_{0},1,1,1,\lambda_{\#}}(S^{pm}(P),S^{pm}(Q)) 
= \olsum_{{x \in \mathcal{X}}} 
\Big[ 
- \log\left(\frac{p_{P}(x)}{p_{Q}(x)}\right) + \frac{p_{P}(x)}{p_{Q}(x)} - 1 
\Big] \ , 
\nonumber
\\
& & \hspace{-0.2cm} \textstyle
0 \leq D_{\phi_{2},1,1,1,\lambda_{\#}}(S^{pm}(P),S^{pm}(Q)) 
= \sum_{{x \in \mathcal{X}}} 
\frac{1}{2} \cdot
\Big[ p_{P}(x) - p_{Q}(x) \Big]^2 \, ,  
\nonumber
\end{eqnarray}
where again 
one should exclude $\alpha \leq 0$ whenever $p_{P}(x) =0$ 
for all $x$ in some $A$ with $\lambda_{\#}[A]>0$,
respectively $\alpha \leq 1$ whenever $p_{Q}(x) =0$ 
for all $x$ in some $\tilde{A}$ with $\lambda_{\#}[\tilde{A}]>0$.
For example, take the context from the paragraph right 
after \eqref{BroStuHB22:fo.def.20c},
with discrete random variable $Y$, $p_{Q}(x) = Q[Y=x]$, 
$p_{P}(x) = p_{P_{N}^{emp}}(x)$.
Then, the divergences $2N \cdot D_{\phi_{\alpha},1,1,1,\lambda_{\#}}(S^{pm}(P_{N}^{emp}),S^{pm}(Q))$
\ (for $\alpha \in \mathbb{R}$) can be used as goodness-of-fit test statistics;
see e.g. Ki{\ss}linger \& Stummer~\cite{Kis:16} 
for their limit behaviour as the sample size $N$
tends to infinity.

\vspace{0.2cm}
\noindent
\textbf{Classical quantile functions.} 
The divergence \eqref{BroStuHB22:fo.def.40var} with $S(P)=S^{qu}(P)$, $S(Q)=S^{qu}(Q)$
can be interpreted as a quantitative measure of tail risk
of $P$, relative to some pregiven reference distribution Q
\footnote{
hence, such a divergence 
represents an alternative to Faugeras \& R\"uschendorf~\cite{Fau:18} 
where they use hemimetrics (which e.g. have only a weak identity-property,
but satisfiy triangle inequality)
rather than divergences
}.

\vspace{0.2cm}
\noindent
Especially, for $\mathcal{Y}=\mathbb{R}$ and $\mathcal{X}=(0,1)$, 
$S_{x}(P) = S_{x}^{qu}(P) = F_{P}^{\leftarrow}(x)$,
$S_{x}(Q) = S_{x}^{qu}(Q) = F_{Q}^{\leftarrow}(x)$, 
and the Lebesgue measure $\lambda = \lambda_{L}$ 
(with the usual $\mathrm{d}\lambda_{L}(x)=dx$),
we get from 
\eqref{BroStuHB22:fo.def.140d}
the special case
\begin{eqnarray} 
& & \hspace{-0.2cm} \textstyle
0 \leq D_{\phi_{2},\mathbb{1},\mathbb{1}, r\cdot \mathbb{1},\lambda}(S^{qu}(P),S^{qu}(Q)) 
= \olint_{(0,1)} 
\left( F_{P}^{\leftarrow}(x) - F_{Q}^{\leftarrow}(x) \right)^2 
\, \mathrm{d}\lambda_{L}(x) \ 
\label{BroStuHB22:fo.def.667L2Wass}
\end{eqnarray}
which is nothing but the $2-$Wasserstein distance between the two probability measures
$P$ and $Q$. 
Corresponding connections with optimal transport are discussed in Section \ref{subsec.2.7new} below.
Notice that \eqref{BroStuHB22:fo.def.667L2Wass} does generally
not coincide with its analogue
\begin{eqnarray} 
& & \hspace{-0.2cm} \textstyle
D_{\phi_{2},\mathbb{1},\mathbb{1}, r\cdot \mathbb{1},\lambda}(S^{cd}(P),S^{d}(Q)) 
= \olint_{\mathbb{R}} 
\left( F_{P}(x) - F_{Q}(x) \right)^2 
\, \mathrm{d}\lambda_{L}(x) \ ;
\label{BroStuHB22:fo.def.667L2Wass2}
\end{eqnarray}
to see this, take e.g. $0 < c_{2} < c_{1}$ (e.g. $c_{1}=2$, $c_{2}=1$) and the 
exponential quantile functions
$F_{P}^{\leftarrow}(x) = - \frac{1}{c_{1}} \cdot \log(1- x)$, 
$F_{Q}^{\leftarrow}(x) = - \frac{1}{c_{2}} \cdot \log(1- x)$
for which \eqref{BroStuHB22:fo.def.667L2Wass} becomes $2\cdot (\frac{1}{c_{2}} - \frac{1}{c_{1}})^2$,
whereas for the corresponding exponential distribution functions
$F_{P}(x) = 1- \exp(- c_{1} \cdot x)$, $F_{Q}(x) = 1- \exp(- c_{2} \cdot x)$
the divergence \eqref{BroStuHB22:fo.def.667L2Wass2}
becomes $\frac{1}{2 c_{2}} - \frac{2}{c_{1}+c_{2}} + \frac{1}{2 c_{1}}$.

\vspace{0.2cm}
\noindent
\textbf{Depth, outlyingness, centered rank and centered quantile functions.} \\ 
As a special case one gets 
\begin{eqnarray} 
& & \hspace{-0.2cm} 
D^{c}_{\phi,\mathbb{1},\mathbb{1},r\cdot \mathbb{1},\lambda_{L}}(S^{de}(P),S^{de}(Q)),
\nonumber
\\
& & \hspace{-0.2cm} 
D^{c}_{\phi,\mathbb{1},\mathbb{1},r\cdot \mathbb{1},\lambda_{L}}(S^{ou}(P),S^{ou}(Q)),
\nonumber
\\
& & \hspace{-0.2cm} 
\sum_{i=1}^{d} D^{c}_{\phi,\mathbb{1},\mathbb{1},r\cdot \mathbb{1},\lambda_{L}}(S^{cr,i}(P),S^{cr,i}(Q)),
\label{BroStuHB22:fo.centrank}
\\
& & \hspace{-0.2cm} 
\sum_{i=1}^{d} D^{c}_{\phi,\mathbb{1},\mathbb{1},r\cdot \mathbb{1},\lambda_{L}}(S^{cqu,i}(P),S^{cqu,i}(Q)),
\label{BroStuHB22:fo.centquant}
\end{eqnarray} 
all of which have not appeared elsewhere before (up to our knowledge);
recall that the respective domains of $\phi$ have to take care of the ranges 
$\mathcal{R}\left(S^{de}(P)\right)  
\subset [0,\infty]$,
$\mathcal{R}\left(S^{ou}(P)\right)  
\subset [0,\infty]$,
$\mathcal{R}\left(S^{cr,i}(P)\right)  
\subset [-1,1]$,
$\mathcal{R}\left(S^{cqu,i}(P)\right)  
\subset ]-\infty,\infty[$ ($i \in \{1,\ldots,d\}$).
Notice that these divergences differ structurally from
the Bregman distances 
of Hallin \cite{Hal:18} who uses the
centered rank function $R_{P}(\cdot)$ 
(also called center-outward distribution function)
as a \textit{multidimensional}
(in general not additionally separable) \textit{generator} $\boldsymbol{\phi}$,
and not as points between which the distance is to be measured between.

\vspace{0.3cm}

%
\noindent
\textbf{2.5.1.2 \ \
$\mathbf{m_{1}(x) = m_{2}(x) := S_{x}(Q)}$, $\mathbf{m_{3}(x) = r(x) \cdot S_{x}(Q) \in [0, \infty]}$ for some
(measurable) function $\mathbf{r: \mathcal{X} \rightarrow \mathbb{R}}$ satisfying
$\mathbf{r(x) \in ]-\infty,0[ \cup ]0,\infty[}$ for $\mathbf{\lambda -}$a.a. $\mathbf{x \in \mathcal{X}}$ 
}
%

\vspace{0.3cm}
\noindent 
In such a context, we require that 
the function $r(\cdot)$ does not (explicitly) depend
on the functions $S_{\cdot}(P)$ and $S_{\cdot}(Q)$, 
i.e. it is not of the adaptive form $r(\cdot)= h(\cdot, S_{\cdot}(P), S_{\cdot}(Q))$.
The incorporation of the zeros of $S_{\cdot}(P),S_{\cdot}(Q)$
can be adapted from Broniatowski \& Stummer~\cite{Bro:19a}: 
for instance, in a non-negativity set-up  
where for $\lambda$-almost all $x \in \mathcal{X}$ one has  
$r(x) \in ]0,\infty[$
as well as $S_{x}(P) \in [0,\infty[$,
$S_{x}(Q) \in [0,\infty[$
(as it is always the case for $S^{cd}$, $S^{pd}$, $S^{pm}$, $S^{su}$, $S^{mg}$,
$S^{de}$, $S^{ou}$, 
and for nonnegative real-valued random variables also with $S^{qu}$), 
one can take $E=]a,b[=]0,\infty[$ to end up with
the following special case of \eqref{BroStuHB22:fo.def.20}
\begin{eqnarray} 
& & \hspace{-0.2cm} \textstyle
0 \leq D^{c}_{\phi,S(Q),S(Q),r\cdot S(Q),\lambda}(S(P),S(Q)) 
\nonumber\\ 
& & \hspace{-0.2cm} 
=
\olint_{{\mathcal{X}}} 
\Bigg[ \phi \negthinspace \left( {\frac{S_{x}(P)}{S_{x}(Q)}}\right) 
-\phi \negthinspace \left( 1 \right) 
- \phi_{+,c}^{\prime} \negthinspace
\left( 1 \right) \cdot \left( \frac{S_{x}(P)}{S_{x}(Q)}- 1 \right) 
\Bigg] 
S_{x}(Q) \cdot r(x) \, \mathrm{d}\lambda(x) \ ,
\label{BroStuHB22:fo.def.45a}
\\ 
& & \hspace{-0.2cm} 
=
\olint_{{\mathcal{X}}} 
\Bigg[ S_{x}(Q) \cdot \phi \negthinspace \left( {\frac{S_{x}(P)}{S_{x}(Q)}}\right) 
- S_{x}(Q) \cdot \phi \negthinspace \left( 1 \right) 
- \phi_{+,c}^{\prime} \negthinspace
\left( 1 \right) \cdot \left( S_{x}(P) - S_{x}(Q) \right) 
\Bigg] 
r(x) \, \mathrm{d}\lambda(x) \ ,
\nonumber \\ 
\label{BroStuHB22:fo.def.45b} 
\end{eqnarray} 
\begin{eqnarray} 
& & \hspace{-0.2cm}   \textstyle 
= \int_{{\mathcal{X}}} 
r(x)  \cdot \big[ S_{x}(Q) \cdot \phi  \big( {
\frac{S_{x}(P)}{S_{x}(Q)}}\big) - S_{x}(Q) \cdot \phi  \big( 1 \big) 
- \phi_{+,c}^{\prime} 
\big( 1 \big) \cdot \big( S_{x}(P) - S_{x}(Q) \big) 
\big] 
\, 
\nonumber \\ 
& & \hspace{7.2cm} 
\cdot \boldsymbol{1}_{]0,\infty[}\big(S_{x}(P) \cdot S_{x}(Q) 
\big)
\, \mathrm{d}\lambda(x) 
\nonumber\\ 
& & \hspace{-0.2cm}   \textstyle 
+ \big[ \phi^{*}(0) -\phi_{+,c}^{\prime}(1) \big] \cdot
\int_{{\mathcal{X}}} 
r(x) \cdot S_{x}(P) 
\cdot \boldsymbol{1}_{]0,\infty[}\big(S_{x}(P) 
\big)
\cdot \boldsymbol{1}_{\{0\}}\big(S_{x}(Q)\big)
\, \mathrm{d}\lambda(x)
\nonumber\\ 
& & \hspace{-0.2cm}   \textstyle 
+ 
\big[ \phi(0) + \phi_{+,c}^{\prime}(1) - \phi(1) \big] \cdot 
\int_{{\mathcal{X}}} 
r(x) \cdot S_{x}(Q) 
\cdot \boldsymbol{1}_{]0,\infty[}\big(S_{x}(Q) 
\big)
\cdot \boldsymbol{1}_{\{0\}}\big(S_{x}(P)\big)
\, \mathrm{d}\lambda(x)
\nonumber
\\ 
& & \hspace{-0.2cm}   \textstyle 
=
\int_{{\mathcal{X}}} 
r(x)  \cdot \Big[ S_{x}(Q) \cdot \phi  \big( {
\frac{S_{x}(P)}{S_{x}(Q)}}\big) - S_{x}(Q) \cdot \phi  \big( 1 \big) 
- \phi_{+,c}^{\prime} 
\big( 1 \big) \cdot \big( S_{x}(P) - S_{x}(Q) \big) 
\Big] 
\, 
\nonumber \\ 
& & \hspace{7.2cm} 
\cdot \boldsymbol{1}_{]0,\infty[}\big(S_{x}(P) \cdot S_{x}(Q) 
\big)
\, \mathrm{d}\lambda(x) 
\nonumber\\ 
& & \hspace{-0.2cm}   \textstyle 
+ \big[ \phi^{*}(0) -\phi_{+,c}^{\prime}(1) \big] \cdot
\int_{{\mathcal{X}}} 
r(x) \cdot S_{x}(P) 
\cdot \boldsymbol{1}_{\{0\}}\big(S_{x}(Q)\big)
\, \mathrm{d}\lambda(x)
\nonumber\\ 
& & \hspace{-0.2cm}   \textstyle 
+ 
\big[ \phi(0) + \phi_{+,c}^{\prime}(1) - \phi(1) \big] \cdot 
\int_{{\mathcal{X}}} 
r(x) \cdot S_{x}(Q) 
\cdot \boldsymbol{1}_{\{0\}}\big(S_{x}(P)\big)
\, \mathrm{d}\lambda(x)  \, , 
\label{BroStuHB22:fo.def.47bnew}
\end{eqnarray}
with $\phi^{*}(0) := \lim_{u\rightarrow 0} u \cdot \phi\big(\frac{1}{u}\big) = 
\lim_{v\rightarrow \infty} \frac{\phi(v)}{v}$.
In case of 
$\int_{{\mathcal{X}}} S_{x}(Q) \cdot r(x) \, \mathrm{d}\lambda(x) < \infty$,
the divergence \eqref{BroStuHB22:fo.def.47bnew} becomes

\vspace{-0.5cm}

\begin{eqnarray} 
& & \hspace{-0.2cm}   \textstyle \textstyle
0 \leq D^{c}_{\phi,S(Q),S(Q),r\cdot S(Q),\lambda}(S(P),S(Q)) 
\nonumber\\ 
& & \hspace{-0.2cm}   \textstyle 
=
\int_{{\mathcal{X}}} 
r(x)  \cdot \Big[ S_{x}(Q) \cdot \phi  \big( {
\frac{S_{x}(P)}{S_{x}(Q)}}\big) 
- \phi_{+,c}^{\prime} 
\big( 1 \big) \cdot \big( S_{x}(P) - S_{x}(Q) \big) 
\Big] 
\, 
\cdot \boldsymbol{1}_{]0,\infty[}\big(S_{x}(P) \cdot S_{x}(Q) 
\big)
\, \mathrm{d}\lambda(x) 
\nonumber\\ 
& & \hspace{-0.2cm}   \textstyle 
+ \big[ \phi^{*}(0) -\phi_{+,c}^{\prime}(1) \big] \cdot
\int_{{\mathcal{X}}} 
r(x) \cdot S_{x}(P) 
\cdot \boldsymbol{1}_{\{0\}}\big(S_{x}(Q)\big)
\, \mathrm{d}\lambda(x)
\nonumber\\ 
& & \hspace{-0.2cm}   \textstyle 
+ 
\big[ \phi(0) + \phi_{+,c}^{\prime}(1) \big] \negthinspace \cdot \negthinspace
\int_{{\mathcal{X}}} 
r(x) \cdot S_{x}(Q) \negthinspace
\cdot \negthinspace \boldsymbol{1}_{\{0\}}\big(S_{x}(P)\big)
\, \mathrm{d}\lambda(x)  \, 
 - \phi(1) \negthinspace \cdot \negthinspace
\int_{{\mathcal{X}}} 
r(x) \cdot S_{x}(Q) 
\, \mathrm{d}\lambda(x)  \, . 
\label{BroStuHB22:fo.def.47btwo}
\end{eqnarray}

\noindent
Moreover, in case of $\phi \big( 1 \big) = 0$ and 
$\int_{{\mathcal{X}}}  \big( S_{x}(P) - S_{x}(Q) \big) 
\cdot r(x) \, \mathrm{d}\lambda(x) \in ]-\infty, \infty[$
(but not necessarily 
$\int_{{\mathcal{X}}} S_{x}(P) \cdot r(x) \, \mathrm{d}\lambda(x) < \infty$,
$\int_{{\mathcal{X}}} 
S_{x}(Q) \cdot r(x) \, \mathrm{d}\lambda(x) < \infty$),
the divergence \eqref{BroStuHB22:fo.def.47bnew} turns into

\vspace{-0.4cm}

\begin{eqnarray} 
& & \hspace{-0.2cm}   \textstyle 
0 \leq D^{c}_{\phi,S(Q),S(Q),r\cdot S(Q),\lambda}(S(P),S(Q))
\nonumber\\
& & \hspace{-0.2cm}   \textstyle 
= \int_{{\mathcal{X}}} 
r(x)  \cdot S_{x}(Q) \cdot \phi  \big( {
\frac{S_{x}(P)}{S_{x}(Q)}}\big) 
\cdot \boldsymbol{1}_{]0,\infty[}\big(S_{x}(P) \cdot S_{x}(Q) 
\big)
\, \mathrm{d}\lambda(x) 
\nonumber\\ 
& & \hspace{-0.2cm}   \textstyle 
+ \phi^{*}(0) \cdot
\int_{{\mathcal{X}}} \negthinspace
r(x) \negthinspace \cdot \negthinspace S_{x}(P) \negthinspace
\cdot \negthinspace \boldsymbol{1}_{\{0\}}\big(S_{x}(Q)\big)
\, \mathrm{d}\lambda(x)
+ \phi(0) \cdot 
\int_{{\mathcal{X}}} \negthinspace
r(x) \negthinspace \cdot \negthinspace S_{x}(Q) \negthinspace
\cdot \negthinspace \boldsymbol{1}_{\{0\}}\big(S_{x}(P)\big)
\, \mathrm{d}\lambda(x)  
\nonumber \\ 
& & \hspace{-0.2cm}   \textstyle
- \phi_{+,c}^{\prime}  \big( 1 \big)  \cdot
\int_{{\mathcal{X}}} 
r(x) \cdot \big( S_{x}(P) - S_{x}(Q) \big)
\, \mathrm{d}\lambda(x) \, .
\label{BroStuHB22:fo.def.47cnew}
\end{eqnarray} 

\noindent To obtain the sharp identifiability (reflexivity)
of the divergence\\
 $D^{c}_{\phi,S(Q),S(Q),r\cdot S(Q),\lambda}(S(P),S(Q))$
of \eqref{BroStuHB22:fo.def.47bnew}, one can either use the conditions 
formulated after \eqref{BroStuHB22:fo.reflex1}
in terms of $s \in \mathcal{R}\big(\frac{S(P)}{S(Q)}\big)$
and $t \in \mathcal{R}\big(\frac{S(Q)}{S(Q)}\big) = \{1\}$,
or the strict convexity of $\phi$ at $t=1$ together with
\begin{eqnarray} 
& & \hspace{-0.2cm}   \textstyle
\int_{{\mathcal{X}}} 
\big( S_{x}(P) - S_{x}(Q) \big) \cdot r(x) \, \mathrm{d}\lambda(x) \ = \ 0  
\label{BroStuHB22:fo.AssuCASD1}
\end{eqnarray}
\footnote{and thus, $c$ becomes obsolete}
see Broniatowski \& Stummer~\cite{Bro:19a} for corresponding details.
Additionally, in the light of \eqref{BroStuHB22:fo.def.47btwo}
let us indicate that if one wants to use 
$\Xi := \olint_{{\mathcal{X}}} 
S_{x}(Q) \cdot \phi  \big( {
\frac{S_{x}(P)}{S_{x}(Q)}}\big)
\cdot r(x)  
\, \mathrm{d}\lambda(x) 
$ (with appropriate zero-conventions) as a divergence, then 
one should employ generators $\phi$ satisfying $\phi(1)=\phi_{+,c}^{\prime}(1)=0$,
or employ models fulfilling the assumption 
\eqref{BroStuHB22:fo.AssuCASD1}
together with 
generators $\phi$ with $\phi(1)=0$. On the other hand, if this integral $\Xi$
appears in your application context ``naturally'',  
then one should be aware that 
$\Xi$ may become negative depending on the involved set-up;
for a counter-example, see e.g. Stummer \& Vajda~\cite{Stu:10}.

\vspace{0.2cm}
\noindent
An important generator-concerning example is 
the power-function (limit) case  
$\phi = \phi_{\alpha}$ with 
$\alpha \in \mathbb{R}$
(cf. \eqref{BroStuHB22:fo.def.32}, \eqref{BroStuHB22:fo.def.120b},
\eqref{BroStuHB22:fo.def.120d}, \eqref{BroStuHB22:fo.def.34b})
under the constraint 
$\int_{{\mathcal{X}}}  \big( S_{x}(P) - S_{x}(Q) \big) \cdot r(x) 
\, \mathrm{d}\lambda(x) \in ]-\infty, \infty[$.
Accordingly, the 
``implicit-boundary-describing'' divergence \eqref{BroStuHB22:fo.def.45b}
resp. the corresponding  ``explicit-boundary'' version 
\eqref{BroStuHB22:fo.def.47cnew}
turn into
the generalized power divergences of order $\alpha$
(cf. Stummer \& Vajda~\cite{Stu:10} for $\mathbbm{r}(x) \equiv 1$)

\vspace{-0.6cm}
\begin{eqnarray} 
& & \hspace{-0.2cm}  \textstyle
0 \leq D_{\phi_{\alpha},S(Q),S(Q),r\cdot S(Q),\lambda}(S(P),S(Q)) 
\nonumber\\ 
& & \hspace{-0.2cm}   \textstyle 
=
\olint_{{\mathcal{X}}}
\frac{1}{\alpha \cdot (\alpha-1)} \cdot
\Big[ \big( {
\frac{S_{x}(P)}{S_{x}(Q)}}\big)^{\alpha} - \alpha \cdot \frac{S_{x}(P)}{S_{x}(Q)} + \alpha -1 
\Big] 
\cdot S_{x}(Q) \cdot r(x) \, \mathrm{d}\lambda(x) 
\label{BroStuHB22:fo.def.793a} \\ 
& & \hspace{-0.2cm}   \textstyle 
= \frac{1}{\alpha \cdot (\alpha-1)}  \negthinspace \cdot 
\int_{{\mathcal{X}}} \negthinspace
r(x) \negthinspace \cdot \negthinspace S_{x}(Q) \negthinspace \cdot \negthinspace \Big[ \big( {
\frac{S_{x}(P)}{S_{x}(Q)}}\big)^{\alpha} \negthinspace - \negthinspace
\alpha \negthinspace \cdot \negthinspace \frac{S_{x}(P)}{S_{x}(Q)} \negthinspace + \negthinspace \alpha -1 
\Big] 
\cdot \negthinspace \boldsymbol{1}_{]0,\infty[}\big(S_{x}(P)\negthinspace  \cdot \negthinspace S_{x}(Q) 
\big)
\, \mathrm{d}\lambda(x) 
\nonumber\\ 
& & \hspace{-0.2cm}   \textstyle 
+   \phi_{\alpha}^{*}(0) \negthinspace \cdot \negthinspace
\int_{{\mathcal{X}}} \negthinspace
r(x) \negthinspace \cdot \negthinspace S_{x}(P) 
\negthinspace \cdot \negthinspace \boldsymbol{1}_{\{0\}}\big(S_{x}(Q)\big)
\, \mathrm{d}\lambda(x) 
\negthinspace + \negthinspace \phi_{\alpha}(0) \negthinspace \cdot \negthinspace
\int_{{\mathcal{X}}} \negthinspace
r(x) \negthinspace \cdot \negthinspace S_{x}(Q) 
\negthinspace \cdot \negthinspace \boldsymbol{1}_{\{0\}}\big(S_{x}(P)\big)
\, \mathrm{d}\lambda(x)  
\nonumber \\ 
& & \hspace{-0.2cm}   \textstyle 
= \frac{1}{\alpha \cdot (\alpha-1)} 
\int_{{\mathcal{X}}} 
r(x)  \cdot \Big[ 
S_{x}(P)^{\alpha} \cdot S_{x}(Q)^{1-\alpha} -  S_{x}(Q)
\Big] 
\cdot \boldsymbol{1}_{]0,\infty[}\big(S_{x}(P) \cdot S_{x}(Q) 
\big)
\, \mathrm{d}\lambda(x)
\nonumber\\ 
& & \hspace{-0.2cm}   \textstyle 
+ \frac{1}{1-\alpha} \cdot \int_{{\mathcal{X}}} \negthinspace
\negthinspace r(x) \negthinspace \cdot \negthinspace 
(S_{x}(P) \negthinspace - \negthinspace S_{x}(Q)) \, \mathrm{d}\lambda(x)
\negthinspace + \negthinspace \infty \negthinspace \cdot \negthinspace 
\boldsymbol{1}_{]1,\infty[}(\alpha) \negthinspace \cdot \negthinspace
\int_{{\mathcal{X}}} \negthinspace
\negthinspace r(x) \negthinspace \cdot \negthinspace S_{x}(P) 
\negthinspace \cdot \negthinspace \boldsymbol{1}_{\{0\}}\big(S_{x}(Q)\big)
\, \mathrm{d}\lambda(x)
\nonumber \\ 
& & \hspace{-0.2cm}   \textstyle 
+ \big(\frac{1}{\alpha \cdot (1-\alpha)} \cdot \boldsymbol{1}_{]0,1] \cup ]1,\infty[}(\alpha)
+ \infty \cdot \boldsymbol{1}_{]-\infty,0[}(\alpha) \big) \cdot 
\int_{{\mathcal{X}}} 
r(x) \cdot S_{x}(Q) 
\cdot \boldsymbol{1}_{\{0\}}\big(S_{x}(P)\big)
\, \mathrm{d}\lambda(x) , \footnotemark \qquad \ 
\label{BroStuHB22:fo.def.793b}  \\ 
& & \hspace{-0.2cm}   \textstyle 
0 \leq D_{\phi_{1},S(Q),S(Q),r\cdot S(Q),\lambda}(S(P),S(Q)) 
\nonumber \\
& & \hspace{-0.2cm}   \textstyle
= \olint_{{\mathcal{X}}}
\Big[ 
\frac{S_{x}(P)}{S_{x}(Q)} \cdot \log \big( {\frac{S_{x}(P)}{S_{x}(Q)}} \big) 
+ 1 - \frac{S_{x}(P)}{S_{x}(Q)} 
\Big] 
\cdot S_{x}(Q) \cdot r(x) \, \mathrm{d}\lambda(x) \qquad \ 
\label{BroStuHB22:fo.def.793c} 
\\ 
& & \hspace{-0.2cm}   \textstyle 
= \int_{{\mathcal{X}}} 
r(x)  \cdot S_{x}(P) \cdot \log \big( {\frac{S_{x}(P)}{S_{x}(Q)}} \big) 
\cdot \boldsymbol{1}_{]0,\infty[}\big(S_{x}(P) \cdot S_{x}(Q) 
\big)
\, \mathrm{d}\lambda(x) 
\nonumber\\ 
& & \hspace{-0.2cm}   \textstyle 
+ \int_{{\mathcal{X}}} 
r(x) \cdot (S_{x}(Q) - S_{x}(P)) \, \mathrm{d}\lambda(x) 
+ \infty \cdot
\int_{{\mathcal{X}}} 
r(x) \cdot S_{x}(P) \cdot \boldsymbol{1}_{\{0\}}\big(S_{x}(Q)\big)
\, \mathrm{d}\lambda(x) ,  \qquad \ 
\label{BroStuHB22:fo.def.793d}
 \\
& & \hspace{-0.2cm}   \textstyle 
0 \leq D_{\phi_{0},S(Q),S(Q),r\cdot S(Q),\lambda}(S(P),S(Q)) 
= \olint_{{\mathcal{X}}}
\big[ 
- \log \big( {\frac{S_{x}(P)}{S_{x}(Q)}} \big) 
+ \frac{S_{x}(P)}{S_{x}(Q)} - 1 
\big] 
\cdot S_{x}(Q) \cdot r(x) \, \mathrm{d}\lambda(x)  \qquad \ 
\label{BroStuHB22:fo.def.793e}  
\\ 
& & \hspace{-0.2cm}   \textstyle 
= 
\int_{{\mathcal{X}}} 
r(x)  \cdot S_{x}(Q) \cdot \log \big( {\frac{S_{x}(Q)}{S_{x}(P)}} \big)  
\cdot \boldsymbol{1}_{]0,\infty[}\big(S_{x}(P) \cdot S_{x}(Q) 
\big)
\, \mathrm{d}\lambda(x) 
\nonumber\\ 
& & \hspace{-0.2cm}   \textstyle 
+ \int_{{\mathcal{X}}} 
r(x) \cdot (S_{x}(P) - S_{x}(Q)) \, \mathrm{d}\lambda(x)
+ \infty \cdot
\int_{{\mathcal{X}}} 
r(x) \cdot S_{x}(Q) \cdot \boldsymbol{1}_{\{0\}}\big(S_{x}(P)\big)
\, \mathrm{d}\lambda(x) , \qquad \
\label{BroStuHB22:fo.def.793f}  \\
& & \hspace{-0.2cm}   \textstyle 
0 \leq D_{\phi_{2},S(Q),S(Q),r\cdot S(Q),\lambda}(S(P),S(Q)) 
= \olint_{{\mathcal{X}}}
\frac{1}{2} \cdot
\frac{(S_{x}(P) - S_{x}(Q))^2}{S_{x}(Q)} \cdot r(x) \, \mathrm{d}\lambda(x) 
\label{BroStuHB22:fo.def.793g}
\\ 
& & \hspace{-0.2cm}   \textstyle 
= \frac{1}{2} 
\int_{{\mathcal{X}}} 
r(x)  \cdot \frac{(S_{x}(P) - S_{x}(Q))^2}{S_{x}(Q)}
\cdot \boldsymbol{1}_{[0,\infty[}(S_{x}(P))
\cdot \boldsymbol{1}_{]0,\infty[}(S_{x}(Q))
\, \mathrm{d}\lambda(x) 
\nonumber\\ 
& & \hspace{-0.2cm}   \textstyle 
+ \infty \cdot
\int_{{\mathcal{X}}} 
r(x) \cdot S_{x}(P) 
\cdot \boldsymbol{1}_{\{0\}}\big(S_{x}(Q)\big)
\, \mathrm{d}\lambda(x) \ ,
\label{BroStuHB22:fo.def.793h}
\end{eqnarray}

\vspace{-0.2cm} 

\noindent
which is an adaption of a result of Broniatowski \& Stummer~\cite{Bro:19a}.

Another important generator-concerning example is 
the total variation case  
$\phi_{TV}(t):= |t-1|$
(cf. \eqref{BroStuHB22:fo.def.699})
together with $c=\frac{1}{2}$.
Accordingly, the 
``implicit-boundary-describing'' divergence \eqref{BroStuHB22:fo.def.45b}
resp. the corresponding  ``explicit-boundary'' version 
\eqref{BroStuHB22:fo.def.47cnew}
turn into

\begin{eqnarray} 
& & \hspace{-0.2cm} \textstyle
0 \leq D^{1/2}_{\phi_{TV},S(Q),S(Q),r\cdot S(Q),\lambda}(S(P),S(Q)) 
=
\olint_{{\mathcal{X}}} 
S_{x}(Q) \cdot \left| {\frac{S_{x}(P)}{S_{x}(Q)}} -1 \right|   
\cdot r(x) \, \mathrm{d}\lambda(x) \ 
\nonumber \\ 
& & \hspace{-0.2cm} 
= \int_{{\mathcal{X}}} 
\left| S_{x}(P) - S_{x}(Q) \right| 
\cdot r(x) \, \mathrm{d}\lambda(x) \ ,
\label{BroStuHB22:fo.def.665}
\end{eqnarray}
which is also an adaption of a result of Broniatowski \& Stummer~\cite{Bro:19a}.
Notice that \eqref{BroStuHB22:fo.def.665} -- 
which is nothing but the $r-$weighted $L_{1}-$distance between the two statistical
functionals $S(P)$ and $S(Q)$ -- 
can be used also in cases where 
$S_{x}(P) \in \mathbb{R}$,
$S_{x}(Q) \in \mathbb{R}$, 
and thus e.g. for $S^{qu}$ for arbitrary real-valued random variables.  

\vspace{0.2cm}
\noindent
As usual, for arbitrary discrete setup 
$(\mathcal{X},\lambda) = (\mathcal{X}_{\#},\lambda_{\#})$
all the terms $\int_{{\mathcal{X}}} \ldots \, \mathrm{d}\lambda(x)$ 
(respectively $\olint_{{\mathcal{X}}} \ldots \, \mathrm{d}\lambda(x)$)
in the divergences \eqref{BroStuHB22:fo.def.47bnew} to \eqref{BroStuHB22:fo.def.665} 
turn into $\sum_{x \in \mathcal{X}} \ldots $
(respectively $\olsum_{x \in \mathcal{X}} \ldots $).

\vspace{0.2cm}
\noindent
As far as concrete statistical functionals is concerned, let us briefly discuss several important sub-cases.\\

\noindent
\textbf{$\lambda-$probability-density functions.}
First, in the ``$\lambda-$probability-density functions'' context of Remark \ref{BroStuHB22:rem.40}
one has for general 
$\mathcal{X}$
the statistical functionals $S_{x}^{\lambda pd}(P) := f_{P}(x) \geq 0$, 
$S_{x}^{\lambda pd}(Q) := f_{Q}(x) \geq 0$,
and under the constraints $\phi(1)=0$, the corresponding special case 
$D_{\phi,S^{\lambda pd}(Q),S^{\lambda pd}(Q),r\cdot S^{\lambda pd}(Q),\lambda}(S^{\lambda pd}(P),S^{\lambda pd}(Q))$
of \eqref{BroStuHB22:fo.def.47bnew}
turns out to be the ($r-$)``local $\phi-$divergence of Avlogiaris et al.~\cite{Avl:16a,Avl:16b};
in case of $r(x) \equiv 1$ (where \eqref{BroStuHB22:fo.AssuCASD1} is satisfied),
this reduces to the classical Csiszar-Ali-Silvey-Morimoto~\cite{Csi:63},\cite{Ali:66},\cite{Mori:63} $\phi-$divergence
\footnote{
see e.g. Liese \& Vajda~\cite{Lie:87}, Vajda~\cite{Vaj:89} on comprehensive studies thereupon
} 
\begin{eqnarray} 
& & \hspace{-0.2cm} 
0 \leq D_{\phi,S^{\lambda pd}(Q),S^{\lambda pd}(Q),1\cdot 
S^{\lambda pd}(Q),\lambda}(S^{\lambda pd}(P),S^{\lambda pd}(Q)) 
\nonumber\\ 
& & \hspace{-0.2cm} 
= \int_{{\mathcal{X}}} 
 f_{Q}(x) \cdot \phi \negthinspace \left( {\frac{f_{P}(x)}{f_{Q}(x)}}\right) 
\cdot \boldsymbol{1}_{]0,\infty[}\left(f_{P}(x) \cdot f_{Q}(x)\right)
\, \mathrm{d}\lambda(x) 
\nonumber\\ 
& & \hspace{-0.2cm} 
+ \phi^{*}(0) \cdot
\int_{{\mathcal{X}}} f_{P}(x) 
\cdot \boldsymbol{1}_{\{0\}}\left(f_{Q}(x)\right)
\, \mathrm{d}\lambda(x)
+ \phi(0) \cdot \int_{{\mathcal{X}}} f_{Q}(x) 
\cdot \boldsymbol{1}_{\{0\}}\left(f_{P}(x)\right)
\, \mathrm{d}\lambda(x)  
\nonumber \\ 
& & \hspace{-0.2cm}
- \phi_{+,c}^{\prime} \negthinspace \left( 1 \right)  \cdot
\int_{{\mathcal{X}}} 
\left( f_{P}(x) - f_{Q}(x) \right)
\, \mathrm{d}\lambda(x) \, 
\nonumber\\ 
& & \hspace{-0.2cm} 
= \int_{{\mathcal{X}}} 
 f_{Q}(x) \cdot \phi \negthinspace \left( {\frac{f_{P}(x)}{f_{Q}(x)}}\right) 
\cdot \boldsymbol{1}_{]0,\infty[}\left(f_{P}(x) \cdot f_{Q}(x)\right)
\, \mathrm{d}\lambda(x) 
\nonumber\\ 
& & \hspace{0.2cm} 
+ \phi^{*}(0) \cdot
P[f_{Q}(x) =0] 
+ \phi(0) \cdot 
Q[f_{P}(x) =0] \, 
\label{BroStuHB22:fo.def.47d}
\end{eqnarray}
\footnote{notice that $c$ has become obsolete} 
which coincides with \eqref{BroStuHB22:fo.def.47d.full};
if $\phi(1) \neq 0$ then one has to additionally subtract $\phi(1)$
(cf.\ the corresponding special case of \eqref{BroStuHB22:fo.def.47bnew}).
The corresponding special cases\\ 
$
D_{\phi_{\alpha},S^{\lambda pd}(Q),S^{\lambda pd}(Q),1\cdot 
S^{\lambda pd}(Q),\lambda}(S^{\lambda pd}(P),S^{\lambda pd}(Q))
$
($\alpha \in \mathbb{R}$) of \eqref{BroStuHB22:fo.def.793a} to \eqref{BroStuHB22:fo.def.793h}
are called ``power divergences'' (between
the $\lambda-$density functions $S_{\cdot}^{\lambda pd}(P) := f_{P}(\cdot)$, 
$S_{\cdot}^{\lambda pd}(Q) := f_{Q}(\cdot)$);
if the latter two are strictly positive, the subcase $\alpha=1$ respectively $\alpha=0$ respectively $\alpha=2$
is nothing but the (classical) Kullback-Leibler divergence (relative entropy)
respectively the reverse Kullback-Leibler divergence (reverse relative entropy) 
respectively the Pearson chisquare divergence.
The special case 
\begin{eqnarray} 
& & \hspace{-0.2cm} \textstyle
0 \leq D^{1/2}_{\phi_{TV},S^{\lambda pd}(Q),S^{\lambda pd}(Q),1\cdot S^{\lambda pd}(Q),\lambda}(S^{\lambda pd}(P),S^{\lambda pd}(Q)) =
\int_{{\mathcal{X}}} 
\left| f_{P}(x) - f_{Q}(x) \right|  
\, \mathrm{d}\lambda(x) \ 
\nonumber 
\end{eqnarray}
of \eqref{BroStuHB22:fo.def.665}
is the total variation distance or $L_{1}-$distance
(between
the $\lambda-$density functions $S_{\cdot}^{\lambda pd}(P) := f_{P}(\cdot)$, 
$S_{\cdot}^{\lambda pd}(Q) := f_{Q}(\cdot)$).

\vspace{0.2cm}
\noindent
Analogously to Subsection 2.5.1.1,
for $\mathcal{X}=\mathcal{Y}=\mathbb{R}$ the current context 
subsumes the ``classical density'' functionals
$S^{\lambda pd}(\cdot) = S^{pd}(\cdot)$ with the choice $\lambda = \lambda_{L}$
(and the Riemann integration $\mathrm{d}\lambda_{L}(x)=  \mathrm{d}x$).
In contrast, for the discrete setup  
$\mathcal{Y} = \mathcal{X}= \mathcal{X}_{\#}$ it covers
the ``classical probability mass'' functional $S^{\lambda pd}(\cdot) = S^{pm}(\cdot)$ 
with the choice $\lambda = \lambda_{\#}$
(recall $\lambda_{\#}[\{x\}] =1$ for all $x \in \mathcal{X}_{\#}$); 
accordingly, all the terms $\int_{{\mathcal{X}}} \ldots \, \mathrm{d}\lambda(x)$ in 
the divergences \eqref{BroStuHB22:fo.def.47bnew} to \eqref{BroStuHB22:fo.def.47d} 
turn into $\sum_{x \in \mathcal{X}} \ldots $\ .

\vspace{0.2cm}
\noindent
\textbf{Distribution and survival functions.} \ Let us first consider the context 
$\mathcal{Y}=\mathcal{X}=\mathbb{R}$, 
$S_{x}(P) = S_{x}^{cd}(P) = F_{P}(x)$,
$S_{x}(Q) = S_{x}^{cd}(Q) = F_{Q}(x)$, 
and the Lebesgue measure $\lambda = \lambda_{L}$ 
(with the usual $\mathrm{d}\lambda_{L}(x)=dx$), and $r(x) \equiv 1$.
Therein, the special case 
\begin{eqnarray} 
& & \hspace{-0.2cm} \textstyle
0 \leq D^{1/2}_{\phi_{TV},S^{cd}(Q),S^{cd}(Q),1\cdot S^{cd}(Q),\lambda_{L}}(S^{cd}(P),S^{cd}(Q)) 
= \olint_{\mathbb{R}} 
\left| F_{P}(x) - F_{Q}(x) \right| 
\, \mathrm{d}\lambda_{L}(x) \ 
\nonumber\\
\label{BroStuHB22:fo.def.665kant}
\end{eqnarray}
of \eqref{BroStuHB22:fo.def.665}
is the well-known Kantorovich metric
(between
the distribution functions $F_{P}(\cdot)$,$F_{Q}(\cdot)$).
It is known 
that the integral in \eqref{BroStuHB22:fo.def.665kant}
is finite provided that 
$\int_{\mathcal{X}} x \, \mathrm{d}F_{P}(x) \in ]-\infty,\infty[$ and
$\int_{\mathcal{X}} x \, \mathrm{d}F_{Q}(x) < ]-\infty,\infty[$ 
(if the distribution $P$ resp. $Q$ is generated by 
some real-valued random variable, say $X$ resp. $Y$, this means
that $E[X]$ resp. $E[Y]$ exists and is finite).
To proceed, let us discuss the special case 
\begin{eqnarray} 
& & \hspace{-0.2cm} \textstyle
0 \leq D_{\phi_{1},S^{cd}(Q),S^{cd}(Q),1\cdot S^{cd}(Q),\lambda_{L}}(S^{cd}(P),S^{cd}(Q)) 
\nonumber \\ 
& & \hspace{-0.2cm} \textstyle
= \olint_{\mathbb{R}}
\Big[ 
\frac{F_{P}(x)}{F_{Q}(x)} \cdot \log \big( {\frac{F_{P}(x)}{F_{Q}(x)}} \big) 
+ 1 - \frac{F_{P}(x)}{F_{Q}(x)} 
\Big] 
\cdot F_{Q}(x) \, \mathrm{d}\lambda_{L}(x)
\label{BroStuHB22:fo.def.793spec}
\end{eqnarray}
of \eqref{BroStuHB22:fo.def.793c}, \eqref{BroStuHB22:fo.def.793d}.
For the special subsetup of nonnegative random variables
(and thus $\mathcal{Y}=\mathcal{X}=]0,\infty[$)
with finite expectations and strictly positive cdf,
\eqref{BroStuHB22:fo.def.793spec} simplifies to
the so-called ``cumulative Kullback-Leibler information'' of
Park et al.~\cite{Park:12}
(see also Park et al. \cite{Park:18} for an extension to the whole real line,
Di Crescenzo \& Longobardi~\cite{Dic:15} for an adaption to
possibly smaller support as well as for an adaption to a dynamic form
analogously to the explanations in the following lines).
In contrast, we illuminate the special case
\begin{eqnarray} 
& & \hspace{-0.2cm} \textstyle
0 \leq D_{\phi_{1},S^{su}(Q),S^{su}(Q),1\cdot S^{su}(Q),\lambda_{L}}(S^{su}(P),S^{su}(Q)) 
\nonumber \\ 
& & \hspace{-0.2cm} \textstyle
= \olint_{\mathbb{R}}
\Big[ 
\frac{1-F_{P}(x)}{1-F_{Q}(x)} \cdot \log \big( {\frac{1-F_{P}(x)}{1-F_{Q}(x)}} \big) 
+ 1 - \frac{1-F_{P}(x)}{1-F_{Q}(x)} 
\Big] 
\cdot (1-F_{Q}(x)) \, \mathrm{d}\lambda_{L}(x)
\label{BroStuHB22:fo.def.793spec2}
\end{eqnarray}
of \eqref{BroStuHB22:fo.def.793c}, \eqref{BroStuHB22:fo.def.793d}.
This has been 
employed by Liu \cite{Liu2:07} for the special case of 
$P=P_{N}^{emp}$ and $Q=Q_{\theta}$
in order to obtain corresponding minimum-divergence parameter estimator of $\theta$
(see e.g. also Yari \& Saghafi \cite{Yari:12}, Yari et. al \cite{Yari:13},  
and Mehrali \& Asadi \cite{Meh:21} for follow-up papers).
For the general 
context of nonnegative, absolutely continuous random variables
(and thus $\mathcal{Y}=\mathcal{X}=]0,\infty[$)
with finite expectations and strictly positive cdf,
\eqref{BroStuHB22:fo.def.793spec2} simplifies to
the so-called ``cumulative (residual) Kullback-Leibler information'' of
Baratpour \& Habibi Rad~\cite{Bara:12} (see also Park et al.~\cite{Park:12}
for further properties\footnote{
in this subsetup, they also introduce an alternative
with $\tilde{\phi}_{1}(t)$ of \eqref{BroStuHB22:fo.def.120a}
together with  $S^{su,var}_{x}(P) := \frac{1-F_{P}(x)}{\int_{0}^{\infty}(1-F_{P}(\xi))\mathrm{d}\xi}$  
-- rather than with $\phi_{1}(t)$ of \eqref{BroStuHB22:fo.def.120b} 
together with $S^{su}_{x}(P) := 1-F_{P}(x)$ -- 
(and analogously for $Q$)
} and Park et al. \cite{Park:18} for an extension to the whole real line); the latter has been
adapted to a dynamic form by Chamany \& Baratpour~\cite{Cha:14} as
follows (adapted to our terminolgy): take arbitrarily fixed ``instance'' $t \geq 0$,
$\mathcal{Y}=\mathcal{X} =]t,\infty[$ and
replace in \eqref{BroStuHB22:fo.def.793spec2} the survival function 
$S_{x}^{su}(P) := \left\{1-F_{P}(x)\right\}_{x \in \mathbb{R}}$ by 
$S_{x}^{su,t}(P) := \left\{\frac{1-F_{P}(x)}{1-F_{P}(t)}\right\}_{x \in ]t,\infty[}$
being essentially the survival function of a random variable 
(e.g. residual lifetime) $[X-t | X > t]$ under $P$,
and analogously for $Q$; accordingly, the integral range is $]t,\infty[$.
We can generalize this by simply plugging in $S^{su,t}(P)$, $S^{su,t}(Q)$ into our
general divergences \eqref{BroStuHB22:fo.def.45a}
and \eqref{BroStuHB22:fo.def.40var} --- and even \eqref{BroStuHB22:fo.def.1} 
--- (with $\lambda = \lambda_{L}$).
An analogous dynamization can be done for density-functionals,
by plugging  $S^{\lambda pd, t} := \left\{\frac{f_{P}(x)}{1-F_{P}(t)}\right\}_{x \in ]t,\infty[}$
instead of $S^{\lambda pd} = \left\{f_{P}(x)\right\}_{x \in \mathbb{R}}$ 
into \eqref{BroStuHB22:fo.def.45a}
and \eqref{BroStuHB22:fo.def.40var} --- and even \eqref{BroStuHB22:fo.def.1} ---
and thus covering the corresponding dynamic Kullback-Leibler divergence
of Ebrahimi \& Kirmani~\cite{Ebr:96} as well as the more general \textit{$\phi-$divergences
between residual lifetimes} of Vonta \& Karagrigoriou \cite{Vonta:10} as special cases;
notice that $S^{\lambda pd, t}$ is essentially the density function of the 
random variable $X_{t}:= [X-t | X > t]$ under $P$, where e.g.
$X$ is typically a (non-negative) absolutely continuous random variable
which describes the residual lifetime of a person or an item or a ``process''
and hence, $X_{t}$ is called \textit{residual lifetime} (at $t$) which is fundamentally used in
survival analysis and systems reliability engineering. In risk management
and extreme value theory, $X_{t}$ describes the important notion of random excess 
(e.g. of a loss $X$) over the threshold $t$, which is e.g. employed in the
well-known peaks-over-threshold method.

Analogously, we can plug in 
$\widetilde{S}^{\lambda pd, t} := \left\{\frac{f_{P}(x)}{F_{P}(t)}\right\}_{x \in ]t,\infty[}$
instead of $S^{\lambda pd} = \left\{f_{P}(x)\right\}_{x \in \mathbb{R}}$ 
into \eqref{BroStuHB22:fo.def.45a}
and \eqref{BroStuHB22:fo.def.40var} --- and even \eqref{BroStuHB22:fo.def.1} ---
and thus covering the corresponding dynamic Kullback-Leibler divergence
of Di Crescenzo \& Longobardi~\cite{Dic:04} as well as the more general 
\textit{$\phi-$divergences
between past lifetimes} of Vonta \& Karagrigoriou \cite{Vonta:10} as special cases;
notice that $\widetilde{S}^{\lambda pd, t}$ is essentially the density function of the 
random variable $[X | X \leq t]$ under $P$.

\vspace{0.2cm}
\noindent
\textbf{Classical quantile functions.} 
The divergence \eqref{BroStuHB22:fo.def.45a} with $S(P)=S^{qu}(P)$, $S(Q)=S^{qu}(Q)$
can be interpreted as a quantitative measure of tail risk
of $P$, relative to some pregiven reference distribution Q
\footnote{
hence, such a divergence 
represents an alternative to Faugeras \& R\"uschendorf~\cite{Fau:18} 
where they use hemimetrics rather than divergences
}.

\vspace{0.2cm}
\noindent
For $\mathcal{Y}=\mathbb{R}$ and $\mathcal{X}=(0,1)$, we get for the quantiles context
\begin{eqnarray} 
& & \hspace{-0.2cm} \textstyle
0 \leq D^{1/2}_{\phi_{TV},S^{qu}(Q),S^{qu}(Q),1\cdot S^{qu}(Q),\lambda_{L}}(S^{qu}(P),S^{qu}(Q)) 
= \olint_{{\mathcal{X}}} 
\left| F_{P}^{\leftarrow}(x) - F_{Q}^{\leftarrow}(x) \right| 
\, \mathrm{d}\lambda_{L}(x) \  
\nonumber\\
\label{BroStuHB22:fo.def.666L1Wass}
\end{eqnarray}
which is nothing but the $1-$Wasserstein distance between the two probability measures
$P$ and $Q$. It is well-known that the right-hand sides of 
\eqref{BroStuHB22:fo.def.665kant} and \eqref{BroStuHB22:fo.def.666L1Wass}
coincide, in contrast to the 
discussion on the ``$L_{2}-$case'' right after \eqref{BroStuHB22:fo.def.667L2Wass2}.
Corresponding connections with optimal transport are discussed 
in Section \ref{subsec.2.7new} below.

\vspace{0.2cm}
\noindent
Let us briefly discuss some other connections between $\phi-$divergences and quantile functions.
In the above-mentioned setup of Baratpour \& Habibi Rad~\cite{Bara:12}
(under the existence of strictly positive probability density functions), 
Sunoj et al.~\cite{Sun:18} rewrite the
cumulative Kullback-Leibler information (cf. the special case of \eqref{BroStuHB22:fo.def.793spec2})
equivalently in terms of quantile functions.
In contrast, in a context of absolutely continuous probability distributions
$P$ and $Q$
on $\mathcal{X}=\mathbb{R}$ with strictly positive density functions
$f_{P}$ and $f_{Q}$,
Sankaran et al.~\cite{San:16}
rewrite the classical Kullback-Leibler divergence
$\int_{{\mathcal{X}}} 
\Big[ 
f_{P}(x) \cdot \log\left(\frac{f_{P}(x)}{f_{Q}(x)}\right) + f_{Q}(x) - f_{P}(x) 
\Big]
\, \mathrm{d}\lambda_{L}(x) \ =
D_{\phi_{1},S^{\lambda_{\negthinspace L} pd}(Q),S^{\lambda_{\negthinspace L} pd}(Q),
1\cdot S^{\lambda_{\negthinspace L} pd}(Q),\lambda_{\negthinspace L}}
(S^{\lambda_{\negthinspace L} pd}(P),S^{\lambda_{\negthinspace L} pd}(Q))
$
(cf. \eqref{BroStuHB22:fo.def.793c})
equivalently in terms of quantile functions;
in the same setup, for $\alpha \in ]0,1[ \cup ]1,\infty[$ 
Kayal \& Tripathy \cite{Kay:18} rewrite the
classical $\alpha-$order power divergences
(in fact, the classical $\alpha-$order Tsallis cross-entropies which are 
multiples thereof)
$
\int_{{\mathcal{X}}}
\frac{1}{\alpha \cdot (\alpha-1)} \cdot
\Big[ \big( {
\frac{f_{P}(x)}{f_{Q}(x)}}\big)^{\alpha} - 
\alpha \cdot \frac{f_{P}(x)}{f_{Q}(x)} + \alpha -1 
\Big] 
\cdot f_{Q}(x) \, \mathrm{d}\lambda(x)
=
D_{\phi_{\alpha},S^{\lambda_{\negthinspace L} pd}(Q),S^{\lambda_{\negthinspace L} pd}(Q),
1\cdot S^{\lambda_{\negthinspace L} pd}(Q),\lambda_{\negthinspace L}}
(S^{\lambda_{\negthinspace L} pd}(P),S^{\lambda_{\negthinspace L} pd}(Q))$
(cf. \eqref{BroStuHB22:fo.def.793a})
equivalently in terms of quantile functions, where they also emphasize the
advantage for distributions $P$ and $Q$ having closed-form quantile functions
but non-closed-form distribution functions.

The above-mentioned contexts differ considerably from that of
Broniatowski \& Decurninge \cite{Bro:16}, who basically
employ $\phi-$divergences $D_{\phi}(\mathfrak{Q}_{Q},\mathfrak{Q}_{P})$ 
between special \textit{quantile measures}
(rather than quantile functions) $\mathfrak{Q}_{Q}$ and $\mathfrak{Q}_{P}$;
recall that for any probability measure $P$ on $\mathbb{R}$,
one can associate a (signed) quantile measure $\mathfrak{Q}_{P}$ on $]0,1[$
having as its generalized distribution function nothing else but the 
quantile function $F_{P}^{\leftarrow}$ of $P$.
In more detail, similarly to the above-mentioned empirical likelihood
principle, \cite{Bro:16} consider 
--- in an i.i.d. context --- the minimization 
\begin{equation}
D_{\phi}(\Omega_{N}^{dis}, \mathfrak{Q}_{P_{N}^{emp}} ) := 
\inf_{\mathfrak{Q}_{Q}\in\Omega_{N}^{dis}} 
D_{\phi}( \mathfrak{Q}_{Q}, \mathfrak{Q}_{P_{N}^{emp}} )
\nonumber
\end{equation}
of the $\phi-$divergences 
$D_{\phi}(\mathfrak{Q}_{Q},\mathfrak{Q}_{P_{N}^{emp}})$,
where $\Omega_{N}^{dis}$ is the subclass of quantile measures
$\mathfrak{Q}_{Q}$ having support on $\left\{ \frac{i}{n},1\leq i\leq n\right\} $
of a desired model $\Omega$ of quantile measures
$\mathfrak{Q}_{\widetilde{Q}}$ having support on $\mathbb{R}$;
for example, the $\widetilde{Q}$'s may be taken from 
a tubular neighborhood 
$\Lambda$ --- constructed through a finite collection of conditions
on $L-$moments (cf. e.g. Hosking \cite{Hos:90}) ---
of some class of distributions on $\mathbb{R}^{+}$, such as the
Pareto- or Weibull-distribution class.
Such tasks have numerous applications in climate sciences or hydrology.
As a side remark, let us mention that for the general context of quantile measures
$\mathfrak{Q}_{Q}$ and $\mathfrak{Q}_{P}$ being absolutely continuous
(with respect to the Lebesgue measure $\lambda_{L}$ on $[0,1]$),
the $\phi-$divergence $D_{\phi}(\mathfrak{Q}_{Q},\mathfrak{Q}_{P})$
turns into to the divergence 
$D^{c}_{\phi,S^{qd}(P),S^{qd}(P),S^{qd}(P),\lambda_{L}}(S^{qd}(Q),S^{qd}(P))$
(cf. \eqref{BroStuHB22:fo.def.45a})
between the \textit{quantile density functions}
$S^{qd}(P):= \left\{S_{x}^{qd}(P)\right\}_{x \in ]0,1[} 
:= \left\{\big(F_{P}^{\leftarrow}\big)^{\prime}(x)\right\}_{x \in ]0,1[}$
and $S^{qd}(Q)$. Thus, by applying our general divergences 
\eqref{BroStuHB22:fo.def.1}
to $S^{qd}(Q)$ and $S^{qd}(P)$ we end up
with a completely new framework 
$D^{c}_{\phi,m_{1},m_{2},m_{3},\lambda}(S^{qd}(Q),S^{qd}(P))$
(and many interesting special cases) 
for quantifying dissimilarities between
quantile density functions.

\vspace{0.3cm}
\noindent
\textbf{Depth, outlyingness, centered rank and centered quantile functions.} \\ 
As a special case one gets 
$D^{c}_{\phi,S^{de}(Q),S^{de}(Q),r\cdot S^{de}(Q),\lambda_{L}}(S^{de}(P),S^{de}(Q))$,\\
$D^{c}_{\phi,S^{ou}(Q),S^{ou}(Q),r\cdot S^{ou}(Q),\lambda_{L}}(S^{ou}(P),S^{ou}(Q))$,\\
$\sum_{i=1}^{d} D^{c}_{\phi,S^{cr,i}(Q),S^{cr,i}(Q),r\cdot S^{cr,i}(Q),
\lambda_{L}}(S^{cr,i}(P),S^{cr,i}(Q))$,\\
$\sum_{i=1}^{d} D^{c}_{\phi,S^{cqu,i},S^{cqu,i},r\cdot S^{cqu,i},
\lambda_{L}}(S^{cqu,i}(P),S^{cqu,i}(Q))$,\\
all of which have not appeared elsewhere before (up to to our knowledge);
recall that the respective domains of $\phi$ have to take care of the ranges 
$\mathcal{R}\left(S^{de}(P)\right)  
\subset [0,\infty]$,
$\mathcal{R}\left(S^{ou}(P)\right)  
\subset [0,\infty]$,
$\mathcal{R}\left(S^{cr,i}(P)\right)  
\subset [-1,1]$,
$\mathcal{R}\left(S^{cqu,i}(P)\right)  
\subset ]-\infty,\infty[$ ($i \in \{1,\ldots,d\}$).

\vspace{0.4cm}

%
\noindent
\textbf{2.5.1.3 \ \
$\mathbf{m_{1}(x) = m_{2}(x) := w(S_{x}(P),S_{x}(Q))}$, 
$\mathbf{m_{3}(x) = r(x) \cdot w(S_{x}(P),S_{x}(Q)) \in [0, \infty[}$ for some
(measurable) functions
$\mathbf{w: \mathcal{R}(S(P)) \times \mathcal{R}(S(Q)) \rightarrow \mathbb{R}}$
and $\mathbf{r: \mathcal{X} \rightarrow \mathbb{R}}$
}
%

Such a choice extends the contexts of the previous Subsections 2.5.1.1 resp. 2.5.1.2
(where the ``connector function'' $w$ took the simple form $w(u,v) = 1$ resp. $w(u,v) = v$).
This introduces a wide adaptive modeling flexibility, 
where \eqref{BroStuHB22:fo.def.20} specializes to
\begin{eqnarray} 
& & \hspace{-0.2cm} \textstyle
0 \leq D^{c}_{\phi,w(S(P),S(Q)),w(S(P),S(Q)),r\cdot w(S(P),S(Q)),\lambda}(S(P),S(Q)) 
\nonumber\\ 
& & \hspace{-0.2cm} \textstyle
: =
\olint_{{\mathcal{X}}} 
\Bigg[ \phi \negthinspace \left( {
\frac{S_{x}(P)}{w(S_{x}(P),S_{x}(Q))}}\right) -\phi \negthinspace \left( 
{\frac{S_{x}(Q)}{w(S_{x}(P),S_{x}(Q))}}\right)
\nonumber\\ 
& & \hspace{-0.2cm} \textstyle 
- \phi_{+,c}^{\prime} \negthinspace
\left( {\frac{S_{x}(Q)}{w(S_{x}(P),S_{x}(Q))}}\right) \cdot \left( 
\frac{S_{x}(P)}{w(S_{x}(P),S_{x}(Q))}-\frac{S_{x}(Q)}{w(S_{x}(P),S_{x}(Q))}\right) 
\Bigg] \cdot
w(S_{x}(P),S_{x}(Q)) \cdot r(x) \, \mathrm{d}\lambda(x) \ ,
\nonumber \\ 
\label{BroStuHB22:fo.def.20w}
\end{eqnarray}
which for the discrete setup 
$(\mathcal{X},\lambda) = (\mathcal{X}_{\#},\lambda_{\#})$ 
(recall $\lambda_{\#}[\{x\}] =1$ for all $x \in \mathcal{X}_{\#}$) 
simplifies to
\begin{eqnarray}
& & \hspace{-0.2cm} \textstyle
0 \leq D^{c}_{\phi,w(S(P),S(Q)),w(S(P),S(Q)),r\cdot w(S(P),S(Q)),\lambda_{\#}}(S(P),S(Q))
\nonumber\\ 
& & \hspace{-0.2cm} \textstyle
 =
\olsum_{{x \in \mathcal{X}}} 
\Bigg[ \phi \negthinspace \left( {
\frac{S_{x}(P)}{w(S_{x}(P),S_{x}(Q))}}\right) -
\phi \negthinspace \left( {\frac{S_{x}(Q)}{w(S_{x}(P),S_{x}(Q))}}\right)
\nonumber\\ 
& & \hspace{-0.2cm}  \textstyle
- \phi_{+,c}^{\prime} \negthinspace
\left( {\frac{S_{x}(Q)}{w(S_{x}(P),S_{x}(Q))}}\right) \cdot 
\left( \frac{S_{x}(P)}{w(S_{x}(P),S_{x}(Q))}-\frac{S_{x}(Q)}{w(S_{x}(P),S_{x}(Q))}\right) 
\Bigg] \cdot
w(S_{x}(P),S_{x}(Q)) \cdot r(x) \ .
\nonumber \\ 
\label{BroStuHB22:fo.def.21w}
\end{eqnarray}
As a side remark, let us mention that  
by appropriate choices of $w(\cdot,\cdot)$ and $\phi$ in \eqref{BroStuHB22:fo.def.20w}
we can even derive 
divergences of the form \eqref{BroStuHB22:fo.def.47cnew}
but with non-convex non-concave $\phi$: see e.g. 
the ``perturbed'' power divergences of
Roensch \& Stummer~\cite{Roe:17}.

\vspace{0.2cm}
\noindent
In the following, let us illuminate 
the important special case of \eqref{BroStuHB22:fo.def.21w} with 
$\phi = \phi_{\alpha}$ ($\alpha \in \mathbb{R}$,
cf. \eqref{BroStuHB22:fo.def.32}, \eqref{BroStuHB22:fo.def.120b},
\eqref{BroStuHB22:fo.def.120d}, \eqref{BroStuHB22:fo.def.34b})
together with 
$S_{x}(P) \geq 0$,
$S_{x}(Q) \geq 0$ 
(as it is always the case for $S^{cd}$, $S^{pd}$, $S^{pm}$, $S^{su}$, $S^{mg}$,
$S^{de}$, $S^{ou}$, 
and for nonnegative real-valued random variables also with $S^{qu}$): 

\begin{eqnarray} 
& & \hspace{-0.2cm}   \textstyle 
0 \leq D_{\phi_{\alpha},w(S(P),S(Q)),w(S(P),S(Q)),r\cdot w(S(P),S(Q)),\lambda}(S(P),S(Q))
\nonumber\\ 
& & \hspace{-0.2cm}   \textstyle 
= \olint_{{\mathcal{X}}} 
\frac{r(x) \cdot (w(S_{x}(P),S_{x}(Q)))^{1-\alpha}}{\alpha \cdot (\alpha-1)} \cdot
\Big[ 
\left(S_{x}(P)\right)^{\alpha} + (\alpha-1) \cdot \left(S_{x}(Q)\right)^{\alpha} 
\nonumber\\
& & \hspace{-0.2cm}   \textstyle 
- \alpha 
\cdot S_{x}(P) \cdot \left(S_{x}(Q)\right)^{\alpha-1} 
\Big]
 \, \mathrm{d}\lambda(x) \ ,
 \hspace{2.0cm} \textrm{ for } \alpha \in \mathbb{R}\backslash\{0,1\},
\label{BroStuHB22:fo.def.513}
\\
& & \hspace{-0.2cm}   \textstyle 
0 \leq D_{\phi_{1},w(S(P),S(Q)),w(S(P),S(Q)),r\cdot w(S(P),S(Q)),\lambda}(S(P),S(Q)) 
\nonumber\\ 
& & \hspace{-0.2cm}   \textstyle 
= \olint_{{\mathcal{X}}} 
r(x)  \cdot \big[ 
S_{x}(P) \cdot \log\big(\frac{S_{x}(P)}{S_{x}(Q)}\big) + S_{x}(Q) - S_{x}(P) 
\big] \, \mathrm{d}\lambda(x) ,
\label{BroStuHB22:fo.def.514} 
\\ 
& & \hspace{-0.2cm}   \textstyle 
0 \leq D_{\phi_{0},w(S(P),S(Q)),w(S(P),S(Q)),r\cdot w(S(P),S(Q)),\lambda}(S(P),S(Q)) 
\nonumber\\ 
& & \hspace{-0.2cm}   \textstyle 
= \olint_{{\mathcal{X}}} 
r(x)  \cdot w(S_{x}(P),S_{x}(Q)) \cdot \Big[ 
- \log\big(\frac{S_{x}(P)}{S_{x}(Q)}\big) + \frac{S_{x}(P)}{S_{x}(Q)} - 1 
\Big] \, \mathrm{d}\lambda(x), 
\label{BroStuHB22:fo.def.515a}
\\ 
& & \hspace{-0.2cm}   \textstyle 
0 \leq D_{\phi_{2},w(S(P),S(Q)),w(S(P),S(Q)),r\cdot w(S(P),S(Q)),\lambda}(S(P),S(Q)) 
\nonumber\\ 
& & \hspace{-0.2cm}   \textstyle 
= \int_{{\mathcal{X}}} 
\frac{r(x)}{2} \cdot
\frac{( S_{x}(P) - S_{x}(Q) )^2}{w(S_{x}(P),S_{x}(Q))}
 \, \mathrm{d}\lambda(x) \ . 
\label{BroStuHB22:fo.def.516a}
\end{eqnarray} 

\vspace{0.2cm}
\noindent
\textbf{$\lambda-$probability-density functions.} 
For general $\mathcal{X}$, $r(x) =1$, and 
(cf. Remark \ref{BroStuHB22:rem.40}(c))
$S_{x}(P) = S_{x}^{\lambda pd}(P) := f_{P}(x) \geq 0$, 
$S_{x}(Q) = S_{x}^{\lambda pd}(Q) = f_{Q}(x) \geq 0$,
the divergences \eqref{BroStuHB22:fo.def.20w}, \eqref{BroStuHB22:fo.def.21w},
\eqref{BroStuHB22:fo.def.513} to \eqref{BroStuHB22:fo.def.516a}
are due to Kisslinger \& Stummer~\cite{Kis:13},~\cite{Kis:15a}, ~\cite{Kis:16}
(where they also gave indications on non-probability measures).
Recall that this directly subsumes for $\mathcal{X}=\mathcal{Y}=\mathbb{R}$ the 
``classical density'' functional
$S^{\lambda pd}(\cdot) = S^{pd}(\cdot)$ with the choice $\lambda = \lambda_{L}$
(and the Riemann integration $\mathrm{d}\lambda_{L}(x)=  \mathrm{d}x$),
as well as for the discrete setup  
$\mathcal{Y} = \mathcal{X}= \mathcal{X}_{\#}$
the ``classical probability mass'' functional $S^{\lambda pd}(\cdot) 
= S^{pm}(\cdot)$ with the choice $\lambda = \lambda_{\#}$.\\

\noindent
\textbf{Distribution functions.} Recall that $\mathcal{Y}=\mathcal{X} = \mathbb{R}$, 
$S_{x}(P) = S_{x}^{cd}(P) = F_{P}(x)$,
$S_{x}(Q) = S_{x}^{cd}(Q) = F_{Q}(x)$. Let us illuminate 
\eqref{BroStuHB22:fo.def.516a}, 
for the setup of a real-valued random variable $Y$, 
$F_{Q}(x) = Q[Y\leq x]$ 
under a hypothetical/candidate law $Q$,
$F_{P}(x) = 
\frac{1}{N} \cdot \# \{ i \in \{ 1, \ldots, N\}: Y_i \leq x \} =: F_{P_{N}^{emp}}(x)$
as the distribution function of the 
corresponding data-derived ``empirical distribution'' 
$P:= P_{N}^{emp} := \frac{1}{N} \cdot \sum_{i=1}^{N}   \delta_{Y_{i}}[\cdot]$
of an $N-$size i.i.d. sample $Y_1, \ldots, Y_N$ of $Y$.
In such a set-up, the choice $\lambda = Q$
in \eqref{BroStuHB22:fo.def.516a} and multiplication with $2N$
lead to 
\begin{eqnarray}
& & \hspace{-0.2cm}   \textstyle 
0 \leq 2N \cdot D_{\phi_{2},w(S^{cd}(P_{N}^{emp}),S^{cd}(Q)),w(S^{cd}(P_{N}^{emp}),S^{cd}(Q)),
r\cdot w(S^{cd}(P_{N}^{emp}),S^{cd}(Q)),Q}(S^{cd}(P_{N}^{emp}),S^{cd}(Q)) 
\nonumber\\ 
& & \hspace{-0.2cm}   \textstyle 
= N \cdot \int_{{\mathbb{R}}} r(x) \cdot
\frac{( F_{P_{N}^{emp}}(x) - F_{Q}(x) )^2}{w(F_{P_{N}^{emp}}(x),F_{Q}(x))}
 \, \mathrm{d}Q(x) \ . 
\label{BroStuHB22:fo.def.516b}
\end{eqnarray} 
The special case $w(u,v) =1$ reduces to the 
\textit{Cramer-von Mises (test statistics) family}
\eqref{BroStuHB22:fo.def.140e},
and the choice $r(x)=1$, $w(u,v) = v \cdot (1-v)$
gives the \textit{Anderson-Darling~\cite{And:52} test statistics}.
With \eqref{BroStuHB22:fo.def.516b}, we can also imbed as special cases (together with $r(x)=1$) 
some other known divergences which emphasize the \textit{upper tails}:
$w(u,v) = 1-v$ (cf. Ahmad et al. \cite{Ahm:88}), 
$w(u,v) = 1-v^{2}$ (cf. Rodriguez \& Viollaz \cite{Rod:95}, 
see also Shin et al. \cite{Shin:12}
for applications in environmental extreme-value theory),
$w(u,v) = (1-v)^{\beta}$ with $\beta >0$ 
(cf. Deheuvels \& Martynov \cite{Deh:03}, see also Chernobai et al. \cite{Chernobai:15}
for the case $\beta=2$ together with a left-truncated version of the
empirical distribution function).
Moreover, \eqref{BroStuHB22:fo.def.516b} covers as special cases (together with $r(x)=1$) 
some other known divergences which emphasize the \textit{lower tails}:
$w(u,v) = v$ (cf. Ahmad et al. \cite{Ahm:88}, Scott \cite{Sco:99}),
$w(u,v) = v^{\beta}$ with $\beta >0$ (cf. Deheuvels \& Martynov \cite{Deh:03}),
$w(u,v) = v \cdot (2-v)$ (cf. Rodriguez \& Viollaz \cite{Rod:95},
see also Shin et al. \cite{Shin:12}).
In contrast, in a two-sample-test situation where
$Q$ is replaced by the empirical distribution 
$\widetilde{P}_{L}^{emp} := \frac{1}{L} \cdot \sum_{i=1}^{L}   
\delta_{\widetilde{Y}_{i}}[\cdot]$
of a $L-$size i.i.d. sample $\widetilde{Y}_1, \ldots, \widetilde{Y}_N$ of $Y$ (under $Q$),
some authors (e.g. Rosenblatt \cite{Ros:52}, Hajek et al. \cite{Haj:99}) 
choose divergences which can be imbedded 
(with the choice $w(u,v) =1$, $r(x) =1$)
in our framework as multiple of
$D_{\phi_{2},\mathbb{1},\mathbb{1},\mathbb{1}, \lambda}(S^{cd}(P_{N}^{emp}),
S^{cd}(\widetilde{P}_{L}^{emp}))$
where $\lambda = c_{1} \cdot P_{N}^{emp} + (1-c_{1}) \cdot \widetilde{P}_{L}^{emp}$ 
is an appropriate mixture with $c_{1} \in ]0,1[$. 
In further contrast, if one chooses the Lebesgue measure
$\lambda = \lambda_{L}$ 
(with the usual Riemann integration $\mathrm{d}\lambda_{L}(x)=  \mathrm{d}x$)
and $r(x) \equiv 1$ in \eqref{BroStuHB22:fo.def.516b},
then one ends up with an adaptively weighted extension of \eqref{BroStuHB22:fo.def.140f}.

\vspace{0.2cm}
\noindent
\textbf{Classical quantile functions.} \\
The divergence \eqref{BroStuHB22:fo.def.20w} with $S(P)=S^{qu}(P)$, $S(Q)=S^{qu}(Q)$,
$\lambda = \lambda_{L}$,\\
i.e. $D^{c}_{\phi,w(S^{qu}(P),S^{qu}(Q)),w(S^{qu}(P),S^{qu}(Q)),r\cdot w(S^{qu}(P),S^{qu}(Q)),
\lambda_{L}}(S^{qu}(P),S^{qu}(Q))$
-- which has first been given in Stummer~\cite{Stu:21} in an even more flexible form --
can be interpreted as a quantitative measure of tail risk
of $P$, relative to some pregiven reference distribution Q;
corresponding connections with optimal transport are discussed 
in Section \ref{subsec.2.7new} below.

\vspace{0.2cm}
\noindent
\textbf{Depth, outlyingness, centered rank and centered quantile functions.} \\ 
As a special case of \eqref{BroStuHB22:fo.def.20w} one gets\\ 
$D^{c}_{\phi,w(S^{de}(Q),S^{de}(Q)),w(S^{de}(Q),S^{de}(Q)),r\cdot w(S^{de}(Q),S^{de}(Q)),\lambda_{L}}(S^{de}(P),S^{de}(Q))$,\\
$D^{c}_{\phi,w(S^{ou}(Q),S^{ou}(Q)),w(S^{ou}(Q),S^{ou}(Q)),r\cdot w(S^{ou}(Q),S^{ou}(Q)),\lambda_{L}}(S^{ou}(P),S^{ou}(Q))$,\\
$\sum_{i=1}^{d} D^{c}_{\phi,w(S^{cr,i}(Q),S^{cr,i}(Q)),w(S^{cr,i}(Q),S^{cr,i}(Q)),r\cdot 
w(S^{cr,i}(Q),S^{cr,i}(Q)),\lambda_{L}}(S^{cr,i}(P),S^{cr,i}(Q))$,\\
$\sum_{i=1}^{d} D^{c}_{\phi,w(S^{cqu,i},S^{cqu,i}),w(S^{cqu,i},S^{cqu,i}),r\cdot 
w(S^{cqu,i},S^{cqu,i}),\lambda_{L}}(S^{cqu,i}(P),S^{cqu,i}(Q))$,\\
all of which have not appeared elsewhere before (up to to our knowledge);
recall that the respective domains of $\phi$ have to take care of the ranges 
$\mathcal{R}\left(S^{de}(P)\right)  
\subset [0,\infty]$,
$\mathcal{R}\left(S^{ou}(P)\right)  
\subset [0,\infty]$,
$\mathcal{R}\left(S^{cr,i}(P)\right)  
\subset [-1,1]$,
$\mathcal{R}\left(S^{cqu,i}(P)\right)  
\subset ]-\infty,\infty[$ ($i \in \{1,\ldots,d\}$).

\vspace{0.3cm}

%
\noindent
\textbf{2.5.2 \ \
$\mathbf{m_{1}(x) =  \tilde{S}_{x}(P)}$ and
$\mathbf{m_{2}(x) =  \tilde{S}_{x}(Q)}$ with statistical functional $\tilde{S} \ne S$, \ 
$\mathbf{m_{3}(x) \geq 0}$ 
}
%

\vspace{0.3cm}
\noindent
Recall $S(P) := \left\{S_{x}(P)\right\}_{x \in \mathcal{X}}$, $S(Q) := \left\{S_{x}(Q)\right\}_{x \in \mathcal{X}}$,
and let $\tilde{S}(P) := \left\{\tilde{S}_{x}(P)\right\}_{x \in \mathcal{X}}$, 
$\tilde{S}(Q) := \left\{\tilde{S}_{x}(Q)\right\}_{x \in \mathcal{X}}$
for (typically) $\tilde{S}$ being ``essentially different'' to $S$
(e.g., take $\tilde{S}$ and $S$ as different choices from 
$S^{cd}$, $S^{pd}$, $S^{pm}$, $S^{su}$, $S^{mg}$, $S^{qu}$, $S^{de}$, $S^{ou}$).

\vspace{0.2cm}
\noindent
For this special case, from  \eqref{BroStuHB22:fo.def.1} one can deduce
\begin{eqnarray} 
& & \hspace{-0.2cm} \textstyle
0 \leq D_{\phi,\tilde{S}(P),\tilde{S}(Q),m_{3},\lambda}(S(P),S(Q)) 
\nonumber\\ 
& & \hspace{-0.2cm} 
= \int_{{\mathcal{X}}} 
\Bigg[ \phi \negthinspace \left( {\frac{S_{x}(P)}{\tilde{S}_{x}(P)}}\right) -
\phi \negthinspace \left( {\frac{S_{x}(Q)}{\tilde{S}_{x}(Q)}}\right) 
- \phi_{+,c}^{\prime} \negthinspace
\left( {\frac{S_{x}(Q)}{\tilde{S}_{x}(Q)}}\right) \cdot 
\left( \frac{S_{x}(P)}{\tilde{S}_{x}(P)}-\frac{S_{x}(Q)}{\tilde{S}_{x}(Q)}\right) 
\Bigg] 
m_{3}(x) \, \mathrm{d}\lambda(x) ,
\nonumber \\ 
\label{BroStuHB22:fo.def.55a}
\end{eqnarray}
which for the discrete setup 
$(\mathcal{X},\lambda) = (\mathcal{X}_{\#},\lambda_{\#})$ 
simplifies to
\begin{eqnarray} 
& & \hspace{-0.2cm} \textstyle
0 \leq D_{\phi,\tilde{S}(P),\tilde{S}(Q),m_{3},\lambda_{\#}}(S(P),S(Q)) 
\nonumber\\ 
& & \hspace{-0.2cm} 
 = \sum_{{x \in \mathcal{X}}} 
\Bigg[ \phi \negthinspace \left( {
\frac{S_{x}(P)}{\tilde{S}_{x}(P)}}\right) -
\phi \negthinspace \left( {\frac{S_{x}(Q)}{\tilde{S}_{x}(Q)}}\right)  
- \phi_{+,c}^{\prime} \negthinspace
\left( {\frac{S_{x}(Q)}{\tilde{S}_{x}(Q)}}\right) \cdot 
\left( \frac{S_{x}(P)}{\tilde{S}_{x}(P)}-\frac{S_{x}(Q)}{\tilde{S}_{x}(Q)}\right) 
\Bigg] 
m_{3}(x) \, .
\nonumber 
\end{eqnarray}

\noindent
As an example, take 
$\mathcal{Y}= \mathcal{X} = 
[0,\infty[$, $\lambda = \lambda_{L}$,
the probability (Lebesgue-) density functions $S=S^{pd}$, i.e.
$S(P) = \left\{S_{x}(P)\right\}_{x \in [0,\infty[}
= \left\{f_{P}(x)\right\}_{x \in [0,\infty[}
= \left\{\frac{\mathrm{d} F_{P}(x)}{\mathrm{d} x}\right\}_{x \in [0,\infty[}$,
as well as the survival (reliability, tail) functions $\tilde{S}= S^{su}$,\\
\noindent i.e.
$\tilde{S}(P)
= \left\{\tilde{S}_{x}(P)\right\}_{x \in [0,\infty[} 
= \left\{1-F_{P}(x)\right\}_{x \in [0,\infty[}
= \left\{P[(x,\infty)]\right\}_{x \in [0,\infty[}$.
Accordingly, the function 
$x \rightarrow \frac{S_{x}(P)}{\tilde{S}_{x}(P)} =
\frac{f_{P}(x)}{1-F_{P}(x)}
$
-- with the convention $\frac{c}{0} = \infty$ for all $c \in \mathbb{R}$ -- 
can be interpreted as the hazard rate function 
(failure rate function, 
force of mortality)
under the model distribution $P$ (and analogously under the alternative model distribution $Q$)
of a nonnegative random variable $Y$. Hence,
\eqref{BroStuHB22:fo.def.55a} turns into
\begin{eqnarray} 
& & \hspace{-0.2cm} \textstyle
0 \leq D_{\phi,S^{su}(P),S^{su}(Q),m_{3},\lambda_{L}}\left(S^{pd}(P),S^{pd}(Q) \right) 
\nonumber\\ 
& & \hspace{-0.2cm} 
 =
\int_{{\mathcal{X}}} 
\Bigg[ \phi \negthinspace \left( {
\frac{f_{P}(x)}{1-F_{P}(x)}}\right) -\phi \negthinspace \left( {\frac{f_{Q}(x)}{1-F_{Q}(x)}}\right)
\nonumber\\ 
& & \hspace{-0.2cm}  
-
\phi_{+,c}^{\prime} \negthinspace
\left( {\frac{f_{Q}(x)}{1-F_{Q}(x)}}\right) \cdot \left(\frac{f_{P}(x)}{1-F_{P}(x)} - \frac{f_{Q}(x)}{1-F_{Q}(x)}\right) 
\Bigg] 
m_{3}(x) \, \mathrm{d}\lambda_{L}(x) ,
\nonumber
\end{eqnarray}
which can be interpreted as divergence between the two modeling hazard rate functions at 
stake.

%
\subsection{Auto-Divergences}
%

The main-stream of this paper deals with divergences/distances
between (families of) real-valued ``statistical
functionals'' $S(\cdot)$ of the form
$S(P) := \left\{S_{x}(P)\right\}_{x \in \mathcal{X}}$
and $S(Q) := \left\{S_{x}(Q)\right\}_{x \in \mathcal{X}}$
stemming from two different distributions $P$ and $Q$.
In quite some meaningful situations, $P$ and $Q$ can stem from the \textit{same}
fundamental underlying random mechanism $\breve{P}$. Take for instance
the situation where $\mathcal{Y} = \mathcal{X} = \mathbb{R}$, $\lambda = \lambda_{L}$
 and $Y_{1}, \ldots Y_{N}$ are i.i.d.
observations from a random variable $Y$ with distribution $\breve{P}$ having 
(with a slight abuse of notation 
$\breve{P} = \breve{P} \circ Y^{-1}$)
distribution function $F_{\breve{P}}(x) = \breve{P}[Y \leq x]$ which is differentiable with a density
$f_{\breve{P}}(x) = \frac{\mathrm{d}F_{\breve{P}}(x)}{\mathrm{d}x}$
being positive in an interval and zero elsewhere.
The corresponding order statistics are denoted by 
$Y_{1:N} < Y_{2:N} < \ldots < Y_{N:N}$
where $Y_{k:N}$ is the $k-$th largest observation
and in particular
$Y_{1:N} := \min \{Y_{1}, \ldots  Y_{N} \}$,
$Y_{N:N} := \max \{Y_{1}, \ldots  Y_{N} \}$;
the distribution $\breve{P}_{k}$ of $Y_{k:N}$ \, ($k\in \{1,\ldots,N\}$) has
distribution function $F_{\breve{P}_{k}}(x) := \breve{P}[Y_{k:N} \leq x]$ 
with well-known density function
\begin{eqnarray} 
& & \hspace{-0.2cm} 
f_{\breve{P}_{k}}(x) := 
\frac{N!}{(N-k)! \cdot (k-1)!} 
\cdot
\left( F_{\breve{P}}(x) \right)^{k-1} \cdot 
\left( 1- F_{\breve{P}}(x) \right)^{n-k} \cdot
f_{\breve{P}}(x) \, .
\label{BroStuHB22:fo.def.381a}
\end{eqnarray}
(see e.g. 
Reiss \cite{Reiss:89}, Arnold et al. \cite{Arn:08},
David \& Nagaraja \cite{David:03}
for comprehensive treatments of order statistics).
In such a context, it makes sense to 
take $P:= \breve{P}_{j}$, $Q:= \breve{P}_{k}$ ($j,k\in \{1,\ldots,N\}$)
respectively $P:= \breve{P}$, $Q:= \breve{P}_{k}$ (or vice versa)
and study the divergences
\begin{eqnarray} 
& & \hspace{-0.2cm} \textstyle
0 \leq D_{\phi,m_{1},m_{2},m_{3},\lambda_{L}}\left(S^{pd}(\breve{P}_{j}),S^{pd}(\breve{P}_{k})\right) 
\nonumber\\ 
& & \hspace{-0.2cm} 
: = \int_{{\mathcal{X}}} 
\Bigg[ \phi \negthinspace \left( {\frac{f_{\breve{P}_{j}}(x)}{m_{1}(x)}}\right) 
-\phi \negthinspace \left( {\frac{f_{\breve{P}_{k}}(x)}{m_{2}(x)}}\right)
- \phi_{+,c}^{\prime} \negthinspace
\left( {\frac{f_{\breve{P}_{k}}(x)}{m_{2}(x)}}\right) \cdot \left( \frac{f_{\breve{P}_{j}}(x)}{m_{1}(x)}
-\frac{f_{\breve{P}_{k}}(x)}{m_{2}(x)}\right) 
\Bigg] 
m_{3}(x) \, \mathrm{d}\lambda_{L}(x)
\nonumber 
\end{eqnarray}
respectively
\begin{eqnarray} 
& & \hspace{-0.2cm} \textstyle
0 \leq D_{\phi,m_{1},m_{2},m_{3},\lambda_{L}}\left(S^{pd}(\breve{P}),S^{pd}(\breve{P}_{k})\right) 
\nonumber\\ 
& & \hspace{-0.2cm} 
: = \int_{{\mathcal{X}}} 
\Bigg[ \phi \negthinspace \left( {\frac{f_{\breve{P}}(x)}{m_{1}(x)}}\right) -
\phi \negthinspace \left( {\frac{f_{\breve{P}_{k}}(x)}{m_{2}(x)}}\right)
- \phi_{+,c}^{\prime} \negthinspace
\left( {\frac{f_{\breve{P}_{k}}(x)}{m_{2}(x)}}\right) \cdot 
\left( \frac{f_{\breve{P}}(x)}{m_{1}(x)}
-\frac{f_{\breve{P}_{k}}(x)}{m_{2}(x)}\right) 
\Bigg] 
m_{3}(x) \, \mathrm{d}\lambda_{L}(x) \, ,
\nonumber \\ 
\label{BroStuHB22:fo.def.381c}
\end{eqnarray}
or deterministic transformations thereof.

\vspace{0.2cm}
For instance, (some of) the divergences in 
Ebrahimi et al.~\cite{Ebr:04}, Asadi et al.~\cite{Asa:06}
can be imbedded here as the special cases
$D_{\phi_{1},1,1,1,\lambda_{L}}\left(S^{pd}(\breve{P}_{j}),S^{pd}(\breve{P}_{k})\right)$,
$D_{\phi_{1},1,1,1,\lambda_{L}}\left(S^{pd}(\breve{P}),S^{pd}(\breve{P}_{k})\right)$,\\
$\frac{1}{\alpha-1} \log \left[ 1 + \alpha \cdot (\alpha-1) \cdot
D_{\phi_{\alpha},S^{pd}(\breve{P}_{k}),S^{pd}(\breve{P}_{k}),S^{pd}(\breve{P}_{k}),\lambda}
\left(S^{pd}(\breve{P}_{j}),S^{pd}(\breve{P}_{k})\right)
\right]
$,\\
$\frac{1}{\alpha-1} \log \left[ 1 + \alpha \cdot (\alpha-1) \cdot
D_{\phi_{\alpha},S^{pd}(\breve{P}_{k}),S^{pd}(\breve{P}_{k}),S^{pd}(\breve{P}_{k}),\lambda}
\left(S^{pd}(\breve{P}),S^{pd}(\breve{P}_{k})\right)
\right]
$,\\
\noindent
for $\alpha \in \mathbb{R}\backslash\{0,1\}$.

\vspace{0.3cm}
\noindent
For other (non-auto type) scaled Bregman divergences involving distributions 
of certain transforms of spacings
between observations (i.e., differences of order statistics), the reader is e.g. referred 
to Roensch \& Stummer \cite{Roe:19b}.

\vspace{0.3cm}
\noindent
Vaughan \& Venables \cite{Vau:72}, Bapat \& Beg \cite{Bap:89} and Hande \cite{Hande:94} 
give some extensions of \eqref{BroStuHB22:fo.def.381a} for random observations $Y_{1}, \ldots Y_{N}$ 
which are independent but non-identically distributed, e.g. 
their distributions may be linked by a common (scalar or multidimensional) parameter;
this is a common situation in contemporary statistical applications
e.g. in data analytics, artificial intelligence and machine learning (which employ GLM models, etc.). 
By employing \eqref{BroStuHB22:fo.def.381c} for these extensions of \eqref{BroStuHB22:fo.def.381a},
we end up with an even wider new toolkit for auto-divergences between (distributions of) order statistics.

%
\subsection{Connections with optimal transport and coupling}
\label{subsec.2.7new}
%

In this section we consider the context of Subsection 2.5.1.3 with $\mathcal{X}= \, ]0,1[$,
Lebesgue measure $\lambda=\lambda_{L}$ as well $r(x)=1$ for all $x \in \mathcal{X}$, 
and apply this
to the quantile functions 
$S^{qu}(P) = 
\left\{S_{x}(P)\right\}_{x \in ]0,1[} 
:= \left\{F_{P}^{\leftarrow}(x)\right\}_{x \in ]0,1[}$\\
$:= \left\{\inf\{z \in \mathbb{R}: F_{P}(z) \geq x \} \right\}_{x \in ]0,1[}$ 
respectively $S^{qu}(Q)$ of two random variables $X$ respectively $Y$ on $\mathcal{Y}=\mathbb{R}$
having distribution $P$ respectively $Q$; recall from Section \ref{subsec.2.0averynew}
that for $\mathcal{Y}=[0,\infty)$ we take
$S^{qu}(P) = 
\left\{S_{x}(P)\right\}_{x \in ]0,1[} 
:= \left\{F_{P}^{\leftarrow}(x)\right\}_{x \in ]0,1[}
:= \left\{\inf\{z \in [0,\infty): F_{P}(z) \geq x \} \right\}_{x \in ]0,1[}
$ instead. Accordingly, we quantify the corresponding dissimilarity
as the divergence (directed distance)
\begin{eqnarray} 
& & \hspace{-0.2cm} \textstyle
D^{c}_{\phi,w(S^{qu}(P),S^{qu}(Q)),w(S^{qu}(P),S^{qu}(Q)),
w(S^{qu}(P),S^{qu}(Q)),\lambda_{L}}(S^{qu}(P),S^{qu}(Q)) 
\nonumber\\ 
& & \hspace{-0.2cm} \textstyle
: =
\olint_{]0,1[} 
\Bigg[ \phi \negthinspace \left( {\frac{F_{P}^{\leftarrow}(x)}{w(F_{P}^{\leftarrow}(x),F_{Q}^{\leftarrow}(x))}}\right) 
-\phi \negthinspace \left( {\frac{F_{Q}^{\leftarrow}(x)}{w(F_{P}^{\leftarrow}(x),F_{Q}^{\leftarrow}(x))}}\right)
\nonumber\\ 
& & \hspace{-0.2cm} \textstyle 
- \phi_{+,c}^{\prime} \negthinspace
\left( {\frac{F_{Q}^{\leftarrow}(x)}{w(F_{P}^{\leftarrow}(x),F_{Q}^{\leftarrow}(x))}}\right) \cdot \left( \frac{F_{P}^{\leftarrow}(x)}{w(F_{P}^{\leftarrow}(x),F_{Q}^{\leftarrow}(x))}
-\frac{F_{Q}^{\leftarrow}(x)}{w(F_{P}^{\leftarrow}(x),F_{Q}^{\leftarrow}(x))}\right)  
\Bigg] \nonumber\\ 
& & \hspace{-0.2cm} \textstyle
\cdot
w(F_{P}^{\leftarrow}(x),F_{Q}^{\leftarrow}(x)) \cdot r(x) \, \mathrm{d}\lambda_{L}(x) \ 
\label{BroStuHB22:fo.def.347ot.new}
\\ 
& & \hspace{-0.2cm} \textstyle
= \int_{]0,1[} 
\overline{\widetilde{\psi}}(F_{P}^{\leftarrow}(x),F_{Q}^{\leftarrow}(x)) \, \mathrm{d}\lambda_{L}(x)
\nonumber
\end{eqnarray}
with \noindent $\overline{\widetilde{\psi}}: 
\mathcal{R}\big(F_{P}^{\leftarrow}\big) \times \mathcal{R}\big(F_{Q}^{\leftarrow}\big) 
\mapsto [0,\infty]$
defined by (cf. (I2) and \eqref{BroStuHB22:fo.def.3})

\vspace{-0.5cm}

\begin{eqnarray}
& & \hspace{-0.2cm} 
\textstyle
\overline{\widetilde{\psi}}(u,v) := 
W(u,v) \cdot
\overline{\psi_{\phi,c}}
\Big(\frac{u}{W\big(u,v\big)},
\frac{v}{W\big(u,v\big)}\Big) \geq 0 \qquad \textrm{with} 
\nonumber
\\[-0.1cm]
& & \hspace{-0.2cm} 
\textstyle
\psi_{\phi,c}
\Big(\frac{u}{W\big(u,v\big)},
\frac{v}{W\big(u,v\big)}\Big) \negthinspace := \negthinspace
\Big[ \phi \negthinspace \big( {\frac{u}{W(u,v)}}\big) \negthinspace - \negthinspace
\phi \negthinspace \big( {\frac{v}{W(u,v)}}\big) 
\negthinspace - \negthinspace \phi_{+,c}^{\prime} \negthinspace
\big( {\frac{v}{W(u,v)}}\big) \cdot \big( {\frac{u}{W(u,v)}}\negthinspace 
 - \negthinspace{\frac{v}{W(u,v)}} \big)
\Big]. \qquad \ 
\nonumber
\end{eqnarray}
Under Assumption 1 (and hence under the more restrictive Assumption 2) of Stummer~\cite{Stu:21}
-- who deals even with a more general context where the scaling and the aggregation function
need not coincide --
one can adapt Theorem 4 and Corollary 1 of Broniatwoski \& Stummer~\cite{Bro:19a}
to obtain the desired basic divergence properties (D1) and (D2) in the form of 
\begin{eqnarray}
& & \hspace{-0.2cm}
(NN) D^{c}_{\phi,w(S^{qu}(P),S^{qu}(Q)),w(S^{qu}(P),S^{qu}(Q)),
w(S^{qu}(P),S^{qu}(Q)),\lambda_{L}}(S^{qu}(P),S^{qu}(Q)) \geq 0
\nonumber \\
& & \hspace{-0.2cm}
\textstyle
\textit{(RE)} \  D^{c}_{\phi,w(S^{qu}(P),S^{qu}(Q)),w(S^{qu}(P),S^{qu}(Q)),
w(S^{qu}(P),S^{qu}(Q)),\lambda_{L}}(S^{qu}(P),S^{qu}(Q)) = 0 
\nonumber \\
& & \hspace{4.2cm}
\textstyle 
\textrm{if and only if} \ \  
F_{P}^{\leftarrow}(x) =
F_{Q}^{\leftarrow}(x)
\ \textrm{for $\lambda$-a.a. $x \in \mathcal{X}$.
} 
\qquad \ 
\nonumber 
\end{eqnarray} 

\noindent In order to establish a connection between the divergence \eqref{BroStuHB22:fo.def.347ot.new}
and optimal transport problems, we impose for the rest of this section 
the additional requirement that the 
function $\overline{\widetilde{\psi}}$
is continuous (except for the point $(u,v)=(0,0)$) and quasi-antitone\footnote{other names are: 
submodular, Lattice-subadditive, 2-antitone, 2-negative, $\Delta-$antitone, supernegative,
``satisfying the (continuous) Monge property/condition''}
in the sense
\begin{eqnarray}
& & \hspace{-0.2cm} \textstyle
\overline{\widetilde{\psi}}(u_{1},v_{1}) + \overline{\widetilde{\psi}}(u_{2},v_{2}) 
\leq \overline{\widetilde{\psi}}(u_{2},v_{1}) + \overline{\widetilde{\psi}}(u_{1},v_{2})
\quad \textrm{for all $u_{1} \leq u_{2}$, $v_{1} \leq v_{2}$}; \qquad \ 
\nonumber
\end{eqnarray}
in other words, $-\overline{\widetilde{\psi}}(\cdot,\cdot)$
is assumed to be continuous (except for the point $(u,v)=(0,0)$)  
and quasi-monotone\footnote{other names are: 
supermodular, Lattice-superadditive, 2-increasing, 2-positive, $\Delta-$monotone, 
2-monotone, 
``fulfilling the moderate growth property'', ``satisfying the measure property'',
``satisfying the twist condition''}
\footnote{
a comprehensive discussion on general quasi-monotone functions can be found e.g. in Chapter 6.C of
Marshall et al.~\cite{Mar:11}}. 
For such a setup, one can consider the 
Kantorovich transportation problem (KTP) 
with the \textit{pointwise-BS-distance-type} (pBS-type) cost function 
$(u,v) \mapsto \overline{\widetilde{\psi}}(u,v)$;
indeed,  Stummer~\cite{Stu:21} recently obtained (an even more general version of) the following 

\vspace{-0.2cm}

\begin{theorem}
\label{Stu:thm1}
Let $\widetilde{\Gamma}(P,Q)$ be
the family of all probability distributions $\mathfrak{P}$ on 
$\mathbb{R} \times \mathbb{R}$
which have marginal distributions $\mathfrak{P}[\, \cdot \times 
\mathbb{R} ] = P[\cdot]$ 
and $\mathfrak{P}[ \mathbb{R} \times \cdot \, ] = Q[\, \cdot \, ]$. 
Moreover, we denote
the corresponding upper Hoeffding-Fr\'echet bound (cf. e.g. Theorem 3.1.1 of
Rachev \& R\"uschendorf~\cite{Rac:98a}) by
$\mathfrak{P}^{com}$ having ``comonotonic'' 
distribution function $F_{\mathfrak{P}^{com}}(u,v) := \min\{F_{P}(u),F_{Q}(v)\}$
($u,v \in \mathbb{R}$). Then
\begin{eqnarray}
& & \hspace{-0.2cm} \textstyle
\min_{\{ X \sim P, \, Y \sim Q \}} \  \mathbb{E}\big[ \    
\overline{\widetilde{\psi}}(X,Y) \, \big] 
\label{BroStuHB22:fo.def.351ot}
\\ 
& & \hspace{-0.2cm} \textstyle 
= \min_{\{ \mathfrak{P} \in \widetilde{\Gamma}(P,Q) \}} \  
\int\displaylimits_{\mathbb{R} \times \mathbb{R}}
\overline{\widetilde{\psi}} \, \mathrm{d}\mathfrak{P}(u,v)
\label{BroStuHB22:fo.def.350ot}
\\ 
& & \hspace{-0.2cm} \textstyle
= \int\displaylimits_{\mathbb{R} \times \mathbb{R}}
\overline{\widetilde{\psi}}(u,v) \, \mathrm{d}\mathfrak{P}^{com}(u,v)
\label{BroStuHB22:fo.def.349ot}
\\ 
& & \hspace{-0.2cm} \textstyle
=  \int_{]0,1[}
\overline{\widetilde{\psi}}(F_{P}^{\leftarrow}(x),F_{Q}^{\leftarrow}(x)) \, \mathrm{d}\lambda_{L}(x)
\nonumber
\\ 
& & \hspace{-0.2cm} \textstyle
= D^{c}_{\phi,w(S^{qu}(P),S^{qu}(Q)),w(S^{qu}(P),S^{qu}(Q)),
w(S^{qu}(P),S^{qu}(Q)),\lambda_{L}}(S^{qu}(P),S^{qu}(Q)) 
\geq 0 , \qquad  \ 
\label{BroStuHB22:fo.def.XXXot}
\end{eqnarray}

\vspace{-0.2cm}

\noindent
where the minimum in \eqref{BroStuHB22:fo.def.351ot} is taken over
all $\mathbb{R}-$valued
random variables $X$, $Y$ (on an arbitrary probability space $(\Omega, \mathcal{A}, \mathfrak{S})$)
such that $\mathfrak{P}[X \in \cdot \, ] = P[ \, \cdot \, ]$,
$\mathfrak{P}[Y \in \cdot \, ] = Q[ \, \cdot \, ]$. As usual, 
$\mathbb{E}$ denotes the expectation with respect to $\mathfrak{P}$.
\end{theorem}

\begin{remark}
\label{Stu.rem1}
(i) \, Notice that $\mathfrak{P}^{com}$
is $\overline{\widetilde{\psi}}-$independent, and may not be 
the unique minimizer in \eqref{BroStuHB22:fo.def.350ot}. 
As a (not necessarily unique) minimizer
in \eqref{BroStuHB22:fo.def.351ot}, one can take $X := F_{P}^{\leftarrow}(U)$,
$Y := F_{Q}^{\leftarrow}(U)$ for some 
uniform random variable $U$ on $[0,1]$.\\
(ii) \, In Theorem \ref{Stu:thm1}
we have shown that 
$\mathfrak{P}^{com}$ (cf. \eqref{BroStuHB22:fo.def.349ot}) is an optimal transport plan of
the KTP \eqref{BroStuHB22:fo.def.350ot}
with the \textit{pointwise-BS-distance-type} (pBS-type)
cost function 
$\overline{\widetilde{\psi}}(u,v)$. The outcoming minimal value
is equal to \\
$D^{c}_{\phi,w(S^{qu}(P),S^{qu}(Q)),w(S^{qu}(P),S^{qu}(Q)),
w(S^{qu}(P),S^{qu}(Q)),\lambda_{L}}(S^{qu}(P),S^{qu}(Q))$
which is typically straightforward to compute (resp. approximate).\\
(iii) Depending on the chosen divergence, one may have to restrict the support
of $P$ respectively $Q$, for instance to (subsets of) $[0,\infty[$.

\end{remark}

\vspace{-0.1cm}

\noindent
Remark \ref{Stu.rem1}(ii) generally contrasts to those prominently used 
KTP whose cost function is a power
$d(u,v)^{p}$ of 
a metric $d(u,v)$ (denoted as POM-type cost function) which leads to the well-known Wasserstein distances.
(Apart from technicalities) 
There are some overlaps, though:

\vspace{-0.1cm}

\begin{example}
(i) Take $\mathcal{Y} \subset [0,\infty)$ (and thus the support of $P$ and $Q$
is contained in $[0,\infty[$)
together with the \textit{non-smooth} $\phi(t) := \phi_{TV}(t):= |t-1|$
($t \in [0,\infty[$), $c=\frac{1}{2}$,  
$W(u,v) := v \in [0,\infty[$
to obtain
$\overline{\widetilde{\psi}}(u,v) = |u-v| =: d(u,v)$
($u,v \in [0,\infty[$). For an extension to $\mathcal{Y}=\mathbb{R}$
see Stummer~\cite{Stu:21}.
\\
(ii) Take $\mathcal{Y}=\mathbb{R}$, $\phi(t) := \phi_{2}(t) := \frac{(t-1)^2}{2}$
($t \in \mathbb{R}$, with obsolete $c$),
$W(u,v) :=1$ to end up with  
$\overline{\widetilde{\psi}}(u,v) = \frac{(u-v)^2}{2}
= \frac{d(u,v)^2}{2}$.\\
(iii) The \textit{symmetric} distances $d(u,v)$ and $\frac{d(u,v)^2}{2}$ 
are convex functions of $u-v$ and thus continuous quasi-antitone
functions. 
The correspondingly outcoming Wasserstein distances are thus considerably flexibilized
by our new much more general distance 
$D^{c}_{\phi,w(S^{qu}(P),S^{qu}(Q)),w(S^{qu}(P),S^{qu}(Q)),
w(S^{qu}(P),S^{qu}(Q)),\lambda_{L}}(S^{qu}(P),S^{qu}(Q))$
of \eqref{BroStuHB22:fo.def.XXXot}.
\end{example} 

\enlargethispage{0.5cm}

\noindent
We give some further special cases of pBS-type cost functions, 
which are continuous and quasi-antitone, 
but which are generally not symmetric and
thus not of POM-type:

\vspace{-0.1cm}

\begin{example}
``smooth'' pointwise \textit{Csiszar-Ali-Silvey-Morimoto divergences} 
(CASM divergences):
take $\phi :[0,\infty[ \mapsto \mathbb{R}$ to be a strictly convex, twice continuously differentiable
function on $]0,\infty[$ with continuous extension on $t=0$, 
together with $W(u,v) := v \in ]0,\infty[$, 
and $c$ is obsolete. Accordingly, $\overline{\widetilde{\psi}}(u,v)  = 
v \cdot \phi \negthinspace \left( \frac{u}{v}\right) - v \cdot \phi \negthinspace \left( 1 \right)  
- \phi^{\prime} \negthinspace
\left( 1 \right) \cdot \left( u - v \right) $,
and hence the second mixed derivative satisfies 
$\frac{\partial^{2} \overline{\widetilde{\psi}}}{\partial u \partial v}
= - \frac{u}{v^2} \phi^{\prime\prime} \negthinspace \left( \frac{u}{v}\right) < 0$ 
($u,v \in ]0,\infty[$); thus, $\overline{\widetilde{\psi}}$ is
quasi-antitone on $]0,\infty[ \times ]0,\infty[$. 
Accordingly,  
\eqref{BroStuHB22:fo.def.351ot} to 
\eqref{BroStuHB22:fo.def.XXXot}
applies to such kind of (cf. Section 2.5.1.2) CASM divergences 
concerning $P$,$Q$ having support in $[0,\infty[$.
As an example, take e.g. the power function 
$\phi(t):= \frac{t^\gamma-\gamma \cdot t+ \gamma - 1}{\gamma \cdot (\gamma-1)}$
($\gamma \in \mathbb{R}\backslash\{0,1\}$).
A different connection between optimal transport and 
other kind of CASM divergences can be found in Bertrand et al.~\cite{Ber:21}.
\end{example}

\begin{example}
``smooth'' pointwise \textit{classical (i.e. unscaled) Bregman divergences} (CBD): 
take $\phi : \mathbb{R} \mapsto \mathbb{R}$ to be a strictly convex, 
twice continuously differentiable function
$W(u,v) := 1$
and $c$ is obsolete. Accordingly, $\overline{\widetilde{\psi}}(u,v) := 
\phi \negthinspace \left( u\right) -\phi \negthinspace \left( v \right)
- \phi^{\prime} \negthinspace
\left( v \right) \cdot \left( u - v \right) $
and hence $\frac{\partial^{2} \overline{\widetilde{\psi}}(u,v)}{\partial u \partial v}
= - \phi^{\prime\prime} \negthinspace \left( v\right) < 0$ 
($u,v \in \mathbb{R}$); thus, $\Upsilon_{\phi,c,W,W_{3}}$ is
quasi-antitone on $\mathbb{R} \times \mathbb{R}$.
Accordingly, the representation 
\eqref{BroStuHB22:fo.def.351ot} to 
\eqref{BroStuHB22:fo.def.XXXot} 
applies to such kind of (cf. Section 2.5.1.1) CBD. The corresponding special case of
\eqref{BroStuHB22:fo.def.350ot} is called 
``a relaxed Wasserstein distance (parameterized by $\phi$) between
$P$ and $Q$'' in the recent papers of Lin et al.~\cite{Lin:19} and Guo et al.~\cite{Guo:21}
for a \textit{restrictive} setup where $P$ and $Q$ are supposed to have 
\textit{compact} support;
the latter two references do 
not give connections to divergences of quantile functions,
but substantially concentrate on applications to topic sparsity
for analyzing user-generated web content and social media, 
respectively, to Generative Adversarial Networks (GANs).
\end{example}

\vspace{-0.4cm}

\begin{example}
``smooth'' pointwise \textit{Scaled Bregman Distances}: for instance, 
consider $P$ and $Q$ with support in $[0,\infty[$.
One gets that $\overline{\widetilde{\psi}}$ is quasi-antitone on 
$]0,\infty[ \times ]0,\infty[$
if the generator function $\phi$ is strictly convex and
thrice continuously differentiable on $]0,\infty[$ (and hence, c is obsolete) and 
the so-called scale connector $W$ is twice continuously differentiable such that --
on $]0,\infty[ \times ]0,\infty[$ -- 
$\overline{\widetilde{\psi}}$ is twice continuously differentiable 
and $\frac{\partial^2 \overline{\widetilde{\psi}}}{\partial u \partial v} \leq 0$
(an explicit formula of the latter is given in the appendix of
Ki{\ss}linger \& Stummer~\cite{Kis:18}, who also
give applications to robust change detection in data streams). 
Illustrative examples of suitable $\phi$ and $W$ 
can be found e.g. in Ki{\ss}linger \& Stummer \cite{Kis:16}.
\end{example}

\noindent
Returning to the general context,  
it is straightforward to see that if $P$ does not give mass to points
(i.e. it has continuous distribution function $F_{P}$) then there exists even a deterministic
optimal transportation plan: indeed, for the map $T^{com} :=  F_{Q}^{\leftarrow} \circ F_{P}$ 
one has $\mathfrak{P}^{com}[ \, \cdot \, ] = P[(id,T^{com}) \in \cdot \, ]$
and thus \eqref{BroStuHB22:fo.def.349ot} is equal to

\vspace{-0.3cm}

\begin{eqnarray}
& & \hspace{-0.2cm} \textstyle
\int\displaylimits_{\mathbb{R}} 
\overline{\widetilde{\psi}}(u,T^{com}(u)) \, \mathrm{d}P(u)
\nonumber
\\  
& & \hspace{-0.2cm} \textstyle
= \min_{\{ T \in \widehat{\Gamma}(P,Q) \}} \ 
\int\displaylimits_{\mathbb{R}} 
\overline{\widetilde{\psi}}(u,T(u)) \, \mathrm{d}P(u)
\label{BroStuHB22:fo.def.353ot}
\\ 
& & \hspace{-0.2cm} \textstyle
= \min_{\{ X \sim P, \, T(X) \sim Q \}} \  \mathbb{E}\big[ \    
\overline{\widetilde{\psi}}(X,T(X)) \, \big]
\nonumber
\end{eqnarray}

\vspace{-0.3cm}

\noindent
where \eqref{BroStuHB22:fo.def.353ot} is called Monge transportation problem (MTP).
Here, $\widehat{\Gamma}(P,Q)$ denotes the family of all measurable maps
$T: \mathbb{R} \mapsto \mathbb{R}$
such that
$P[ T \in \cdot \, ] = Q[\, \cdot \, ]$.

%
\section{Aggregated/Integrated Divergences}
%

Suppose that $\phi = \phi_{z}$, $P =P_{z}$, $Q=Q_{z}$, $m_{1} = m_{1,z}$ 
$m_{2} = m_{2,z}$, $m_{3} = m_{3,z}$, $\lambda = \lambda_{z}$ depend on the same (!!) 
``parameter/quantity''  $z\in \mathcal{Z}$.
Then it makes sense to study
the aggregated/integrated divergence  
$\int_{\mathcal{Z}}  D_{\phi_{z},m_{1,z},m_{2,z},m_{3,z},\lambda_{z}}(S(P_{z}),S(Q_{z})) \, 
\mathrm{d}\breve{\lambda}(z)$
where $\breve{\lambda}$ is a $\sigma-$finite measure on $\mathcal{Z}$
(e.g. the Lebesgue-measure $\lambda_{L}$, the counting measure $\lambda_{\#}$
or a probability measure, where in case of the latter one also uses the terminology
``expected divergence'').\\[0.1cm]

\noindent
Another interesting special case is the following family:
recall first that for the two-element space $\mathcal{Y}=\mathcal{X} = \{0,1 \}$   
we denote the corresponding probability mass functions as
$S^{pm}(P) = \left\{P[\{x\}]\right\}_{x \in \mathcal{X}} =
 \left\{1-P[\{1\}],P[\{1\}]\right\} $;
in other words, $P$ is a Bernoulli distribution $Ber(\theta)$
which is completely determined by its parameter $\theta  \in [0,1]$ with interpretation
$\theta_{P} = P[\{1\}]$. Now suppose that $\theta_{P} = \theta_{P}(z)$ depends on a real-valued
parameter $z\in \mathbb{R}$. In such a situation it makes sense to study the 
the aggregated (integrated) divergence for $\phi \in \Phi_{C_{1}}(]a,b[)$
\begin{eqnarray} 
& & \hspace{-0.2cm} 
\textstyle
0 \leq \olint_{\mathbb{R}} D_{\phi,m_{1,z},m_{2,z},m_{3,z},\lambda_{\#}}(S^{pm}(Ber(\theta_{P}(z))),S^{pm}(Ber(\theta_{Q}(z))))
\, \mathrm{d}\breve{\lambda}(z)
\nonumber\\
& & \hspace{-0.2cm} = \ \olint_{\mathbb{R}} 
\Bigg\{
\Bigg[ \phi \negthinspace \left( {\frac{1-\theta_{P}(z)}{m_{1,z}(0)}}\right) 
-\phi \negthinspace \left( 
{\frac{1-\theta_{Q}(z)}{m_{2,z}(0)}}\right)
- \phi^{\prime} \negthinspace
\left( {\frac{1-\theta_{Q}(z)}{m_{2,z}(0)}}\right) \cdot \left( 
\frac{1-\theta_{P}(z)}{m_{1,z}(0)}-\frac{1-\theta_{Q}(z)}{m_{2,z}(0)}\right) 
\Bigg] \cdot
m_{3,z}(0)
\nonumber\\
& & \hspace{-0.2cm} +
\Bigg[ \phi \negthinspace \left( {\frac{\theta_{P}(z)}{m_{1,z}(1)}}\right)
 -\phi \negthinspace \left( 
{\frac{\theta_{Q}(z)}{m_{2,z}(1)}}\right) 
- \phi^{\prime} \negthinspace
\left( {\frac{\theta_{Q}(z)}{m_{2,z}(1)}}\right) \cdot \left( 
\frac{\theta_{P}(z)}{m_{1,z}(1)}-\frac{\theta_{Q}(z)}{m_{2,z}(1)}\right) 
\Bigg] \cdot
m_{3,z}(1)\Bigg\} \, \mathrm{d}\breve{\lambda}(z)
\label{BroStuHB22:fo.def.500}
\end{eqnarray}  
where $\breve{\lambda}$ is a $\sigma-$finite measure on $\mathbb{R}$
(e.g. the Lebesgue-measure $\lambda_{L}$, the counting measure $\lambda_{\#}$
 or a probability measure) and the scaling functions $m_1$, $m_2$ as well
as the aggregating function $m_3$ are allowed to depend (in a measurable way) on $z$
(which is denoted by extending their indices with $z$).
For the non-differentiable case $\phi \in \Phi_{0}(]a,b[)$, the derivative 
$\phi^{\prime}$ has to be replaced 
by $\phi_{+,c}^{\prime}$.

\vspace{0.2cm}
In adaption of the discussion after formula \eqref{BroStuHB22:fo.def.17},
by defining the integral functional 
$\tilde{g}_{\phi,m_{3},\breve{\lambda}}(\tilde{\xi}) := \int_{\mathbb{R}} \Big[ \int_{\{0,1\}} 
\phi(\tilde{\xi}(x,z)) \cdot m_{3}(x) \, \mathrm{d}\lambda_{\#}(x) \Big] \, \mathrm{d}\breve{\lambda}(z)$
and plugging in e.g. 
\begin{eqnarray} 
& & \hspace{-0.2cm} \textstyle
\tilde{g}_{\phi,m_{3},\breve{\lambda}}\negthinspace \negthinspace \left( {\frac{
S^{pm}(Ber(\theta_{P}(\cdot)))
}{m_{1,\cdot}}}\right) 
= 
\int_{\mathbb{R}} \Big\{
\phi \negthinspace \left( {\frac{1-\theta_{P}(z)}{m_{1,z}(0)}}\right)  \cdot m_{3,z}(0)
+ \phi \negthinspace \left( {\frac{\theta_{P}(z)}{m_{1,z}(1)}}\right) \cdot
m_{3,z}(1)\Big\} \, \mathrm{d}\breve{\lambda}(z),
\nonumber \\
\label{BroStuHB22:fo.def.502}
\end{eqnarray} 
the divergence in \eqref{BroStuHB22:fo.def.500} can be (formally) interpreted as 
\begin{eqnarray} 
& & \hspace{-0.2cm} \textstyle
0 \leq 0 \leq \int_{\mathbb{R}} D_{\phi,m_{1,z},m_{2,z},m_{3,z},
\lambda_{\#}}(S^{pm}(Ber(\theta_{P}(z))),S^{pm}(Ber(\theta_{Q}(z))))
\, \mathrm{d}\breve{\lambda}(z) 
\nonumber\\ 
& & \hspace{-0.2cm} 
= \tilde{g}_{\phi,m_{3},\breve{\lambda}}\negthinspace \negthinspace \left( 
{\frac{S^{pm}(Ber(\theta_{P}(\cdot)))}{m_{1,\cdot}}}\right) 
- \tilde{g}_{\phi,m_{3},\breve{\lambda}}\negthinspace \negthinspace \left( 
{\frac{S^{pm}(Ber(\theta_{Q}(\cdot)))}{m_{2,\cdot}}}\right) \nonumber\\
& & \hspace{-0.2cm} \textstyle
- \tilde{g}_{\phi,m_{3},\breve{\lambda}}^{\prime} \negthinspace \negthinspace 
\left( {\frac{S^{pm}(Ber(\theta_{Q}(\cdot)))}{m_{2,\cdot}}},
{\frac{S^{pm}(Ber(\theta_{P}(\cdot)))}{m_{1,\cdot}}} - {\frac{S^{pm}(Ber(\theta_{Q}(\cdot)))}{m_{2,\cdot}}}\right) \ .
\nonumber
\end{eqnarray}

\vspace{0.2cm}
\noindent
As an important special case, take $\breve{\lambda} := \lambda_{L}$
(and we formally identify the Lebesgue-integral with the Riemann-integral over $\mathrm{d}z$),
$\theta_{P}(z) := F_{P}(z) = P[(-\infty,z]] = S_{z}^{cd}(P)$,
$\theta_{Q}(z) := F_{Q}(z) = Q[(-\infty,z]] = S_{z}^{cd}(Q)$,
$m_{1,z}(0)= m_{2,z}(0)= m_{3,z}(0) = 1- \theta_{Q}(z)$,
$m_{1,z}(1)= m_{2,z}(1)= m_{3,z}(1) = \theta_{Q}(z)$,
and accordingly \eqref{BroStuHB22:fo.def.500} simplifies to
\begin{eqnarray} 
& & \hspace{-0.2cm} 
\textstyle
0 \leq \olint_{\mathbb{R}} D_{\phi,m_{1,z},m_{2,z},m_{3,z},
\lambda}(S^{pm}(Ber(\theta_{P}(z))),S^{pm}(Ber(\theta_{Q}(z))))
\, \mathrm{d}\breve{\lambda}(z)
\nonumber\\
& & \hspace{-0.2cm} = \ \olint_{\mathbb{R}} 
\Bigg\{
\Bigg[ \phi \negthinspace \left( {\frac{1-F_{P}(z)}{1-F_{Q}(z)}}\right) 
-\phi \negthinspace \left( 1 \right)
- \phi^{\prime} \negthinspace
\left( 1 \right) \cdot \left( \frac{1-F_{P}(z)}{1-F_{Q}(z)}-1\right) 
\Bigg] \cdot
(1-F_{Q}(z))
\nonumber\\
& & \hspace{-0.2cm} +
\Bigg[ \phi \negthinspace \left( {
\frac{F_{P}(z)}{F_{Q}(z)}}\right) -\phi \negthinspace \left( 1 \right) 
- \phi^{\prime} \negthinspace
\left( 1 \right) \cdot \left( \frac{F_{P}(z)}{F_{Q}(z)}-1 \right) 
\Bigg] \cdot
F_{Q}(z) \Bigg\} \, 
\mathrm{d}z \nonumber\\
& & \hspace{-0.2cm} 
=: CPD_{\phi}(P,Q) \ ,
\nonumber
\end{eqnarray} 
which in case of $\phi(1)=0$ becomes
\begin{eqnarray} 
& & \hspace{-0.2cm} 
\textstyle
0 \leq CPD_{\phi}(P,Q)
\nonumber\\
& & \hspace{-0.2cm} = \ \olint_{\mathbb{R}} 
\Bigg\{
\phi \negthinspace \left( {
\frac{1-F_{P}(z)}{1-F_{Q}(z)}}\right) \cdot
(1-F_{Q}(z)) +
\phi \negthinspace \left( {
\frac{F_{P}(z)}{F_{Q}(z)}}\right)  \cdot
F_{Q}(z) \Bigg\} \, \mathrm{d}z \ .
\label{BroStuHB22:fo.def.515}
\end{eqnarray}

\noindent
If basically $\phi(0) = \phi(1)=0$ 
and $P$, $Q$ are generated by random variables,
say $P=Pr[X \in \cdot \, ]$, \, $Q=Pr[Y \in \cdot \, ]$
-- and thus $F_{P}(z) = Pr[X\leq z]$, $F_{Q}(z) = Pr[Y\leq z]$ -- 
then according to \eqref{BroStuHB22:fo.def.515} the $CPD_{\phi}(P,Q)$   
coincides with the cumulative paired $\phi-$divergence $CPD_{\phi}(X,Y)$ of Klein et al.~\cite{Kle:16};
the special case $CPD_{\phi_{\alpha}}(X,Y)$ with $\phi = \phi_{\alpha}$ from \eqref{BroStuHB22:fo.def.32}
was employed by Jager \& Wellner~\cite{Jag:07}.
Notice that without the assumption $\phi(1)=0=\phi^{\prime} \negthinspace
\left( 1 \right)$, the right-hand side of \eqref{BroStuHB22:fo.def.515}
may become negative and thus is not a divergence anymore.

\noindent
As a side remark, notice that in the ``unscaled setup'' 
$\breve{\lambda} := \lambda_{L}$,
$\theta_{P}(z) := F_{P}(z)$, 
$m_{1,z}(0)= m_{3,z}(0) =m_{1,z}(1)= m_{3,z}(1) = 1$,
the formula \eqref{BroStuHB22:fo.def.502} becomes
\begin{eqnarray} 
& & \hspace{-0.2cm} \textstyle
\tilde{g}_{\phi,m_{3},\breve{\lambda}}\negthinspace \negthinspace \left( {\frac{
S^{pm}(Ber(\theta_{P}(\cdot)))
}{m_{1,\cdot}}}\right) 
= \int_{\mathbb{R}} 
\Big\{ \phi \negthinspace \left( 1-F_{P}(z) \right) 
+ \phi \negthinspace \left( F_{P}(z) \right) \Big\} \, \mathrm{d}z
\nonumber 
\end{eqnarray} 
which corresponds to the \textit{cumulative $\phi-$entropy of $P$}
introduced by Klein et al.~\cite{Kle:16}.\\[0.2cm]


\section{Dependence expressing divergences}

Let the data take values in some product space 
$\mathcal{Y} = \bigtimes_{i=1}^{d}  \mathcal{Y}_{i}$
with product-$\sigma-$algebra
$\mathcal{A} = \bigotimes_{i=1}^{d}  \mathcal{A}_{i}$.
On this, we consider probability distributions $P$ 
having marginals $P_{i}$
determined by 
$P_{i}[A_{i}] := P[\mathcal{Y}_{1} \cdots \times A_{i} \cdots \times \mathcal{Y}_{d}]$
($i \in \{1,\ldots,d \}$, $A_i \in \mathcal{A}_{i}$).
Furthermore, let $Q := \bigotimes_{i=1}^{d} P_{i}$ be the 
product measure having the same marginals
as $P$.
Typically, $P[\cdot] := Pr[(Y_{1}, \ldots, Y_{d}) \in \cdot \, ]$ is the joint distribution 
of some random variables $Y_{1}, \ldots, Y_{d}$; on the other hand, the latter are independent
under the (generally different) probability measure $Q$. \\
As usual, we also involve statistical functionals $S(P) := \left\{S_{x}(P)\right\}_{x \in \mathcal{X}}$
and $S(Q) := \left\{S_{x}(Q)\right\}_{x \in \mathcal{X}}$,
where $\mathcal{X}$ is an index space equipped with a $\sigma-$algebra
$\mathcal{F}$ and a $\sigma-$finite measure $\lambda$
(e.g. a probability measure, the Lebesgue measure, a counting measure, etc.).
Accordingly, any of the above divergences (cf. \eqref{BroStuHB22:fo.def.1})
\begin{eqnarray} 
& & \hspace{-0.2cm} \textstyle
0 \leq D^{c}_{\phi,m_{1},m_{2},m_{3},\lambda}(S(P),S(Q)) 
\nonumber\\ 
& & \hspace{-0.2cm} 
: = \olint_{{\mathcal{X}}} 
\Bigg[ \phi \negthinspace \left( {\frac{S_{x}(P)}{m_{1}(x)}}\right) 
-\phi \negthinspace \left( 
{\frac{S_{x}(Q)}{m_{2}(x)}}\right) 
- \phi_{+,c}^{\prime} \negthinspace
\left( {\frac{S_{x}(Q)}{m_{2}(x)}}\right) \cdot \left( 
\frac{S_{x}(P)}{m_{1}(x)}-\frac{S_{x}(Q)}{m_{2}(x)}\right) 
\Bigg] 
m_{3}(x) \, \mathrm{d}\lambda(x)
\nonumber \\ 
\label{BroStuHB22:fo.def.1dep}
\end{eqnarray}
can be interpreted as a \textit{directed degree of dependence of $P$}
(e.g. of the above-mentioned random variables $Y_{1}, \ldots, Y_{d}$), 
since it measures the amount of
dissimilarity between the same statistical functional of $P$ 
and of the independence-expressing $Q$.
Some special cases of \eqref{BroStuHB22:fo.def.1dep} have already 
appeared in literature (which we put into our notation):\\

\noindent
(1) Micheas \& Zografos~\cite{Mic:06} consider 
Csiszar-Ali-Silvey-Morimoto (CASM) $\phi-$divergences
between $\lambda-$density functions, i.e. they take $\mathcal{X}: =\mathcal{Y}$,
a real continuous convex function on $[0,\infty[$,
a product measure $\lambda := \bigotimes_{i=1}^{d} \lambda_{i}$,
$S_{x}^{\lambda pd}(Q) := f_{Q}(x) := \prod_{i=1}^{d} f_{P_{i}}(x_{i}) \geq 0$
where $x=(x_{1},\ldots,x_{d})$ and $f_{P_{i}}$ is the $\lambda_{i}-$density 
function of the marginal distribution $P_{i}$,
as well as $S_{x}^{\lambda pd}(P) := f_{P}(x) \geq 0$ to be the $\lambda-$density
function of $P$,
to end up with the following special case of \eqref{BroStuHB22:fo.def.47d}:
\begin{eqnarray} 
& & \hspace{-0.2cm} 
0 \leq D_{\phi,S^{\lambda pd}(Q),S^{\lambda pd}(Q),1\cdot 
S^{\lambda pd}(Q),\lambda}(S^{\lambda pd}(P),S^{\lambda pd}(Q)) 
\nonumber\\  
& & \hspace{-0.2cm} 
= \int_{{\mathcal{X}}} 
\prod_{i=1}^{d} f_{P_{i}}(x_{i}) \cdot \phi \negthinspace \left( {
\frac{f_{P}(x)}{\prod_{i=1}^{d} f_{P_{i}}(x_{i})}}\right) 
\cdot \boldsymbol{1}_{]0,\infty[}\left(f_{P}(x) \cdot 
\prod_{i=1}^{d} f_{P_{i}}(x_{i})\right)
\, \mathrm{d}\lambda(x) 
\nonumber\\ 
& & \hspace{0.2cm} 
+ \phi^{*}(0) \cdot
P\left[\prod_{i=1}^{d} f_{P_{i}}(x_{i}) =0\right] 
+ \phi(0) \cdot 
Q[f_{P}(x) =0] - \phi(1) \, .
\label{BroStuHB22:fo.def.47ddep}
\end{eqnarray} 
In applications, one often takes $\mathcal{X} =\mathcal{Y} = \mathbb{R}^{d}$,
$\mathcal{Y}_{i}=\mathbb{R}$, $\lambda_{i} := \lambda_{L}$ to be the Lebesgue measure
on $\mathbb{R}$ and thus $\lambda = \lambda_{L}$ is the Lebesgue measure on 
$\mathbb{R}^{d}$ (with a slight abuse of notation),
$f_{P}$ to be the classical joint (Lebesgue) density function of $Y_{1}, \ldots, Y_{d}$,
and  $f_{P_{i}}$ to be the classical (Lebesgue) density function of $Y_{i}$. \\
By plugging $\phi(t) = \phi_{1}(t)= t \cdot \log t + 1 - t \ \in [0, \infty[$ 
$t \in ]0,\infty[$ (cf. \eqref{BroStuHB22:fo.def.120b}) into 
\eqref{BroStuHB22:fo.def.47ddep}, one obtains the prominent mutual information.
References to further subcases of \eqref{BroStuHB22:fo.def.47ddep} can be found e.g. 
in~\cite{Mic:06}.\\
For $d=2$,  $\mathcal{X} =\mathcal{Y} = \mathbb{R}^{2}$, $\lambda := \lambda_{L}$,
 continuous marginal density functions $f_{P_{1}}$, $f_{P_{2}}$,
by Sklar's theorem~\cite{Skl:59} one can uniquely rewrite the joint distribution
function $F_{P}(x_{1},x_{2}) = C(F_{P_{1}}(x_{1}),F_{P_{2}}(x_{1})$
in terms of a copula $C(\cdot,\cdot)$.
Suppose further that $C(\cdot,\cdot)$ is absolutely continuous (with respect to the
Lebesgue measure on $[0,1] \times [0,1]$, and hence for its (Lebesgue) density 
function $c(\cdot,\cdot)$
-- called copula density -- one gets
$c(u_{1},u_{2}) = \frac{\partial^2 C(u_{1},u_{2})}{\partial u_{1} \partial u_{2}}$ 
for almost all $u_{1},u_{2} \in [0,1] \times [0,1]$ (see e.g. p.83 in
Durante \& Sempi~\cite{Dura:16}
and the there-mentioned references).
Accordingly, $f_{P}(x_{1},x_{2}) = f_{P_{1}}(x_{1}) \cdot f_{P_{2}}(x_{2})
\cdot c(F_{P_{1}}(x_{1}),F_{P_{2}}(x_{2}))$
and thus, in case of strictly positive $f_{P_{1}}(\cdot)>0$, $f_{P_{2}}(\cdot)>0$
the divergence \eqref{BroStuHB22:fo.def.47ddep} rewrites as 
\begin{eqnarray} 
& & \hspace{-0.2cm} 
0 \leq D_{\phi,S^{\lambda pd}(Q),S^{\lambda pd}(Q),1\cdot 
S^{\lambda pd}(Q),\lambda}(S^{\lambda pd}(P),S^{\lambda pd}(Q)) 
\nonumber\\ 
& & \hspace{-0.2cm} 
=
\int_{\mathbb{R}} \int_{\mathbb{R}} 
f_{P_{1}}(x_{1}) \cdot f_{P_{2}}(x_{2}) \cdot \phi \negthinspace \left( 
\frac{f_{P}(x_{1},x_{2})}{f_{P_{1}}(x_{1}) \cdot f_{P_{2}}(x_{2})}\right) 
\, \mathrm{d}\lambda_{L}(x_{1}) \, \mathrm{d}\lambda_{L}(x_{2}) 
- \phi(1) \, 
\nonumber\\ 
& & \hspace{-0.2cm} 
= \int_{0}^{1} \int_{0}^{1}  
\phi \negthinspace \left( c(u_{1},u_{2}) \right) 
\, \mathrm{d}\lambda_{L}(u_{1}) \, \mathrm{d}\lambda_{L}(u_{2}) 
- \phi(1) , 
\nonumber
\end{eqnarray} 
which solely depends on the copula (density) and not on the marginals. 
For $\phi(1)=0$ formula (114) was established basically in Durrani \& Zeng~\cite{Dur:09}
without assumptions and without a proof; they also give some 
examples including $\phi = \phi_{\alpha}$ ($\alpha \in \mathbb{R}\backslash\{0,1\}$) 
of \eqref{BroStuHB22:fo.def.32}, as well as the KL-generator $\phi = \tilde{\phi}_{1}(t)$ of
\eqref{BroStuHB22:fo.def.120a} leading to the ``copula-representation of mutual information''.
The latter also appears in the earlier work of Davy \& Doucet~\cite{Dav:03},
as well as e.g. in Zeng \& Durrani~\cite{Zen:11}, Zeng et al.~\cite{Zen:14} and Tran~\cite{Tra:18}; 
in contrast, Tran
also gives a copula-representation of the Kullback-Leibler information divergence
between two general $d-$dimensional Lebesgue density functions
$S_{\cdot}^{\lambda_{L} pd}(P) := f_{P}(\cdot)$ and 
$S_{\cdot}^{\lambda_{L} pd}(Q) := f_{Q}(\cdot)$ 
where $P$ and $Q$ are allowed to have different marginals,
and $Q$ need not be of independence-expressing product type.
\\
(2) For the special case $\mathcal{X}: =\mathcal{Y} = \mathbb{R}^2$, 
continuous marginal distribution functions 
$F_{P_{1}}$ and $F_{P_{2}}$,
product measure $\lambda := P_{1} \otimes P_{2}$,  
$S_{x}^{cd}(Q) := F_{Q}(x) = F_{P_{1}}(x_{1}) \cdot F_{P_{2}}(x_{2}) \in [0,1]$,
as well as joint distribution function $S_{x}^{cd}(P) := F_{P}(x) \in [0,1]$,
one gets the following special cases of \eqref{BroStuHB22:fo.def.140d}
respectively \eqref{BroStuHB22:fo.def.665}:
\begin{eqnarray} 
& & 0 \leq D_{\phi_{2},\mathbb{1},\mathbb{1}, 1\cdot \mathbb{1},\lambda}(S^{cd}(P),S^{cd}(Q)) 
\nonumber\\ 
& & \hspace{-0.2cm} 
= \int_{\mathbb{R}} \int_{\mathbb{R}}  
\frac{1}{2} \cdot
\Big[ F_{P}(x_{1},x_{2}) - F_{P_{1}}(x_{1}) \cdot F_{P_{2}}(x_{2}) \Big]^2
 \, \mathrm{d}P_{1}(x_{1})  \, \mathrm{d}P_{2}(x_{2}) \  
\nonumber
\\ 
& & \hspace{-0.2cm} 
= \int_{0}^{1} \int_{0}^{1}
\Big[ C(u_{1},u_{2}) - u_{1} \cdot u_{2} \Big]^2
\, \mathrm{d}\lambda_{L}(u_{1}) \, \mathrm{d}\lambda_{L}(u_{2}) 
\nonumber
\end{eqnarray}
(cf. Blum et al.~\cite{Blu:61}, Schweizer \& Wolff~\cite{Schw:81}, up to constants and squares) and
\begin{eqnarray} 
& & \hspace{-0.2cm} \textstyle
0 \leq D^{1/2}_{\phi_{TV},S^{cd}(Q),S^{cd}(Q),1\cdot S^{cd}(Q),\lambda}(S^{cd}(P),S^{cd}(Q)) 
= \nonumber \\ 
& & \hspace{-0.2cm} 
=
\int_{\mathbb{R}} \int_{\mathbb{R}}
\left| F_{P}(x_{1},x_{2}) - F_{P_{1}}(x_{1}) \cdot F_{P_{2}}(x_{2}) \right| 
 \, \mathrm{d}P_{1}(x_{1})  \, \mathrm{d}P_{2}(x_{2}) \  
\nonumber
\\ 
& & \hspace{-0.2cm} 
= \int_{0}^{1} \int_{0}^{1}
\left| C(u_{1},u_{2}) - u_{1} \cdot u_{2}  \right|
\, \mathrm{d}\lambda_{L}(u_{1}) \, \mathrm{d}\lambda_{L}(u_{2}) 
\nonumber
\end{eqnarray}
(cf. Schweizer \& Wolff~\cite{Schw:81}, up to constants).

\vspace{0.3cm} 
\noindent
As a side remark, let us mention that other interplays between 
divergences and copula functions can be constructed.
For instance, suppose that $P$ and $Q$ are two probability distributions
on the $d-$dimensional product (measurable) space $(\mathcal{Y},\mathcal{A})$
having copula density functions $c_{P}$ respectively $c_{Q}$; the latter can be 
interpreted as special statistical functionals $S^{cop}(P)$ of $P$ respectively 
$S^{cop}(Q)$ of $Q$, and thus, by employing the divergences
\eqref{BroStuHB22:fo.def.1} we obtain 
\begin{eqnarray} 
& & \hspace{-0.2cm} \textstyle
0 \leq D^{c}_{\phi,m_{1},m_{2},m_{3},\lambda_{{L}^{d}}}(S^{cop}(P),S^{cop}(Q)) 
\nonumber\\ 
& & \hspace{-0.2cm} 
: = \olint_{{\mathcal{Y}}} 
\Bigg[ \phi \negthinspace \left( {
\frac{c_{P}(x)}{m_{1}(x)}}\right) -\phi \negthinspace \left( {\frac{c_{Q}(x)}{m_{2}(x)}}\right)
- \phi_{+,c}^{\prime} \negthinspace
\left( {\frac{c_{Q}(x)}{m_{2}(x)}}\right) \cdot \left( 
\frac{c_{P}(x)}{m_{1}(x)}-\frac{c_{Q}(x)}{m_{2}(x)}\right) 
\Bigg] 
m_{3}(x) \, \mathrm{d}\lambda_{{L}^{d}}(x)
\nonumber \\ 
\label{BroStuHB22:fo.def.divcop}
\end{eqnarray} 
where $\lambda_{{L}^{d}}$ denotes the
$d-$dimensional Lebesgue measure and thus the integral
in \eqref{BroStuHB22:fo.def.divcop} turns out to be (with some rare exceptions)
of $d-$dimensional Riemann-type with $\mathrm{d}\lambda_{{L}^{d}}(x) = \mathrm{d}x$.
The (CASM $\phi-$divergences type) special case 
$D^{c}_{\phi,S^{cop}(Q),S^{cop}(Q),S^{cop}(Q),\lambda_{{L}^{d}}}(S^{cop}(P),S^{cop}(Q))$
leads to a divergence which has been used by 
Bouzebda \& Keziou \cite{Bou:10} in order to obtain
new estimates and tests of independence in semiparametric copula models
with the help of variational methods.

%
\section{Bayesian contexts}
%

There are various different ways how divergences can be used in Bayesian frameworks: 

\vspace{0.2cm}
\noindent
\textit{(1)} \ as ``direct'' quantifiers of dissimilarities between
statistical functionals of various parameter distributions:\\
for instance, consider a $n-$dimensional vector of observable random quantities 
$\mathbf{z}= (Z_{1}, \ldots  Z_{n})$ whose distribution depends on an unobservable 
(and hence, also random) 
multivariate parameter $\mathbf{\Theta}
:= (\Theta_{1}, \ldots, \Theta_{d}) $, 
as well as a real-valued quantity $Z_{n+1}$ 
(whose distribution also depends on $\mathbf{\Theta}$)
to be predicted.  
Corresponding candidates for distributions $P$, $Q$ 
-- to be used in $D(S(P),S(Q))$ -- are for example the following: 
the prior distribution 
$Pr_{\mathbf{\Theta}}[\cdot]  := Pr[\mathbf{\Theta} \in \cdot \, ]$ 
of $\mathbf{\Theta}$
(under some underlying probability measure $Pr$),
the posterior distribution 
$Pr_{\mathbf{\Theta}| \mathbf{z}=\mathbf{z}}[\cdot] := 
Pr[\mathbf{\Theta} \in \cdot \, \Big| \, \mathbf{z}=\mathbf{z} ]$
of $\mathbf{\Theta}$ given the data observation $\mathbf{z}=\mathbf{z}$,
the predictive prior distribution  $Pr_{Z_{n+1}}[\cdot] = Pr[Z_{n+1} \in \cdot \, ]
=  \int_{\mathbb{R}^{d}} Pr_{Z_{n+1}| \mathbf{\Theta}=\mathbf{\theta}}[\cdot]
\, \mathrm{d}Pr_{\mathbf{\Theta}}(\mathbf{\theta})$
of $Z_{n+1}$, and
the predictive posterior distribution  
$Pr_{Z_{n+1}| \mathbf{z}=\mathbf{z}}[\cdot] = Pr[Z_{n+1} \in \cdot \,  \Big| \, \mathbf{z}=\mathbf{z} ]
=  \int_{\mathbb{R}^{d}} Pr_{Z_{n+1}| \mathbf{\Theta}=\mathbf{\theta}}[\cdot]
\, \mathrm{d}Pr_{\mathbf{\Theta}| \mathbf{z}=\mathbf{z}}(\mathbf{\theta})$
of $Z_{n+1}$.
For instance, the divergence 
$D(S^{\lambda pd}(Pr_{\mathbf{\Theta}}),S^{\lambda pd}(Pr_{\mathbf{\Theta}| \mathbf{z}=\mathbf{z}}))$
serves as ``degree of informativity of the new data-point
observation on the learning of the true unknown parameter''.
Analogously, one can also consider more complex setups 
like e.g. a continuum $\mathbf{z}= \{Z_{t} : t \in [0,T]\}$
of observations, parameters $\mathbf{\Theta}$ of function type, 
and $Z_{u}$ ($u >T$) rather than $Z_{n+1}$.

\vspace{0.2cm}
\noindent
\textit{(2)} \ as ``decision risk reduction'' (``model risk reduction'', ``information gain''):
in a dichotomous Bayesian decision problem 
between the two alternative probability distributions $P:=P_{\mathcal{H}}$ and 
$Q:=P_{\mathcal{A}}$, one takes $\mathbf{\Theta} = \{\mathcal{H},\mathcal{A} \}$,
$Pr_{\mathbf{\Theta}}[\cdot] := \pi_{\mathcal{H}} \cdot \delta_{\mathcal{H}}[\cdot]
+ (1-\pi_{\mathcal{H}}) \cdot \delta_{\mathcal{A}}[\cdot]$ for some $\pi_{\mathcal{H}} \in ]0,1[$.
Within this context, suppose we want to make decisions/actions $\mathbb{d}$ taking values in a space $\mathbb{D}$.
Furthermore, for the case that $\mathcal{H}$ were true we attribute a real-valued loss 
$\mathbb{L}_{\mathcal{H}}(\mathbb{d}) \geq 0$ to each $d$; $\mathbb{L}_{\mathcal{H}}(\mathbb{d}) =0$ 
corresponds to a ``right'' decision $d$, $\mathbb{L}_{\mathcal{H}}(\mathbb{d}) > 0$ to the amount of loss
taking the ``wrong'' decision $d$. In the same way, for the case that $\mathcal{A}$ were true
we use $\mathbb{L}_{\mathcal{A}}(\mathbb{d}) \geq 0$.
Prior to random observations $\mathbf{Z}$, the corresponding prior minimal mean decision loss
(prior Bayes loss, prior Bayes risk) is given by
$$
\mathcal{B}(\pi_{\mathcal{H}}) := 
\inf_{\mathbb{d} \in \mathbb{D}} \left\{ \pi_{\mathcal{H}} \cdot \mathbb{L}_{\mathcal{H}}(\mathbb{d})
+ (1-\pi_{\mathcal{H}}) \cdot \mathbb{L}_{\mathcal{A}}(\mathbb{d})
\right\} .
$$
Based upon a concrete observation $\mathbf{z}$, we decide for some ``action'' $\mathbb{d} \in \mathbb{D}$,
operationalized by a decision rule $\mathfrak{d}$ from the space of all possible observations to $\mathbb{D}$
(i.e. $\mathfrak{d}(\mathbf{z}) \in \mathbb{D}$. 
The corresponding posterior minimal mean decision loss
(posterior Bayes loss, posterior Bayes risk) is defined by
\vspace{-0.2cm} 
\begin{eqnarray} 
& & \hspace{-0.2cm} 
\textstyle
\mathcal{B}(\pi_{\mathcal{H}},P_{\mathcal{H}},P_{\mathcal{A}})
:= \inf_{\mathfrak{d}} \left\{ 
\pi_{\mathcal{H}} \cdot 
\int \mathbb{L}_{\mathcal{H}}(\mathfrak{d}(\mathbf{z})) \, \mathrm{d}P_{\mathcal{H}}(\mathbf{z}) 
+ (1-\pi_{\mathcal{H}}) \cdot 
\int \mathbb{L}_{\mathcal{A}}(\mathfrak{d}(\mathbf{z})) \, \mathrm{d}P_{\mathcal{A}}(\mathbf{z}) 
\right\} 
\nonumber
\end{eqnarray}
where the infimum is taken
amongst all ``admissible'' decision functions $\mathfrak{d}$.
Up to technicalities, one can show that
\vspace{-0.2cm} 
\begin{eqnarray} 
& & \hspace{-0.2cm} 
\textstyle
\mathcal{B}(\pi_{\mathcal{H}},P_{\mathcal{H}},P_{\mathcal{A}})
= \int \mathcal{B}(\pi_{\mathcal{H}}^{post}(\mathbf{z}))
\, \left(\pi_{\mathcal{H}} \cdot  \mathrm{d}P_{\mathcal{H}}(\mathbf{z}) + (1- \pi_{\mathcal{H}}) 
\cdot \mathrm{d}P_{\mathcal{A}}(\mathbf{z}) \right),  \qquad 
\nonumber
\end{eqnarray}
with posterior probability (for $\mathcal{H}$) 
$\pi_{\mathcal{H}}^{post}(\mathbf{z}) := 
\frac{\pi_{\mathcal{H}} \cdot f_{P}(\mathbf{z})}{
\pi_{\mathcal{H}} \cdot f_{P}(\mathbf{z}) 
+ (1-\pi_{\mathcal{H}}) \cdot f_{Q}(\mathbf{z})}$
in terms of the $\lambda-$density functions $f_{P}(\cdot)$ and $f_{Q}(\cdot)$
where $\lambda$ is e.g. $\frac{P+Q}{2}$ (or any measure such that $P$ and $Q$ are absolutely 
continuous w.r.t. $\lambda$).
The difference $\mathcal{I}(\pi_{\mathcal{H}},P_{\mathcal{H}},P_{\mathcal{A}})
:= \mathcal{B}(\pi_{\mathcal{H}}) - \mathcal{B}(\pi_{\mathcal{H}},P_{\mathcal{H}},P_{\mathcal{A}}) \geq 0$
can be interpreted as a statistical information measure 
in the sense of De Groot~\cite{DeG:62}, and 
as degree of reduction of the decision risk due to observation.
Let us first discuss the special case $\mathbb{D} = [0,1]$
with $\mathbb{d}$ interpreted as evidence degree, and
$\mathbb{L}_{\mathcal{H}}(\mathbb{d}) = 1 - \mathbb{d}$,
$\mathbb{L}_{\mathcal{A}}(\mathbb{d}) = \mathbb{d}$ 
(Bayes testing). Hence, $\mathcal{B}(\pi_{\mathcal{H}}) = \pi_{\mathcal{H}} \wedge (1-\pi_{\mathcal{H}})$.
From this, {\"O}sterreicher \& Vajda~\cite{Oes:93}, Liese \& Vajda~\cite{Lie:06}
have shown that Csiszar-Ali-Silvey-Morimoto divergences (i.e., $f-$divergences) 
can be represented as ``average''
statistical information measures, i.e. (in our notation) 
\begin{eqnarray} 
& & \hspace{-0.2cm} \textstyle
\int_{]0,1[} 
\mathcal{I}_{\pi_{\mathcal{H}}}(P_{\mathcal{H}},P_{\mathcal{A}}) \frac{1}{\pi_{\mathcal{H}}} 
\, \mathrm{d}g_{\phi}(\pi_{\mathcal{H}}) =
D_{\phi,S_{x}^{\lambda pd}(P_{\mathcal{A}}),S_{x}^{\lambda pd}(P_{\mathcal{A}}),
1\cdot S_{x}^{\lambda pd}(P_{\mathcal{A}}),\lambda}(S_{x}^{\lambda pd}(P_{\mathcal{H}}),S_{x}^{\lambda pd}(P_{\mathcal{A}}))
\nonumber \\
\label{BroStuHB22:fo.def.633new}
\end{eqnarray}
where $g_{\phi}(\pi) := - \phi_{+}^{\prime}\negthinspace\left(\frac{1-\pi}{\pi} \right) $ 
is nondecreasing
in $\pi \in ]0,1[$. If $\phi$ is twice differentiable, then one can simplify
$\frac{1}{\pi_{\mathcal{H}}} \, \mathrm{d}g_{\phi}(\pi_{\mathcal{H}}) = 
\frac{1}{(\pi_{\mathcal{H}})^3} \, 
\phi^{\prime \prime}\negthinspace\left(\frac{1-\pi_{\mathcal{H}}}{\pi_{\mathcal{H}}}\right) \, 
\mathrm{d}\pi_{\mathcal{H}}$
in \eqref{BroStuHB22:fo.def.633new}. For the divergence generators 
$\phi_{\alpha}$ with $\alpha \in \mathbb{R}$ (cf. \eqref{BroStuHB22:fo.def.32}, \eqref{BroStuHB22:fo.def.120b},
\eqref{BroStuHB22:fo.def.120d}, \eqref{BroStuHB22:fo.def.34b})
one gets $\frac{1}{\pi_{\mathcal{H}}} \, \mathrm{d}g_{\phi}(\pi_{\mathcal{H}}) = 
\frac{(1-\pi_{\mathcal{H}})^{\alpha-2}}{(\pi_{\mathcal{H}})^{\alpha+1}} \, \mathrm{d}\pi_{\mathcal{H}}$;
see also Stummer~\cite{Stu:99a,Stu:01a,Stu:04a} for an adaption to a context of 
path-observations of financial diffusion processes.
In contrast, {\"O}sterreicher \& Vajda~\cite{Oes:93}
have also given a ``direct'' representation (in our notation)
\begin{eqnarray} 
& & \hspace{-0.2cm} \textstyle
\mathcal{I}(\pi_{\mathcal{H}},P_{\mathcal{H}},P_{\mathcal{A}})  =
D_{\phi,S_{x}^{\lambda pd}(P_{\mathcal{A}}),S_{x}^{\lambda pd}(P_{\mathcal{A}}),
1\cdot S_{x}^{\lambda pd}(P_{\mathcal{A}}),\lambda}(S_{x}^{\lambda pd}(P_{\mathcal{H}}),S_{x}^{\lambda pd}(P_{\mathcal{A}}))
\qquad \ \ \ 
\nonumber
\end{eqnarray}
for some appropriately chosen loss functions $\mathbb{L}_{\mathcal{H}}(\cdot)$, $\mathbb{L}_{\mathcal{A}}(\cdot)$
which depend on $\phi$ and $\pi_{\mathcal{H}}$
\footnote{
they also have shown some kind of ``reciprocal''
}; 
see also Stummer~\cite{Stu:01a,Stu:04a} for an adaption of the case 
$\phi := \phi_{\alpha}$ with $\alpha \in \mathbb{R}$ within a context of financial diffusion processes.

\vspace{0.2cm}
\noindent
\textit{(3)} \ as bounds of minimal mean decision losses:\\
In the context of (2), let us now discuss the binary decision space
$\mathbb{D} = \{ \mathbb{d}_{\mathcal{H}}, \mathbb{d}_{\mathcal{A}} \}$
where $\mathbb{d}_{\mathcal{H}}$ stands for an action preferred in the 
case that $P_{\mathcal{H}}$ were true. Furthermore, suppose
that $P_{\mathcal{H}}$ is absolutely continuous with respect to
$\lambda:= P_{\mathcal{A}}$ having density function $f_{P_{\mathcal{H}}}(\cdot)$;
notice that $f_{P_{\mathcal{A}}}(\cdot) \equiv 1$.
For the loss functions 
$\mathbb{L}_{\mathcal{H}}(\mathbb{d}) = c_{\mathcal{H}} \cdot \boldsymbol{1}_{\{\mathbb{d}_{\mathcal{A}}\}}(\mathbb{d})$ and
$\mathbb{L}_{\mathcal{A}}(\mathbb{d}) = c_{\mathcal{A}} \cdot \boldsymbol{1}_{\{\mathbb{d}_{\mathcal{H}}\}}(\mathbb{d})$
with some constants $c_{\mathcal{H}}  > 0$, $c_{\mathcal{A}} >0$,
the posterior minimal mean decision loss (posterior Bayes loss) is
\begin{eqnarray} 
& & \hspace{-0.2cm} \textstyle
\mathcal{B}(\pi_{\mathcal{H}},P_{\mathcal{H}},P_{\mathcal{A}})  =
\int \min\{ \Lambda_{\mathcal{H}} \cdot f_{P_{\mathcal{H}}}(\mathbf{z}) , \Lambda_{\mathcal{A}} \} 
\, \mathrm{d}P_{\mathcal{A}}(\mathbf{z})
\qquad \ \ \ 
\nonumber
\end{eqnarray}
with constants $\Lambda_{\mathcal{H}} := \pi_{\mathcal{H}} \cdot c_{\mathcal{H}} >0$,
$\Lambda_{\mathcal{A}} := (1-\pi_{\mathcal{H}}) \cdot c_{\mathcal{A}} >0$. 
For this,
Stummer \& Vajda~\cite{Stu:07b} have achieved the following bounds
in terms of CASD-type power $D:=D_{\phi_{\chi},S_{x}^{\lambda pd}(P_{\mathcal{A}}),S_{x}^{\lambda pd}(P_{\mathcal{A}}),
1\cdot S_{x}^{\lambda pd}(P_{\mathcal{A}}),\lambda}(S_{x}^{\lambda pd}(P_{\mathcal{H}}),S_{x}^{\lambda pd}(P_{\mathcal{A}}))$
for arbitrary $\chi \in ]0,1[$
\begin{equation}
\mathcal{B}(\pi_{\mathcal{H}},P_{\mathcal{H}},P_{\mathcal{A}})
 \ \left\{\begin{array}{ll}
\geq & 	\frac{\Lambda_{\mathcal{H}}^{\max\{1,\frac{\chi}{1-\chi}\}} \cdot
 \Lambda_{\mathcal{A}}^{\max\{1,\frac{1-\chi}{\chi}\}}}{\left(\Lambda_{\mathcal{H}}
+\Lambda_{\mathcal{A}}\right)^{\max\{\frac{\chi}{1-\chi},\frac{1-\chi}{\chi}\}}}\cdot\Big( 
1-\chi\cdot (1-\chi) \cdot D
\Big)^{\max\{\frac{1}{\chi},\frac{1}{1-\chi}\}}\\
 & \\
\leq &  \Lambda_{\mathcal{H}}^{\chi} \cdot \Lambda_{\mathcal{A}}^{1-\chi}\cdot 
(1-\chi\cdot (1-\chi) \cdot D)
\end{array}
\right.\notag
\end{equation}
(in an even slightly more general form), which can be very useful in case that
the posterior minimal mean decision loss can not be computed explicitly. For instance, 
Stummer \& Vajda~\cite{Stu:07b} 
give applications to decision making of time-continuous, non-stationary
financial stochastic processes.

\vspace{0.2cm}
\noindent
\textit{(4)} \ as auxiliary tools:
for instance, in an i.i.d.-type Bayesian parametric model-misspecification context,
Kleijn \& van der Vaart \cite{Kleijn:12}
employ the reverse-Kullback-Leibler-distance minimizer
$$\widehat{\widehat{\theta}} := \arg \inf_{\theta \in \Theta} D_{\phi_{0}}(Q_{\theta}, P_{tr} )
= \arg \inf_{\theta \in \Theta} 
D_{\phi,S^{\lambda pd}(Q),S^{\lambda pd}(Q),
1\cdot S^{\lambda pd}(Q),\lambda}(S^{\lambda pd}(Q_{\theta}),S^{\lambda pd}(P_{tr}))$$
(cf. \eqref{BroStuHB22:fo.def.47d.full} respectively \eqref{BroStuHB22:fo.def.47d} 
with $\phi=\phi_{0}$) in order to formulate and prove an asymptotic normality --- 
under the unknown true out-of-model-lying
data-generating distribution $P_{tr}$
--- of the involved posterior parameter-distribution.

%
\section{Variational Representations}
\label{sec.7new}
%

Variational representations of (say) $\phi-$divergences, often referred to as dual
representation, transform $\phi-$divergence estimation into an optimization
problem on an infinite dimensional function space, generally, but may also
lead to a simpler optimization problem when some knowledge on the class of
measures $Q$ where $D_{\phi }\left( Q,P\right) $ has to be optimized is
available; moreover, as already mentioned at the end of Section \ref{subsec.1c.zeros} above, 
such variational representations can also be employed to circumvent 
the crossover problem (CO1),(C2),(CO3). 

\vspace{0.2cm}
\noindent
To begin with, in the following we loosely sketch the corresponding \textit{general} setting.
We equip $\mathcal{M}$, the linear space of all finite signed measures
(including all probability measures)  on 
$\left( \mathcal{X}\text{,}\mathcal{B}\right) $ with the so called $\tau $
-topology, the coarsest one which makes the mapping $f\rightarrow \int fdQ$
continuous for all measure $Q$ in $\mathcal{M}$ when $f$ runs in the class $
\mathfrak{M}_{b}$ of all bounded measurable functions on $\left( \mathcal{X}
\text{,}\mathcal{B}\right)$. As an exemplary statistical incentive for the use of
signed measures, let us mention the context where one wants to estimate, respectively test for,
a mixture probability distribution $c\cdot Q_{1} + (1-c) \cdot Q_{2}$
with probability measures $Q_{1}$,$Q_{2}$ and $c \in [0,1]$. In such a situation,
it is sometimes technically useful to extend the range of $c$ beyond $[0,1]$
which leads to a signed finite measure. As a next step, 
since the mapping $Q\rightarrow $ $D_{\phi
}\left( Q,P\right) $ is convex and lower semi-continuous in the $\tau $
-topology we deduce that the following result holds for all $Q$ in $\mathcal{
M}$ and $P$ in $\mathcal{P}:$
\begin{equation}
\widetilde{D}_{\phi }\left( Q,P\right) =\sup_{g\in \mathfrak{M}_{b}}\int_{\mathcal{X}} g(x) \, dQ(x)
-\int_{\mathcal{X}} \phi_{\ast }\left( g(x)\right) \, dP(x)  
\label{variational in tau topo}
\end{equation}
where (cf. Broniatowski \cite{Bro:03} in the Kullback-Leibler divergence
case as well as Broniatowski \& Keziou \cite{Bro:06} for a general formulation)
\begin{eqnarray*}
\widetilde{D}_{\phi }\left( Q,P\right)  &:=&
\begin{cases}
\int_{\mathcal{X}} \phi\left( \frac{dQ}{dP}(x)\right) dP(x), 
\qquad \textrm{for } Q<<P, \\
\infty, \hspace{3.5cm} \textrm{else},
\end{cases}
\end{eqnarray*}
is a slightly adopted version of the $\phi-$divergence
defined in \eqref{BroStuHB22:fo.def.47d.full}
(see also \eqref{BroStuHB22:fo.def.47d})
and $\phi_{\ast }(x):=\sup_{t} (t \cdot x - \phi (t))$ designates the Fenchel-Legendre
transform of the generator $\phi$, see \cite{Bro:06} and Nguyen et al. \cite{Ngu:10}. 
The choice of the $\tau $-topology is motivated
by statistical considerations, since most statistical functionals are
continuous in this topology; see Groeneboom et al. \cite{Gro:79}. This choice is in
contrast with similar representations for the Kullback-Leibler divergences
(see e.g. Dembo \& Zeitouni \cite{Dem:09}, under the weak topology on $\mathcal{P}$,
for which the supremum in \eqref{variational in tau topo} is taken over all
continuous bounded functions on $\left( \mathcal{X}\text{,}\mathcal{B}\right)$.

Representation (\ref{variational in tau topo}) offers a useful mathematical
tool to measure statistical similarity between data collections or to
measure the directed distance between a distribution $P$ (either explicit or
known through sampling), and a class of distributions $\Omega $,  as well as
to compare complex probabilistic models. The main practical advantage of
variational formulas is that an explicit form of the probability
distributions or their likelihood ratio, $dQ/dP$, is not necessary. Only
samples from both distributions are required since the difference of
expected values in (\ref{variational in tau topo}) can be approximated by
statistical averages, in case both $Q$ and $P$ are known through sampling.
In practice, the infinite-dimensional function space has to be approximated
or even restricted. One attempt is the restriction of the function space to
a reproducing kernel Hilbert space (RKHS) and the corresponding kernel-based
approximation in Nguyen et al. \cite{Ngu:10}. In many cases of relevance,
however, some information can be inserted in the description of the
minimization problem of the form $\inf \left\{ \widetilde{D}_{\phi }\left( Q,P\right)
;Q\in \Omega \right\} $ when some relation between $P$ and all members in $
\Omega $ can be assumed. Such is the case in logistic models, or more
globally in two sample problems, when it is assumed that $dQ/dP$ belongs to
some class of functions; for example we may assume that $\Omega $ consists
in all distributions such that $x\rightarrow \left( dQ/dP\right) (x)$
belongs to some parametric class. This requires some analysis around 
(\ref{variational in tau topo}), which is handled now.

The supremum in equation (\ref{variational in tau topo}) may not be reached,
even in elementary cases. Consider the case when $\phi =\phi _{1}$,
hence the case when $\widetilde{D}_{\phi }\left( Q,P\right) $ is the Kullback-Leibler
divergence between $Q$ and $P$, and assume that both $Q$ and $P$ are two
Gaussian probability measures on $\mathbb{R}$ with same variance and
different mean values. Then it is readily checked that the supremum in 
(\ref{variational in tau topo}) is reached on a polynomial with degree $2$, hence
outside of $\mathfrak{M}_{b}.$ For statistical purposes it is relevant that
formula \eqref{variational in tau topo} holds with attainment; indeed the
supremum , in case when $\widetilde{D}_{\phi }\left( Q,P\right) $ is finite, is
reached at $g:= \phi ^{\prime }\left( dQ/dP\right) $, therefore, in
case when $\phi $ is differentiable, on a function which may not be
bounded.

It is also of interest to consider \eqref{variational in tau topo} in the
case when $P$ is atomic and $Q$ is a continuous distribution; for example
let $(X_{1},..,X_{n})$ be an i.i.d. sample under \ some probability measure $R$
on $\mathbb{R}$, and consider $Q$ a probability measure absolutely
continuous with respect to the Lebesgue measure; consider the case when 
$\widetilde{D}_{\phi }=\widetilde{D}_{\phi_{1}}$ is the 
(slightly modified) Kullback-Leibler divergence. Denote by $P_{n}^{emp}$ the
empirical measure of the sample. Taking $g(x):= M \cdot \boldsymbol{1}_{\left\{
X_{1},..,X_{n}\right\} ^{c}}(x)$ for some arbitrary $M$, it holds by 
\eqref{variational in tau topo} that $\widetilde{D}_{\phi }(Q,P_{n}^{emp})\geq M$ proving that
no inference can be performed about $R$ making use of the variational form
as it stands.\ Some more structure and information has to be incorporated in
the variational form of the divergence in order to circumvent this obstacle.
Assuming that  $\phi $ is a differentiable function in its domain, the
supremum in (\ref{variational in tau topo}) is reached at 
$g_{\ast}:=\phi ^{\prime }\left( dQ/dP\right) $ as checked by substitution
\footnote{In case when $\phi$ is not differentiable
at some point, then the supremum in \eqref{Variational form with class F}
should satisfy $g_{\ast }(x)\in \partial \phi \left( dQ/dP\right) (x)$
for all $x$ in $\mathcal{X}$, where $\partial \phi (t)$
is the subdifferential set of the convex function $\phi $
at point $t$, $\partial \phi (t):=\left\{ z\in 
\mathbb{R}:\phi (s)\geq \phi (t)+z\left( s-t\right) ,\forall s\in 
\mathbb{R}\right\}$
}.
Let $\mathcal{F}$ be a class of
functions containing all functions $\phi ^{\prime }\left( dQ/dP\right) (x)
$ as $Q$ runs in a given model $\Omega $.Consider the subspace 
$\mathcal{M}_{\mathcal{F}}$ of all finite signed measures $Q$ such that $\int \left\vert
f\right\vert d\left\vert Q\right\vert $ is finite for all function $f$ in $
\mathcal{F}$, then  similarly as in (\ref{variational in tau topo}) we may
obtained the following variational form of $\widetilde{D}_{\phi }\left( Q,P\right) $,
which is valid when $Q$ belongs to $\mathcal{M}_{\mathcal{F}}$ and $P$
belongs to $\mathcal{P}$
\begin{equation}
\widetilde{D}_{\phi }\left( Q,P\right) =\sup_{g\in \left\langle \mathfrak{M}_{b}\cup 
\mathcal{F}\right\rangle }\int_{\mathcal{X}} g(x) \, dQ(x)
-\int_{\mathcal{X}} \phi_{\ast }\left( g(x)\right) \, dP(x)
\label{Variational form with class F}
\end{equation}
in which we substituted $\mathfrak{M}_{b}$ by the broader class $\mathfrak{M}
_{b}\cup \mathcal{F}$ which may contain unbounded functions; note that 
(\ref{Variational form with class F}) is valid for a smaller class of measures $Q$
than \eqref{variational in tau topo}.  

\noindent
For instance, in the above example pertaining to the
Kullback-Leibler divergence and both $P$ is Gaussian on $\mathbb{R}$ and  $Q$
belongs to the class $\Omega $ of all Gaussian distributions on $\mathbb{R}$
with same variance as $P$,  then $\mathcal{F}$ consists of all polynomial
functions with degree $2$, and the supremum in (\ref{Variational form with
class F}) is attained.\ Looking at the case when $P$ is substituted by $P_{n}
$ and $Q$ is absolutely continuous , and since  $\widetilde{D}_{\phi }\left(
Q,P_{n}\right) $ does not convey any information from the data, we are led
to define a restriction to the supremum operation on the space $\left\langle 
\mathfrak{M}_{b}\cup \mathcal{F}\right\rangle $; since we assumed that  
$\phi ^{\prime }\left( dQ/dP\right) \in \mathcal{F}$ for any $Q$ in 
$\Omega \subset \mathcal{M}_{\mathcal{F}}$ we have   
\begin{equation}
\widetilde{D}_{\phi }\left( Q,P\right) =\sup_{g\in \mathcal{F}}
\int_{\mathcal{X}} g(x) \, dQ(x) - \int_{\mathcal{X}} 
\phi_{\ast }\left( g(x)\right) \, dP(x)  
\label{Var form reduced}
\end{equation}
which is valid only when $Q<<P$. We thus can define a new ``pseudo
divergence'', say \underline{$\widetilde{D}$}$_{\phi }\left( Q,P\right) $ which
coincides with $\widetilde{D}_{\phi }\left( Q,P\right) $ in those cases, and which
takes finite values depending on the data when $P$ is substituted by $P_{n}^{emp}$.
In that case we define 
\begin{equation}
\underline{\widetilde{D}}_{\phi }\left( Q,P_{n}^{emp}\right) := \sup_{g\in \mathcal{F}}
\int_{\mathcal{X}} g(x) \, dQ(x) - \frac{1}{n}\sum_{i=1}^{n}\phi_{\ast }\left( X_{i}\right) ,
\label{variation form empirical}
\end{equation}
which is the starting point of variational divergence-based inference; see 
Broniatowski \& Keziou \cite{Bro:09}. Note that the above formula does not require any grouping
or smoothing.
Also the resulting estimator of the likelihood ratio 
$dQ^{\ast}/dP$ where $Q^{\ast }:=\arg \inf_{Q\in \Omega }\widetilde{D}_{\phi }\left(
Q,P\right) $ results from a double optimization, the inner one pertaining to
the estimation of $g_{\ast }(Q)$ solving (\ref{variation form empirical})
for any $Q$ in $\Omega $. \ Assuming that $\left\{ \phi_{\ast }(g),g\in 
\mathcal{F}\right\} $ is a Glivenko-Cantelli class of functions in some
appropriate metrics provides the ingredients to handle convergence
properties of the estimators. The choice of the divergence $\phi $ may
obey robustness vs efficiency equilibrium, as exemplified in parametric
models; see also Al Mohamad \cite{AlM:18}. 

Formula \eqref{Var form reduced} can be obtained through simple convexity
considerations (see p. 172 of Liese \& Vajda \cite{Lie:87} or 
Theorem 17 of Liese \& Vajda \cite{Lie:06}) 
and is used when $\mathcal{F}$ consists in all the functions 
$\phi ^{\prime }\left( dQ/dP\right) $ as $Q$ and $P$ run in some
parametric model. In a more general (semiparametric or nonparametric
setting), formula (\ref{Variational form with class F}) is adequate for
inference in models consisting in probability distributions $Q$ 
which integrate functions in $\mathcal{F}$, 
and leads to numerical optimization making use of regularity
assumptions on the likelihood ratio $dQ/dP.$

%
\section{Some Further Variants}
%

Extending the (say) $\phi-$divergence definition 
outside the natural context of probability measures
appear as necessary in various situations; for example models defined by
conditions pertaining to expectations of order statistics (or more generally
of L-statistics) are ubiquitous in meteorology, hydrology or in finance
through constraints on the value at risk, for example on the Distortion Risk
Measure (DRM) of index $\alpha $, which is defined in terms of the quantile
function $F^{\leftarrow}$ associated to the distribution function $F$ on 
$\mathbb{R}^{+}$ through $\int_{0}^{1}F^{\leftarrow}(u) \cdot 
\boldsymbol{1}_{\{u>\alpha \}} \, du.$ 
Note that
this class of constraints are not linear with respect to $F$ but with
respect to $F^{\leftarrow}$ only; so characterization of the projection of some
measure $P$ on such sets of measures are not characterized by exponential
family types of distributions.\ Inference on whether a distribution $P$
satisfies this kind of constraints leads to the extension of the definition
of divergences between quantile measures, which might be signed measures.
Variational representations for inference can be defined and projections on
linear constraints pertaining to quantile measures can be characterized; see 
Broniatowski \& Decurninge \cite{Bro:16}. Also in the statistical frame, testing for the number of
components in a finite mixture requires the extension of the definition of
divergences to not-necessarily positive argument, such as occurs for the
Pearson $\chi^{2}-$divergence; this allows to replace the non-regular statistical 
task of estimating (testing) a value of a parameter at the border of its domain into
a regular problem, at the cost of introducing mixtures with negative
weights; an attempt in this direction is made in Broniatowski et al. \cite{Bro:19b}.

\vspace{0.3cm}
\noindent
For large-dimensional spaces $\mathcal{X}$, variational representations of
$\phi-$divergences (i.e. CASM divergences) offer significant theoretical insights and practical advantages
in numerous research areas. Recently, they have gained popularity in machine
learning as a tractable and scalable approach for training probabilistic
models and for statistically differentiating between data distributions; see e.g.
Birrell et al. \cite{Birr:21b}. 

Explicit methods to estimate the $\phi-$divergence and likelihood ratio between two
probability measures known through sampling (hence substituting $Q$ and $P$
in (\ref{Var form reduced}) by their empirical counterparts) have been
considered making some hypothesis on its regularity, or adding some penalty
term in terms of the assumed complexity of the class $\mathcal{F};$ examples
include Sobolev classes of functions or Reproducing Kernel Hilbert Space
approximations; see Nguyen et al. \cite{Ngu:10} for explicit methods and
properties of the estimators.

\vspace{0.3cm}
\noindent
Extensions of the basic divergence formula as given in (\ref{Var form
reduced}) to include some extra inner optimization term have been proposed
in by Birrell et al. \cite{Birr:22} under the name of $(f-\Gamma)-$divergences; 
this new class encompasses both the $\phi-$divergence 
class and many integral probability metrics
(see also Sriperumbudur et al. \cite{Sri:12} on the overlap of the latter two); they provide
uncertainty quantification bounds for misspecified models in terms of the
$\phi-$divergence between the truth and the model, somehow in a similar way as
considered in cryptology (see Arikan \& Merhav \cite{Arik:98} and subsequent extensive
literature). Also, \cite{Birr:22} apply optimization of those divergences to training
Generative Adversial Networks (GAN).

\vspace{0.2cm}
\noindent
Another area where extension of the $\phi-$divergences (i.e. CASM divergences) 
to signed measures is
useful is related to general optimization problems, where one aims at
projecting a vector (or a function) on a class of vectors (or a class of
functions); we refer to Broniatowski \& Stummer \cite{Bro:21a} for an extensive
treatment of such problems in the finite dimensional case.

\vspace{0.2cm}
\noindent
As already indicated above, there are also divergences between stochastic processes 
where $\mathcal{X}$
is the set of all possible paths (i.e. all time-evolution scenarios).
By nature, the analysis of the outcoming (say) $\phi-$divergences
between two distributions on the path space  $\mathcal{X}$
may become very involved. For instance, power divergences between  
diffusion processes  --- and applications to finance, Bayesian decision making, etc. ---
were treated in Stummer~\cite{Stu:99a,Stu:01a,Stu:04a} as well as in
Stummer \& Vajda \cite{Stu:07b} (see also the corresponding binomial-process-approximations in 
Stummer \& Lao \cite{Stu:12b});
in contrast, Kammerer \& Stummer \cite{Kam:20} study power divergences between
Galton-Watson branching processes with immigration and apply the outcomes to
optimal decision making in the presence of a pandemics (such as e.g. COVID-19).

\vspace{0.2cm}
\noindent
For continuous, convex, homogeneous functions $\phi:\mathbb{R}_{+}^{K} \mapsto \mathbb{R}$,
general multivariate $\phi-$dissimilarities of the form
\begin{equation}
\mathcal{D}_{\phi} \left( \mathbf{Q},P\right) =\int_{\mathcal{X}} \phi \left( \frac{dQ_{1}}{dP}(x)
,\ldots,\frac{dQ_{K}}{dP}(x)\right) \, dP(x)
\nonumber
\end{equation}
(which need not necessarily be divergences in the sense of a multivariate
analogue of the above axioms $(D1),(D2)$) have been first introduced by
Gy\"{o}rfi \& Nemetz \cite{Gyo:77} \cite{Gyo:78} 
and later on investigated by e.g. Zografos \cite{Zog:94} for stratified random sampling,
by Zografos \cite{Zog:98} for hypothesis testing, and by Garcia-Garcia \& Williamson
\cite{Gar:12} for multiclass classification problems.
As noticed by \cite{Gyo:77}, the multivariate $\phi-$dissimilarities cover as
special cases Matusita's affinity \cite{Mat:67}, the more general Toussaint's affinity
\cite{Tou:74} \cite{Tou:78}
(which by nature is a multivariate (form of a) Hellinger integral 
being also called Hellinger transform in Liese \& Miescke \cite{Lie:08}).,
and --- in the bivariate case $K=2$ --- also the $\phi-$divergences(i.e the CASM divergences).
Special multivariate $\phi-$divergences $\mathcal{D}_{\phi} \left( \mathbf{Q},P\right)$ were e.g. 
employed by Toussaint \cite{Tou:74} \cite{Tou:78} (see also Menendez et al. \cite{Men:97})
in form of an average over all pairwise
Jeffreys divergences (where the latter are sum-symmetrized Kullback-Leibler divergences),
by Menendez et al. \cite{Men:92} in form of a convex-combination of ``Kullback-Leibler divergences
between each individual probability distribution and the convex-combination of all probability
distributions'' (i.e. multivariate extensions of the Jensen-Shannon divergence), 
and by Werner \& Ye \cite{Wer:17} in form of integrals over the geometric mean of
all the integrands in pairwise $\phi$-divergences (and they even flexibilize 
to components of a $\mathbb{R}_{+}^{K}-$valued function $\mathbf{\phi}$,
and call the outcome a mixed $\mathbf{\phi}-$divergence). 
A \textit{general} ``natural multivariate'' extension of a $\phi-$divergence in the
sense of CASM --- called multidistribution $\phi-$divergence --- 
has been given by Duchi et al. \cite{Duc:18} who employed this
to multiclass classification problems (see also Tan \& Zhang \cite{Tan:21}
for further application to loss functions and regret bounds).
The general multivariate $\phi-$dissimilarity between signed measures
(rather than the more restrictive probability distributions) --- under
assumptions which imply the multivariate analogue of the above axioms $(D1),(D2)$ ---
has been introduced by Keziou \cite{Kez:15}
and used for the analysis of semiparametric 
multisample density ratio models (for the latter, 
see e.g. Keziou \& Leoni-Aubin \cite{Kez:08} and 
Kanamori et al. \cite{Kan:12}).

%
%

\vspace{0.3cm}
\noindent 
\textbf{Acknowledgement.} \, 
W. Stummer is grateful to the  
Sorbonne Universit\'{e} 
Paris for its multiple partial financial support and especially the LPSM 
for its multiple great hospitality.
M. Broniatowski thanks very much the University of Erlangen-N{\"u}rnberg
for its partial financial support and hospitality.
Moreover, W. Stummer would like to thank Ingo Klein and Konstantinos Zografos
for some helpful remarks on a much earlier draft of this paper.

%
%

\end{document}